\documentclass[a4paper,10pt,leqno]{article}
\usepackage[english]{babel}
\usepackage[utf8]{inputenc}
\usepackage{amssymb}
\usepackage{amsmath}
\usepackage{amsthm}
\usepackage{hyperref}
\usepackage{graphicx}
\renewcommand{\theequation}{%
  \ifnum\value{subsection}=0
    \thesection.\arabic{equation}%
  \else
    \ifnum\value{subsubsection}=0
      \thesubsection.\arabic{equation}%
    \else
      \thesubsubsection.\arabic{equation}%
    \fi
  \fi
}
\renewcommand\theequation{\thesubsection.\arabic{equation}}
\makeatletter
\@addtoreset{equation}{subsection}
\makeatother
\input xy
\xyoption{all}

\newtheorem{proposition}[subsection]{Proposition}

\newtheorem{theorem}[subsection]{Theorem}
\newtheorem{theorem2}[subsubsection]{Theorem}

\newtheorem{corollary}[subsection]{Corollary}

\theoremstyle{definition}
\newtheorem{definition}[subsection]{Definition}

\newtheorem{remark}[subsection]{Remark}

\newtheorem{example}[subsection]{Example}

\newcommand{\SP}{\mathrm{Spec}}         

\newcommand{\ori}{{\stackrel{\la}{\times}}} 

\def\commutatif{\ar@{}[rd]|{\circlearrowleft}}
\def\cartesien{\ar@{}[rd]|{\square}}

\newcommand{\ra}{\rightarrow}
\newcommand{\la}{\leftarrow}
\newcommand{\iso}{\stackrel{\sim}{\ra}} 
               
\newcommand{\surj}{\twoheadrightarrow}  
\newcommand{\inj}{\hookrightarrow}      

\def\fl#1{\overleftarrow{#1}} 

\renewcommand{\lim}{{\mathrm{lim}}} 

\def\egazéro#1#2{[{\bf ÉGA}~$0_{\textsc{#1}}$~#2]}         

\begin{document}

\medskip
\centerline{\textbf{Around the Thom-Sebastiani theorem}}

\medskip
\centerline{\textbf{Luc Illusie}\footnote{Laboratoire de Mathématiques d'Orsay, Univ. Paris-Sud, CNRS, Université
Paris-Saclay, 91405 Orsay, France; email: \texttt{Luc.Illusie@math.u-psud.fr}.}}

\medskip
\centerline{\textbf{with an appendix by Weizhe Zheng}\footnote{Morningside Center of Mathematics, Academy of
Mathematics and Systems Science, Chinese Academy of Sciences, Beijing
100190, China; University of the Chinese Academy of Sciences, Beijing 100049, China; email: \texttt{wzheng@math.ac.cn}. Partially supported by
National Natural Science Foundation of China Grant 11321101.}}

\medskip
\textbf{Abstract.} {\small For germs of holomorphic functions $f : (\mathbf{C}^{m+1},0) \to
(\mathbf{C},0)$, $g : (\mathbf{C}^{n+1},0) \to (\mathbf{C},0)$ having an isolated
critical point at 0 with value 0, the classical Thom-Sebastiani
theorem describes the vanishing cycles group $\Phi^{m+n+1}(f \oplus
g)$ (and its monodromy) as a tensor product $\Phi^m(f) \otimes
\Phi^n(g)$, where $(f \oplus g)(x,y) = f(x) + g(y), x = (x_0,...,x_m),
y = (y_0,...,y_n)$. We prove algebraic variants and
generalizations of this result in étale cohomology over fields of any characteristic,
where the tensor product is replaced by a certain local convolution
product, as suggested by Deligne. They generalize \cite{Fu13}. The main ingredient is a Künneth
formula for $R\Psi$ in the framework of Deligne's theory of nearby
cycles over general bases. In the last section, we study the tame case, and the relations between tensor and convolution products, in both global and local situations.}



\medskip
\centerline{\textbf{Contents}}

0. Introduction

1. Review of nearby cycles over general bases

2. External tensor products

3. Interlude: additive convolution

4. Application to Thom-Sebastiani type theorems

5. The tame case

Appendix A. Künneth formula for nearby cycles

\medskip
\noindent \textbf{0. Introduction}

\medskip
If $f : (\mathbf{C}^{m+1},0) \to (\mathbf{C},0)$, $g : (\mathbf{C}^{n+1},0) \to (\mathbf{C},0)$ are germs of holomorphic functions having $0$ as an isolated critical point with value $0$, the germ $f \oplus g : (\mathbf{C}^{m+n+2},0) \to (\mathbf{C},0)$ defined by $(f \oplus g)(x,y) = f(x) +g(y)$ has also $0$ as an isolated critical point, and the classical Thom-Sebastiani theorem \cite{ThomSebastiani71} expresses its  group of vanishing cycles at $0$ as a tensor product:
$$
\Phi^{m}(f) \otimes \Phi^n(g) \iso \Phi^{m+n+1}(f \oplus g). \leqno (0.1)
$$
Here, if $h : (\mathbf{C}^r,0) \to (\mathbf{C},0)$ is a germ of holomorphic function having an isolated critical point at $0$, $\Phi^q(h):= R^q\Phi_h(\mathbf{Z})_0$  is the stalk at $0 \in \mathbf{C}^r$ of $H^qR\Phi_h(\mathbf{Z})$, where $R\Phi_h$ is the vanishing cycles functor of (\cite{SGA7}, XIV); this is also $\widetilde{H}^q(M_h,\mathbf{Z})$, where $M_h$ is a Milnor fiber of $h$ at $0$, and $\widetilde{H}^q = \mathrm{Coker} \,H^q(\mathrm{pt}) \to H^q$. The isomorphism (0.1) is compatible with the monodromy operators $T_f$, $T_g$, $T_{f \oplus g}$, i. e.
$$
T_f \otimes T_g = T_{f \oplus g}\leqno (0.2)
$$
via (0.1). In (\cite{Deligne72}, \cite{Deligne99}), Deligne gave a refinement of (0.2), with the monodromy operators replaced by the \textit{variation} isomorphisms $V_f : R^m\Phi_f(\mathbf{Z})_0 \iso H^m_{\{0\}}(R\Psi_f(\mathbf{Z}))$ (and similarly for $g$ and $f\oplus g$), in the notation of (\cite{SGA7}, XIV), where the group $H^m_{\{0\}}(R\Psi_f(\mathbf{Z}))$, which is dual to $R^m\Phi_f(\mathbf{Z})_0$, is also isomorphic to $H^m_c(M_f - \partial M_f,\mathbf{Z})$, $\partial M_f$ denoting the boundary of $M_f$, namely, 
$$
V_f \otimes V_g = V_{f \oplus g} \leqno (0.3)
$$ 
via (0.1) and its dual. The above groups of vanishing cycles are of an algebraic nature, as it is known that they depend only on a suitable high order jet of the functions. However, the proofs of (0.2) and (0.3) are transcendental. They heavily rely on a description of the Milnor fiber $M_{f \oplus g}$ as homotopic to a join $M_f * M_g$. 

It had been observed by Deligne long ago that, in positive characteristic, an $\ell$-adic analogue of (0.1), compatible with Galois actions, could not hold, as could already be seen in the Picard-Lefschetz situation for quadratic singularities. He suggested that the tensor product on the left hand side should be replaced by a certain local convolution product. More precisely, consider the following set-up. Let $k$ be an algebraically closed field of characteristic $p$. For $i = 1, 2$, let $f_i : X_i \to \mathbf{A}^1_k = \SP \, k[t]$ a flat morphism of finite type with $X_i$ smooth over $k$, of dimension $n_i +1$. Let $n = n_1 + n_2$. Assume that $f_i$ has an isolated critical point at a rational point $x_i$ of the special fibre. Let $a : \mathbf{A}^1_k \times_k \mathbf{A}^1_k \to \mathbf{A}^1_k$  denote the sum map $(u,v) \mapsto u + v$.  Let $f = f_1 \times_k f_2 : X_1 \times_k X_2 \to \mathbf{A}^1_k$. Then the composite morphism
$$
af : X_1 \times_k X_2 \to \mathbf{A}^1_k
$$
(i. e. $f_1 \oplus f_2$) has again an isolated critical point at the point $x = (x_1,x_2)$ of its special fiber. Let $A = \SP \,k\{t\}$ be the henselization at the origin of $\mathbf{A}^1_k$, let $\overline{\eta}$ be a geometric point over the generic point $\eta$ of $A$. Let $\ell$ be a prime number different from $p$. It is known that in this case the vanishing cycles group $R^q\Phi_{f_i}(\mathbf{Z}_{\ell})_{x_i}$ (resp. $R^q\Phi_{af}(\mathbf{Z}_{\ell})_x$) is zero for $q \ne n_i$ (resp. $q \ne n+1$) and $R^{n_i}\Phi_{f_i}(\mathbf{Z}_{\ell})_{x_i}$ (resp. $R^{n+1}\Phi_{sf}(\mathbf{Z}_{\ell})_x$) is a free $\mathbf{Z}_{\ell}$-module of finite type (cf. (\cite{SGA7}, I 4.6) for $p = 0$, (\cite{Illusie03}, 2.10) for the general case). Consider $R^{n_i}\Phi_{f_i}(\mathbf{Z}_{\ell})_{x_i}$ as a sheaf on $\eta$, extended by zero on $A$, and the external tensor product on the henselization $A^2_{(0,0)}$ of $A \times_k A$ (or $\mathbf{A}^2_k$) at $(0,0)$
$$
M := R^{n_1}\Phi_{f_1}(\mathbf{Z}_{\ell})_{x_1} \boxtimes R^{n_2}\Phi_{f_2}(\mathbf{Z}_{\ell})_{x_2} = \mathrm{pr}_1^*R^{n_1}\Phi_{f_1}(\mathbf{Z}_{\ell})_{x_1} \otimes \mathrm{pr}_2^*R^{n_2}\Phi_{f_2}(\mathbf{Z}_{\ell})_{x_2}
$$
Denote again by $a : A^2_{(0,0)} \to A$ the map induced by the sum map. We have $R^q\Phi_a(M)_{(0,0)} = 0$ for $q \ne 1$, and in his seminar \cite{Deligne80}, Deligne sketched a construction of an isomorphism
$$
R^1\Phi_a(M)_{(0,0)} \iso R^{n+1}\Phi_{af}(\mathbf{Z}_{\ell})_{x} , \leqno (0.4)
$$
such that when $k$ is an algebraic closure of a field $k_0$, (0.4) is compatible with the action of $\mathrm{Gal}(k/k_0)$.   This is an analogue of (0.1), where the vanishing cycles group on the left hand side replaces the tensor product. His construction used a suitable compactification of $a$ and $f$. It has not been written up.  The functor associating to a pair $(V_1,V_2)$ of $\overline{\mathbf{Q}}_{\ell}$-sheaves on $\eta$ (i.e., $\mathrm{Gal}(\overline{\eta}/\eta)$-modules), considered as sheaves on $A$ extended by zero, the sheaf
$$
R^1\Phi_a(V_1 \boxtimes V_2)_{(0,0)} \leqno (0.5)
$$
(denoted $V_1 *_1 V_2$ (\ref{*1}) in our paper) was then extensively studied by Laumon \cite{Laumon87} under the name of \textit{local (additive) convolution}. 
Fu \cite{Fu13} revisited the question, and gave a proof of (0.4) (with $\mathbf{Z}_{\ell}$ replaced by $\overline{\mathbf{Q}}_{\ell}$) (and a slight generalization of it), using Laumon's local Fourier transform \cite{Laumon87}. In \cite{Deligne11} Deligne conjectured generalizations taking for coefficients objects of $D^b_c$ and replacing the sum map $a$ by other types of morphisms (see \cite{Fu13} for a precise statement), generalizations which seemed out of reach of the method of \cite{Fu13}.

In this paper, we address Deligne's expectations. Here is the idea. By a basic result of Gabber (a key ingredient in \cite{Fu13}), $R\Psi$ commutes with external tensor products (\cite{Illusie94}, 4.7). However, the external products in question in \textit{loc. cit.} are products $Y_1 \times_S Y_2$, for $S$ a trait, and $S$-schemes $Y_1$, $Y_2$, while here we need external products of the form $f_1 \times_k f_2$, with target of dimension 2. In order to deal with such morphisms, we use Deligne's theory of nearby cycles over general bases, developed in \cite{Illusie12} and \cite{Orgogozo06}. 

For a morphism of finite type $f : X \to Y$ between noetherian schemes, and $K \in D^b_c(X,\Lambda)$, $\Lambda = \mathbf{Z}/n\mathbf{Z}$, $n \ge 1$ invertible on $Y$, the nearby cycles complex $R\Psi_f(K)$ is an object of $D^b(X \ori_Y Y,\Lambda)$, where $X \ori_Y Y$ is the \textit{vanishing topos} of $f$, a topos, not a scheme, which, in the case where $Y$ is a strictly local trait, is a slight enrichment of the topos of sheaves on the special fibre endowed with an action of the inertia group of $Y$. While, if $Y$ is a henselian trait, the formation of $R\Psi_f(K)$ commutes with (surjective) base change of traits (\cite{SGA41/2}, Th. finitude, 3.7), for $\mathrm{dim} Y > 1$, the formation of $R\Psi_fK$ in general does not commute with base change. We say that $(f,K)$ is \textit{$\Psi$-good} when the formation of $R\Psi_fK$ commutes with arbitrary base change $Y' \to Y$ (see \ref{psi-good} for a precise formulation). In this case, we have $R\Psi_fK \in D^b_c(X \ori_Y Y,\Lambda)$, i. e., for all $q$, $R^q\Psi_fK$ is a constructible sheaf of $\Lambda$-modules on $X \ori_Y Y$ (\ref{Orgogozo's theorem}), which vanishes for $q$ sufficiently large. Now, a key result of this paper is the following (Theorem \ref{kunneth theorem}):

\textbf{Theorem}. \textit{Let $S$ be a noetherian scheme, and for $i = 1, 2$, let $f_i : X_i \to Y_i$ be a morphism of $S$-schemes of finite type. Let $X := X_1 \times_S X_2$, $Y := Y_1 \times_S Y_2$, and $ f := f_1 \times_S f_2 : X \to Y$. Let $K_i \in D_{ctf}(X_i,\Lambda)$, with $\Lambda$ as above (or more generally, a noetherian $\mathbf{Z}/n\mathbf{Z}$-algebra), with $n$ invertible on $S$ (\textit{ctf} meaning that $K_i$ is in $D^b_c$ and its stalks are of finite tor-dimension). Assume that $(f_i,K_i)$ is $\Psi$-good. Then the external product $(f,K = K_1 \boxtimes^L K_2)$ is again $\Psi$-good, and the natural map
$$
R\Psi_{f_1}K_1 \boxtimes^L R\Psi_{f_2}K_2 \to R\Psi_f(K_1 \boxtimes^L K_2) \leqno (0.6)
$$
is an isomorphism.}

The proof given in this paper uses the following ingredients: (a) Gabber's theorem quoted above, on the commutation of $R\Psi$ with external products over a trait, (b) Orgogozo's theorem (\cite{Orgogozo06}, 2.1) to the effect that there exists a modification $g : Y' \to Y$ such that base change by $g$ renders $(f,K)$ $\Psi$-good, (c)  Gabber's theorem of oriented cohomological descent (\cite{Orgogozo12b}, 2.2.3).  A new, simpler proof, under slightly relaxed hypotheses, has just been obtained by W. Zheng, see A.3. It needs none of the above ingredients,  and relies only on general nonsense on oriented products and projection formulas.  It is similar in spirit to that of the Künneth formula (\cite{SGA5}, III (1.6.4)).

Isomorphisms of type (0.4) are then more or less formal applications of \ref{kunneth theorem} and general transitivity properties of $R\Psi$.

In §1 we recall basic definitions and facts about nearby cycles on general bases. §2 is devoted to the statement and proof of \ref{kunneth theorem}. Applications of Thom-Sebastiani type are discussed in §4, after reviewing standard material on global and local additive convolutions in §3. Formula (0.4) is obtained in \ref{isolated2}. Generalizations conjectured by Deligne in \cite{Deligne11} are discussed in \ref{Deligne's variant}. In § 5 we study the convolution product $V_1 *_1 V_2$ (\ref{*1}) for $V_1$ and $V_2$ tamely ramified, and analyze its ``arithmetic difference" with $V_1 \otimes V_2$. The results are due to Deligne (\cite{Deligne80} and private communication). Their presentation owes much to discussions with Laumon. In the situation considered at the beginning of this introduction, we give formulas for monodromy and variation in the tame case (\ref{TS - monodromy *}), (\ref{variation5}), recovering (0.1), (0.2) and (0.3).

\medskip
\textit{Conventions}: For a morphism $f : X \to Y$ and a sheaf (or complex) $K$ on $Y$ we will write $K|X$ for $f^*K$ when no confusion can arise. We will sometimes, by abuse, say that a diagram is commutative when it is commutative up to a canonical isomorphism. Rings are assumed to be commutative and unital.

\medskip
\section{Review of nearby cycles over general bases}
\subsection{}\label{oriented product}
Recall the following construction, due to Deligne (see \cite{Illusie12}). Given morphisms of topoi $f : X \to S$, $g : Y \to S$, there is defined a topos $X \ori_S Y$, called the \textit{(left) oriented product} of $X$ and $Y$ over $S$, together with maps $p_1 : X \ori_S Y \to X$, $p_2 : X \ori_S Y \to Y$, and a 2-map $\tau : gp_2 \to fp_1$, which is universal for these data. If $X$, $Y$, $S$ are the categories of sheaves on small sites $C_1$, $C_2$, $D$ satisfying standard exactness properties, $X \ori_S Y$ is the category of sheaves on a site $C$, consisting of pairs of maps $(U \to V \leftarrow W)$ above $(X \overset{f}\to S \overset{g}\leftarrow Y)$, i.e., $U \to f^*V$, $W \to g^*V$, with a topology generated by covering families $(U_i \to V \leftarrow W)$ for $(U_i \to U)$ covering in $C_1$, $(U \to V \leftarrow W_i)$ for $(W_i \to W)$ covering in $C_2$, and
\begin{equation}\label{type c}
\xymatrix{{} & V' \ar[d] & W' \ar[d] \ar[l]\\
U \ar[r] \ar[ur] & V & W \ar[l]},
\end{equation}
where the square is cartesian. If $\mathcal{F}$ is a sheaf on $C$, i.e., an object of $X \ori_S Y$, the restriction map $\mathcal{F}(U\to V \leftarrow W) \to \mathcal{F}(U \to V' \leftarrow W')$ is an isomorphism for any diagram of type (\ref{type c}). If $e_X$, $e_S$, $e_Y$ denotes the final objects of $X$, $S$, $Y$, the maps $p_1$, $p_2$ are given by $p_1^*(U) = (U \to e_S \leftarrow e_Y)$, $p_2^*(W) = (e_X \to e_S \leftarrow W)$.

By the universal property of oriented topoi, the datum of a diagram of topoi
\begin{equation}\label{functoriality1}
\xymatrix{X' \ar[d]_u \ar[r]^{f'} & S' \ar[d]_v & Y' \ar[l]_{g'} \ar[d]_w \\
X \ar[r]^f &S & Y \ar[l]_g}
\end{equation}
and 2-maps $a : gw \to vg'$, $b : vf' \to fu$ defines a morphism of topoi
\begin{equation}\label{functoriality2}
(u,v,w;a,b) : X' \ori_{S'} Y' \to X \ori_S Y,
\end{equation}
called the \textit{functoriality morphism}, denoted sometimes $\fl{u \times_v w}$ (or even $\fl{u \times v}$) for brevity.

An important property, which follows from the behavior of sheaves with respect to coverings of type (\ref{type c}), is that, when $u$, $v$, $w$ are localization morphisms in $X$, $S$, $Y$,  $X'=X$, $u = \mathrm{Id}_X$, and the right hand square is \textit{cartesian} (with $a : gw \iso vg'$), then $(u,v,w;a,b)$ is an equivalence (\cite{Illusie12}, 1.11).

We will be mostly concerned with the case where $X$, $Y$, $S$ are schemes, $f$, $g$ morphisms of schemes, and we consider the corresponding morphisms of topoi of sheaves for the étale topology, still denoted $f$, $g$. The functoriality map (\ref{functoriality2}) will be only used when (\ref{functoriality1}) is a commutative diagram of schemes, except for a crucial construction in (\ref{2-commutative}).

Points of the topos $X \ori_S Y$, i.e., morphisms from the punctual topos $\mathrm{pt}$ to $X \ori_S Y$ are triples $(x, y, c)$, where $x : \mathrm{pt} \to X$, $y : \mathrm{pt} \to Y$ are points of $X$ and $Y$, and $c$ is a morphism $gy \to fx$. Recall (\cite{SGA4}, VIII, 7.9) that, when $T$ is a scheme, points of the étale topos $T$ correspond to usual geometric points of $T$, i.e., morphisms $t$ from the spectrum of a separably closed field to $T$, and morphisms $t \to s$ in $T$ to morphisms of the corresponding strictly localized schemes $T_{(t)} \to T_{(s)}$. 

\subsection{}\label{vanishing topos}
Let $f : X \to S$, $g : Y \to S$ be morphisms of schemes, and let $\mathrm{pr}_1 : X \times_S Y \to X$, $\mathrm{pr}_2 : X \times_S Y \to Y$ be the projections. As $f\mathrm{pr_1 = g\mathrm{pr}_2}$, by the universal property of the oriented product, we get  a diagram where the upper triangles commute:
\begin{equation}\label{extended vanishing topos diagram}
\xymatrix{{} & X \times_S Y \ar[d]^{\Psi_{f,g}} \ar[dl]_{\mathrm{pr}_1} \ar[dr]^{\mathrm{pr}_2} & {} \\
X \ar[dr]_f & X \ori_S Y \ar[l]_{p_1} \ar[r]^{p_2} & Y \ar[dl]^{g} \\
{} & S & {}}. 
\end{equation}
The inverse image functor by $\Psi_{f,g}$ is given by $\Psi_{f,g}^*(U \to V \leftarrow W) = U \times_V W$, in the notation of \ref{oriented product}. 

The most interesting case for us is when $Y = S$, $g = \mathrm{Id}_S$. The oriented product $X \ori_S S$ is then called the \textit{vanishing topos} of $f$. Diagram (\ref{extended vanishing topos diagram}) reduces to
\begin{equation}\label{vanishing topos diagram}
\xymatrix{{} & X \ar[d]^{\Psi_f} \ar[dl]_{\mathrm{Id}_X} \ar[dr]^f & {} \\
X \ar[dr]_f & X \ori_S S \ar[l]_{p_1} \ar[r]^{p_2} & S \ar[dl]^{\mathrm{Id}_S} \\
{} & S & {}}, 
\end{equation}
where $\Psi_f := \Psi_{f,g}$.  

Let $\Lambda$ be a ring. The direct image functor by $\Psi_f$, denoted
\begin{equation}
R\Psi_f := R(\Psi_f)_* : D^+(X,\Lambda) \to D^+(X \ori_S S, \Lambda),
\end{equation}
is called  the \textit{nearby cycles} functor (relative to $f$).

The canonical morphism $p_1^* \to (\Psi_f)_*$ gives rise to a functor, called the \textit{vanishing cycles} functor
\begin{equation}
R\Phi_f : D^+(X,\Lambda) \to D^+(X \ori_S S, \Lambda), 
\end{equation}
with a functorial exact triangle in $K \in D^+(X,\Lambda)$,
\begin{equation}\label{vanishing cycles triangle}
p_1^*K \to R\Psi_fK \to R\Phi_fK \to. 
\end{equation}

Actually, by the method of (\cite{SGA7}, XIII 1.4), one sees that (\ref{vanishing cycles triangle}) underlies an object
\begin{equation}\label{filtered Rpsi}
R\underline{\Psi}_fK
\end{equation}
of the filtered derived category $DF^{[0,1]}(X,\Lambda)$, with $F^1R\underline{\Psi}_fK = p_1^*K$ and $\mathrm{gr}^0_FR\underline{\Psi}_fK = R\Phi_fK$.

The relation between these objects and the classical nearby and vanishing cycles defined in (\cite{SGA7}, XIII) is as follows. Assume that $S$ is the spectrum of a henselian discrete valuation ring, with closed (resp. generic) point $s$ (resp. $\eta$). The topos denoted $X_s \times_s S$ in (\textit{loc. cit}. 1.2) is $X_s \ori_S S$, the union of the open subtopos $X_s \ori_S \eta$ and complementary  closed subtopos $X_s = X_s \ori_S s$. Sheaves on $X_s \ori_S S$ are described by triples $F = (F_s,F_{\eta},\varphi : F_{\overline{s}} \to F_{\overline{\eta}})$ (\cite{SGA7}, XIIII 1.2.4). In particular, we have a natural equivalence (given by $p_2 : s \ori_S S \to S$)
\begin{equation}\label{oriented topos of trait}
s \ori_S S \iso S,
\end{equation}
by which we will generally identify those two topoi (cf. (\cite{SGA7}, XIII 1.2.2 (b)). The functor $R\Psi_{\eta} : D^+(X_{\eta},\Lambda) \to D^+(X_s \ori_S \eta,\Lambda)$ (\cite{SGA7}, XIII 1.3.3, 2.1.2) is 
\begin{equation}\label{classical psi eta}
R\Psi_{\eta} : K \mapsto \fl i^*R\Psi_{f,{\eta}}K = (R\Psi_{f,\eta}K) | X_s \ori_S \eta 
\end{equation}
where $\fl i : X_s \ori_S  \eta \inj X \ori_S \eta$
is the closed subtopos defined by the inclusion $i$ of $X_s$ in $X$. The functor $R\Psi : D^+(X,\Lambda) \to D^+(X_s \ori_S S,\Lambda)$ (\cite{SGA7}, XIII, 1.3.3, 2.1.2) is given by 
\begin{equation}\label{classical psi}
R\Psi : K \mapsto (R\Psi_fK) | X_s \ori_S S,
\end{equation}
where $(R\Psi_fK) | X_s \ori_S S$ defines the triple $(K|X_s, R\Psi_{\eta}(K), \varphi)$. In particular, the functor $R\Phi_{\eta} : D^+(X,\Lambda) \to D^+(X_s \ori_{S} \eta,\Lambda)$ is given by 
\begin{equation}
R\Phi_{\eta} : K \mapsto  (R\Phi_fK) | X_s \ori_S \eta.
\end{equation}

\subsection{}\label{stalks}
We now recall the description of stalks of $R\Psi$ at points of the oriented product. In the situation of (\ref{extended vanishing topos diagram}), let $(x,y, c : g(y) \to f(x))$  be a point of $X \ori_S Y$ (\ref{oriented product}). Neighborhoods of $(x,y,c)$ consisting of $(U \to V \leftarrow W)$ with $U$, $V$, $W$ affine étale neighborhoods of $x$, $f(x)$, $y$ respectively, form a cofinal system, of which $X_{(x)} \times_{S_{(f(x)))}} Y_{(y)}$ is the projective limit, where $Y_{(y)} \to S_{(f(x))}$ is the composition of $Y_{(y)} \to S_{(g(y))}$ and the specialization $c : S_{(g(y))} \to S_{(f(x))}$. Therefore, by (\cite{SGA4}, VII, 5.8), we have, for $K \in D^+(X \times_S T,\Lambda)$,
\begin{equation}\label{Milnor tube}
R\Psi_{f,g}(K)_{(x,y,c)} = R\Gamma(X_{(x)} \times_{S_{(f(x)))}} Y_{(y)},K).
\end{equation}
In particular, for $Y = S$, $g = \mathrm{Id}_S$ in (\ref{extended vanishing topos diagram}), and a point $(x, c : t \to s)$ of $X \ori_S S$, we have
\begin{equation}\label{Milnor tube1}
R\Psi_f(K)_{(x,t \to s)} = R\Gamma(X_{(x)} \times_{S_{(s)}} S_{(t)},K),
\end{equation}
and an exact triangle
\begin{equation}\label{Milnor tube2}
K_x \to R\Gamma(X_{(x)} \times_{S_{(s)}} S_{(t)},K) \to R\Phi_f(K)_{(x,t \to s)} \to.
\end{equation}
Thus, for $X = S$, $f = \mathrm{Id}_S$, $x = s$, $R\Phi_{\mathrm{Id}_S}K$ measures the defect of the specialization maps $K_s \to K_t$ to be isomorphisms.

The fiber product $X_{(x)} \times_{S_{(s)}} S_{(t)}$ is called the \textit{Milnor tube} at $(x, c: t \to s)$, in contrast with the \textit{Milnor fiber} $X_{(x)} \times_{S_{(s)}} t$, a closed subscheme of the Milnor tube.

It is sometimes useful to consider, instead of Milnor tubes and fibers, Orgogozo's \textit{sliced} nearby (resp. vanishing cycles) $\Psi_{s,t}$ (resp. $\Phi_{s,t}$) (\cite{Orgogozo06}, proof of 6.1). Namely, for a specialization map $c : t \to s$ of geometric points of $S$, let $X_{(s)} := X \times_S S_{(s)}$, consider the inclusion $i_s : X_s \inj X_{(s)}$, and the morphism $j_{(t)} : X_{(t)} := X_{(s)} \times_{S_{(s)}} S_{(t)} \to X_{(s)}$ (with $S_{(t)} \to S_{(s)}$ given by $c$). Define the sliced nearby cycles functor 
\begin{equation}\label{slices}
R\Psi_{(f,c)} : D^+(X_{(t)},\Lambda) \to D^+(X_s,\Lambda)
\end{equation}
by
$$
R\Psi_{(f,c)}K = i_s^*Rj_{(t)_*}(K),
$$
and let us write $R\Psi_{(s,t)}$ for short when no confusion can arise. For $K \in D^+(X,\Lambda)$, we have
$$
R\Psi_{(s,t)}K = \fl i_{(s,t)}^*R\Psi_fK,
$$
where $\fl i_{(s,t)} : X_s = X_s \ori_s t \to X \ori_S S$ is the morphism given by $(i_s : X_s \to X, s \to S, t \to S, c : t \to s)$ (\ref{functoriality1}). The stalk of $R\Psi_{(s,t)}K$ at a geometric point $x$ above $s$ is $(R\Psi_fK)_{(x,c)}$ (\ref{Milnor tube1}). The sliced vanishing cycles functor $R\Phi_{s,t}$ is defined by the exact triangle
\begin{equation}\label{slices1}
K|X_s \to R\Psi_{s,t}K \to R\Phi_{s,t}K \to,
\end{equation}
where the first map is the adjunction map. We have $R\Phi_{(s,t)}K = \fl i_{(s,t)}^*R\Phi_fK$.
 
\medskip
\subsection{}\label{base change} 
A commutative diagram of schemes
\begin{equation}\label{base change diagram1}
\xymatrix{X' \ar[d]_{h'} \ar[r]^{f'} & S' \ar[d]_h \\
X \ar[r]^{f} & S}
\end{equation}
defines a morphism of oriented topoi (\ref{functoriality2} )
\begin{equation*}
h' \ori_{h} h : X' \ori_{S'} S' \to X \ori_S S,
\end{equation*}
hence a commutative diagram
\begin{equation}\label{base change diagram}
\xymatrix{X' \ar[d]_{\Psi_{f'}} \ar[r]^{h'} & X \ar[d]_{\Psi_{f}} \\
X' \ori_{S'} S' \ar[r]^{h' \ori_{h} h} & X \ori_S S}
\end{equation}
For $K \in D^+(X,\Lambda)$, we have a base change morphism
\begin{equation}\label{base change morphism}
(h' \ori_{h} h)^* R\Psi_fK \to R\Psi_{f'}h^*K.
\end{equation}

\begin{definition}\label{psi-good}
 We say that the formation of $R\Psi_fK$ commutes with base change, or that the pair $(f,K)$ is \textit{$\Psi$-good}, if, for any cartesian diagram (\ref{base change diagram1}), (\ref{base change morphism}) is an isomorphism.
\end{definition}

Examples will be discussed in \ref{psi-good-examples}, after we have recalled results of Orgogozo in \cite{Orgogozo06}.

\medskip
\subsection{}\label{constructibility}
In the situation of (\ref{extended vanishing topos diagram}), assume that $S$ is noetherian, and $f$ and $g$ are of finite type. It is shown in (\cite{Illusie12}, 2.5) that $X \ori_S Y$ is then a coherent topos and the projections $p_1$, $p_2$ are coherent morphisms (\cite{SGA4}, VI, 2.3, 2.4.5, 3.1). Moreover, by (\cite{Orgogozo06}, 9) one has a good notion of constructibility for sheaves on $X \ori_S Y$. Let $\Lambda$ be a noetherian ring. A sheaf of $\Lambda$-modules $F$ on $X \ori_S Y$ is called \textit{constructible} if it has a presentation $L \to M \to F \to 0$, where $L$ and $M$ are finite sums of $\Lambda$-modules of the form $\Lambda^{(U \to V \leftarrow W)}$, for $(U \to V \leftarrow W)$ objects of the defining site of $X \ori_S Y$, with $U \to X$, $V \to S$, $W \to Y$ separated, étale, and of finite presentation.  It is equivalent to saying that there exist finite partitions of $X$ and $Y$ into locally closed subsets: $X = \cup_{\alpha \in A} X_{\alpha}$, $Y = \cup_{\beta \in B} Y_{\beta}$, such that, for all $(\alpha,\beta)$, the restriction of $F$ to the subtopos $X_{\alpha} \ori_S Y_{\beta}$ is locally constant of finite type. As in the case of noetherian schemes, constructible sheaves of $\Lambda$-modules form a thick subcategory (in the strong sense, i.e., closed under subobjects, quotients and extensions) so that in particular the full subcategory $D^+_c(X \ori_S Y,\Lambda)$ of $D^+(X \ori_S Y,\Lambda)$ consisting of complexes $K$ such that $\mathcal{H}^iK$ is constructible for all $i$ is a triangulated subcategory. 

The following theorem (\cite{Orgogozo06}, 2.1, 3.1, 8.1) is the main result of \cite{Orgogozo06}. 

\medskip
\begin{theorem2}\label{Orgogozo's theorem}
In the situation of (\ref{vanishing topos diagram}), assume that $S$ is noetherian, $f$ of finite type, and $\Lambda$ is a $\mathbf{Z}/n\mathbf{Z}$-algebra, with $n$ invertible on $S$. Let $K \in D^b_c(X,\Lambda)$. Then there exists a modification $a : S' \to S$ such that if $X' = X \times_S S'$ and $K'$ is the inverse image of $K$ on $X'$, the following conditions are satisfied :

(i) $R\Psi_{f'}K'$ belongs to $D^b_c(X' \ori_{S'} S',\Lambda)$ ;

(ii) the formation of $R\Psi_{f'}K'$ commutes with any base change $S'' \to S'$, i.e., $(f',K')$ is $\Psi$-good (\ref{psi-good}).
\end{theorem2}
\medskip
In (\textit{loc. cit.}) $\Lambda = \mathbf{Z}/n\mathbf{Z}$, but the proof goes on without any change for $\Lambda$ as above. Note that (ii) implies (i) by (\cite{Orgogozo06} 8.1, 10.5).

\begin{example}\label{psi-good-examples} Let $\Lambda$ be as in \ref{Orgogozo's theorem}, except in (b) where it can be any ring.

(a) If $S$ is regular of dimension $\le 1$, $(f,K)$ is $\Psi$-good (as any modification $S' \to S$ has a section). This contains, in particular, the universal local acyclicity of pairs $(X,K)$ for $X$ of finite type over a field (\cite{SGA41/2}, Th. finitude, 2.16), and the compatibility of classical nearby cycles with base change by surjective maps of traits (\cite{SGA41/2}, Th. finitude, 3.7).

(b) A pair $(f,K)$ is locally acyclic (\cite{SGA41/2}, Th. finitude, 2.12) if and only if $R\Phi_fK = 0$ and $R\Phi_fK$ commutes with locally quasi-finite base change. In particular, $(f,K)$ is universally locally acyclic if and only if $(f,K)$ is $\Psi$-good and $R\Phi_fK = 0$, i. e. $R\Phi_{f'}(K') = 0$ for all $S' \to S$, where $K' = K|S'$ and $f' = f \times_S S'$. 

Indeed, if $(f,K)$ is locally acyclic, then, by (\cite{SGA41/2}, Th. finitude, Appendice, 2.6), for any geometric point $x$ of $X$ with image $s$ in $S$, if $f_{(x,s)} : X_{(x)} \to S_{(s)}$ is the morphism deduced by localization, the formation of $Rf_{(x,s)*}K$ (where $K$ denotes again, by abuse, its restriction to $X_{(x)}$) commutes with any finite base change (hence $R\Psi_fK$ and $R\Phi_fK$ commute with any locally quasi-finite base change, cf. A.2 (2)). In particular, for any geometric point $t \to S_{(s)}$, the restriction map $R\Gamma(X_{(x)} \times_{S_{(s)}} S_{(t)},K) \to R\Gamma(X_{(x)_t},K)$ from the Milnor tube to the Milnor fiber (cf. (\ref{Milnor tube})) is an isomorphism, as it is the stalk at $t$ of the base change map by the closure of the image of $t$ in $S_{(s)}$. Therefore, as the composition $K_x \to R\Gamma(X_{(x)} \times_{S_{(s)}} S_{(t)},K) \to R\Gamma(X_{(x)_t},K)$ is an isomorphism by the definition of local acyclicity, the first map of this composition is an isomorphism as well, i. e. $(R\Phi_fK)_{(x,t)} = 0$. That proves the ``only if" part of the assertion. The converse is immediate. This argument is copied from (\cite{TSaito15a}, Prop. 1.7).

(c) By (\cite{Orgogozo06}, 6.1), if $(f,K)$ is universally locally acyclic on an open subset $U$ of $X$ whose complement is quasi-finite over $S$, then $(f,K)$ is $\Psi$-good. 

(d) If $f : X \to S$ is the blow-up of the origin in $S = \mathbf{A}^2_k$, $k$ an algebraically closed field, $\Lambda = \mathbf{Z}/n\mathbf{Z}$, then $(f,\Lambda)$ is not $\Psi$-good (\cite{Orgogozo06}, 11).
\end{example}

\medskip
We will see other examples of $\Psi$-goodness in the next section (2.3).

\subsection{}\label{cohomological descent}
The results of the next two subsections will be used only in the proof of \ref{kunneth theorem}.

There are variants of the constructions in \ref{vanishing topos} for morphisms of simplicial schemes $f_{\bullet} : X_{\bullet} \to S_{\bullet}$, $g_{\bullet} : Y_{\bullet} \to S_{\bullet}$, see (\cite{Orgogozo12b}, 2.2). We will use them freely, and sometimes, by abuse, write $X_{\bullet} \ori_{S_{\bullet}} Y_{\bullet}$ for the total topos of this simplicial topos. 

Let $\Lambda$ be as in \ref{Orgogozo's theorem}. It follows from \ref{Orgogozo's theorem} that there exists a hypercovering for the h-topology on $S$ (\cite{Orgogozo12b}, 2.1.3) $a : T_{\bullet} \to S$ such that for all $n \in \mathbf{N}$, $R\Psi_{f_n}(K | X_{T_n})$ belongs to $D^b_c(X_{T_{n}} \ori_{T_n} T_n, \Lambda)$ and is of formation compatible with base change. Such a hypercovering will be called \textit{admissible for $(f,K)$}. Admissible hypercoverings form a cofinal system in the category of hypercoverings of $S$ for the h-topology, up to homotopy.

Note that, for \textit{any} hypercovering $a : T_{\bullet} \to S$ for the h-topology, $R\Psi_f(K)$ is recovered from $a^*K$ by a canonical isomorphism
\begin{equation}\label{recovering psi}
R\Psi_f(K) \iso R\fl a_*R\Psi_{f_{T_{\bullet}}}(K | X_{T_{\bullet}}), 
\end{equation}
where $f_{T_{\bullet}} : X_{T_{\bullet}} \to T_{\bullet}$ is the morphism deduced by base change, and $\fl a : X_{T_{\bullet}} \ori_{T_{\bullet}} T_{\bullet} \to X \ori_S S$ is the corresponding augmentation\footnote{I am indebted to Gabber for this observation.}. Indeed, we have a commutative diagram (cf. (\ref{base change diagram}))
\begin{equation}\label{descent diagram}
\xymatrix{X_{T_{\bullet}} \ar[d]^{\Psi_{f_{T_{\bullet}}}} \ar[r]^{a} & X \ar[d]^{\Psi_f} \\
X_{T_{\bullet}} \ori_{T_{\bullet}} T_{\bullet} \ar[r]^{\fl a} & X \ori_S S},
\end{equation}
and by cohomological descent the adjunction map $K \to Ra_*a^*K$ is an isomorphism. The isomorphism (\ref{recovering psi}) is induced by this isomorphism via (\ref{descent diagram}).  

\subsection{}\label{oriented descent}
By (\cite{Orgogozo12b}, 2.2.3), $\fl a : X_{T_{\bullet}} \ori_{T_{\bullet}} T_{\bullet} \to X \ori_S S$ is of 1-cohomological descent, which means that, for any $L \in D^b(X \ori_S S,\Lambda)$, the adjunction map 
\begin{equation}\label{oriented descent iso}
L \to R\fl a_*\fl a^*L
\end{equation}
is an isomorphism. For $L = R\Psi_fK$, we have a commutative diagram
$$
\xymatrix{{} & R\fl a_*\fl a^*R\Psi_fK \ar[dd] \\
R\Psi_fK \ar[ur] \ar[dr] \\
{} & R\fl a_*R\Psi_{f_{T_{\bullet}}}(K|X_{T_{\bullet}})}
$$
in which the vertical map 
\begin{equation}\label{vertical map}
R\fl a_*\fl a^*R\Psi_f(K) \to R\fl a_*R\Psi_{f_{T_{\bullet}}}(K|X_{T_{\bullet}})
\end{equation}
is $R\fl a_*$ applied to the base change map 
$$
\beta : \fl a^* R\Psi_fK \to R\Psi_{f_{T_{\bullet}}}(K|X_{T_{\bullet}}),
$$
the upper oblique map is the isomorphism (\ref{oriented descent iso}) and the lower one the isomorphism (\ref{recovering psi}). Therefore (\ref{vertical map}) is an isomorphism. However, in general $\beta$ is not an isomorphism, as in general the formation of $R\Psi_fK$ is not compatible with base change.

\medskip
In the rest of this section we collect a few general facts that will be used in sections 3 and 4.

\subsection{}\label{composition}
For morphisms of schemes $f : X \to Y$, $g : Y \to Z$, the diagram (of type (\ref{functoriality1}))
\begin{equation}\label{composition diagram}
\xymatrix{X \ar[dr]_{gf} \ar[r]^f & Y \ar[d]_{g} & Y \ar[l]^{\mathrm{Id}} \ar[d]_g \\
{} & Z & Z \ar[l]^{\mathrm{Id}}}
\end{equation}
induces a morphism
\begin{equation*}
\fl{g} := \fl{\mathrm{Id_X} \times g} : X \ori_Y Y \to X \ori_Z Z,
\end{equation*}
hence a commutative diagram of type (\ref{base change diagram})
\begin{equation}\label{composite psi}
\xymatrix{X \ar[r]^{\Psi_f} \ar[dr]_{\Psi_{gf}} & X \ori_Y Y \ar[d]^{\fl g} \ar[r]^{p_1} & X \\
{} & X \ori_{Z} Z \ar[ur]_{q_1} & {}}, 
\end{equation}
where we denote by $q_1 : X \ori_Z Z \to X$, instead of $p_1$, the canonical projection. 

Let $\Lambda$ be a commutative ring. For $K \in D^+(X,\Lambda)$, (\ref{composition diagram}) induces an isomorphism
\begin{equation}\label{psi of composite}
R\fl g_* R\Psi_fK \iso R\Psi_{gf}K. 
\end{equation}
With this identification the base change map $q_1^*K \to R\fl g_*p_1^*K$ associated with the right triangle of (\ref{composite psi}) gives a map
\begin{equation}\label{phi of composite}
R\Phi_{gf}K \to R\fl g_*R\Phi_fK, 
\end{equation}
which is not an isomorphism in general, as the trivial case where $X = Y$, $f = \mathrm{Id}_X$, $K = \Lambda$ already shows ($R\Phi_{\mathrm{Id_X}}\Lambda = 0$). In \ref{phi of composite-cone}, we give a formula for the cone of (\ref{phi of composite}), under certain hypotheses. 

\subsection{}\label{local sections}
Let $f : X \to Y$ be a morphism of schemes, and $K \in D^+(X,\Lambda)$. We have seen how to calculate $R\Psi_f(K)$ by \textit{slices} (\ref{slices}). Another way to unravel $R\Psi_fK$ is to consider its restrictions to \textit{local sections}. 

Let $x : \SP \, k  \to X$ be a geometric point of $X$, with image the geometric point $y = f(x)$ of $Y$,  and let $X_{(y)} := X \times_Y Y_{(y)}$. The definition of oriented products implies (cf. (\cite{Illusie12}, 1.4, 1.11) that the morphisms of topoi
\begin{equation}\label{stalk of oriented product}
\SP \,k \ori_Y Y \to \SP \, k \ori_{Y_{(y)}} Y_{(y)} \to Y_{(y)}
\end{equation}
are isomorphisms, where the second map is given by $p_2$ and the first one is the inverse of the localization map (\cite{Illusie12}, 1.11). This generalizes the strictly local case of (\ref{oriented topos of trait}). We will write $x \ori_Y Y$ for $\SP \,k \ori_Y Y$, and call it the \textit{stalk of $X \ori_Y Y$ at $x$}. It is constant along $X_y$, of value $Y_{(y)}$.

Recall (\cite{Illusie12} 2.2) that, as the topos $Y_{(y)}$ is local, the point $x$ defines a canonical section
\begin{equation}\label{canonical section}
\sigma_{x,y} : Y_{(y)}  \to X_{(y)} \ori_{Y_{(y)}} Y_{(y)} 
\end{equation} 
of the projection $p_2 : X_{(y)} \ori_{Y_{(y)}} Y_{(y)} \to Y_{(y)}$. 

This section sends a point $(t \to y)$ of $Y_{(y)}$ to the point $(x \to y \la t)$ of $X_{(y)} \ori_{Y_{(y)}} Y_{(y)}$. Formally, $\sigma_{x,y}$ is defined as follows. Let $\varepsilon : Y_{(y)} \to \mathrm{pt}$ be the projection on the punctual topos $\mathrm{pt}$. Consider the 2-map 
\begin{equation*}
c_y : \mathrm{Id}_{Y_{(y)}} \to y\varepsilon
\end{equation*}
given by adjunction map $(y\varepsilon)^*K = (K_y)_{Y_{(y)}} = \Gamma(Y_{(y)},K)_{Y_{(y)}} \to K$, for $K$ a sheaf on $Y_{(y)}$. We have a diagram
\begin{equation}\label{2-commutative}
\xymatrix{\mathrm{pt} \ar[d]_x \ar[r]^{\mathrm{Id}} & \mathrm{pt} \ar[d]_y & Y_{(y)} \ar[l]_{\varepsilon} \ar[d]_{\mathrm{Id}} \\
X_{(y)} \ar[r]^{f_{(y)}} & Y_{(y)} & Y_{(y)} \ar[l]_{\mathrm{Id}}},
\end{equation}
in which the left hand square trivially commutes, and the right hand one is 2-commutative by the datum of the 2-map $c_y$. The canonical section $\sigma_{x,y}$ (\ref{canonical section}) is the functoriality map on the oriented products given by (\ref{2-commutative}) (with the identification of $Y_{(y)}$ with $x \ori_Y Y \iso x \ori_{Y_{(y)}} Y_{(y)}$). It sits in a commutative diagram
\begin{equation}\label{canonical section2}
\xymatrix{\mathrm{pt} \ar[d]_x & \mathrm{pt} \ori_{\mathrm{pt}} Y_{(y)} \ar[l]_{p_1} \ar[r]^{p_2} \ar[d]_{\sigma_{x,y}} & Y_{(y)} \ar[d]_{\mathrm{Id}} \\
X_{(y)} & X_{(y)} \ori_{Y_{(y)}} Y_{(y)} \ar[r]^{p_2} \ar[l]_{p_1} & Y_{(y)}},
\end{equation}
where the upper $p_2$ is one of the isomorphisms (\ref{stalk of oriented product}). When no confusion can arise, we will also denote by $\sigma_{x,y}$ the composition 
\begin{equation}\label{extended canonical section} 
\sigma_{x,y} : Y_{(y)} \to X \ori_Y Y
\end{equation}
of (\ref{canonical section}) and the canonical map  $X_{(y)} \ori_{Y_{(y)}} Y_{(y)} \to X \ori_Y Y$, and we will write $\sigma_x$ for $\sigma_{x,y}$. For $L \in D^+(X \ori_Y Y,\Lambda)$, we will sometimes write
\begin{equation}\label{canonical section1}
L | x \ori_Y Y = L | Y_{(y)} := \sigma_{x,y}^*L \in D^+(Y_{(y)},\Lambda)
\end{equation} 
for the restriction (via $\sigma_{x,y}$) of $L$ to the stalk of $X \ori_Y Y$ at $x$.

By construction, as the point $x : \mathrm{pt} \to X_{(y)}$ factors through $X_{(x)}$, so does $\sigma_{x,y}$: we have a commutative diagram
\begin{equation}\label{canonical section localized}
\xymatrix{{} & X_{(x)} \ori_{Y_{(y)}} Y_{(y)} \ar[d] \ar[dr]_{p_2} & {}\\
Y_{(y)} \ar[ur]^{\widetilde{\sigma}_{x,y}} \ar[r]^{\sigma_{x,y}} & X_{(y)} \ori_{Y(y)} Y_{(y)} \ar[r]^{p_2} &Y_{(y)}},
\end{equation}
where the composition of the horizontal maps is the identity. Recall the following result, due to Gabber (a special case of (\cite{Illusie12}, 2.3.1)):
\begin{proposition}\label{Gabber's isom} The canonical map 
\begin{equation}\label{Gabber's isom1}
\gamma : p_{2*} \to \widetilde{\sigma}_{x,y}^*,
\end{equation}
defined by $p_2 \widetilde{\sigma}_{x,y} = \mathrm{Id}$, is an isomorphism. 
\end{proposition}
It follows that, if $f_{(x,y)} : X_{(x)} \to Y_{(y)}$ is the map deduced from $f$ by localization at $x$ and $y$, by the commutative diagram (upper right triangle of (\ref{vanishing topos diagram})) 
\begin{equation}\label{section iso 2}
\xymatrix{X_{(x)} \ar[d]_{\Psi_{f_{(x,y)}}} \ar[dr]^{f_{(x,y)}} \\
X_{(x)} \ori_{Y_{(y)}} Y_{(y)} \ar[r]^{p_2} & Y_{(y)}},
\end{equation} 
$\gamma$ induces an isomorphism:
\begin{equation}\label{section iso 3}
\gamma_{x,y} : \widetilde{\sigma}_{x,y}^*R\Psi_{f_{(x,y)}}(L) \iso p_{2*}R\Psi_{f_{(x,y)}}(L) \iso Rf_{(x,y)*}(L)
\end{equation}
for any $L \in D^+(X_{(x)},\Lambda)$. On the other hand, for $K \in D^+(X,\Lambda)$, by (\ref{Milnor tube1}) the base change map
\begin{equation}\label{localization}
(R\Psi_f K)| X_{(x)} \ori_{Y_{(y)}} Y_{(y)} \to R\Psi_{f_{(x,y)}}(K |X_{(x)})
\end{equation}
defined by the commutative diagrams
\begin{equation}\label{localization diagram}
\xymatrix{X_{(x)} \ar[r] \ar[d]_{\Psi_{f_{(x,y)}}} & X_{(y)} \ar[d]_{\Psi_{f_{(y)}}} \ar[r] & X \ar[d]_{\Psi_f} \\
X_{(x)} \ori_{Y_{(y)}} Y_{(y)} \ar[r] & X_{(y)} \ori_{Y_{(y)}} Y_{(y)}  \ar[r] & X \ori_Y Y }, 
\end{equation}
is an isomorphism, in which the horizontal arrows are defined by the localization maps. 
By composing $\gamma_{x,y}$ (for $L = K|X_{(x)}$) and (\ref{localization}), we get an isomorphism
\begin{equation}\label{section iso}
\gamma_{x,y} : \sigma_{x,y}^*R\Psi_fK  \iso Rf_{(x,y)*}(K |X_{(x)}), 
\end{equation}
where $\sigma_{x,y}$ is the composition of $\widetilde{\sigma}_{x,y}$ and the localization map $X_{(x)} \ori_{Y_{(y)}} Y_{(y)} \to X \ori_Y Y$. 

At a geometric point $t$ of $Y_{(y)}$ (\ref{section iso}) induces an isomorphism
\begin{equation}
(\sigma_x^*R\Psi_fK )_t \iso R\Gamma(X_{(x)} \times_{Y_{(y)}} Y_{(t)},K) = (R\Psi_fK)_{(x,y \la t)}, 
\end{equation}
and in particular, an isomorphism
\begin{equation}\label{stalk at y section iso}
(\sigma_x^*R\Psi_fK )_y \iso R\Gamma(X_{(x)},K) = K_x. 
\end{equation}
The specialization morphism $K_x \to (\sigma_x^*R\Psi_fK )_t$ is identified with the stalk at $t$ of $\sigma_x^*(p_1^*K \to R\Psi_fK)$, so that, taking into account the fact that $p_1\widetilde{\sigma}_{x,y} : Y_{(y)} \to X_{(x)}$ projects $Y_{(y)}$  onto the closed point of $X_{(x)}$, (\ref{vanishing cycles triangle}) induces a distinguished triangle of $D^+(Y_{(y)},\Lambda)$
\begin{equation}\label{stalk at y vanishing cycles triangle}
(K_x)_{Y_{(y)}} \to \sigma_x^*R\Psi_fK \to \sigma_x^*R\Phi_fK \to, 
\end{equation} 
where $(K_x)_{Y_{(y)}}$ denotes the \textit{constant complex} on $Y_{(y)}$ of value $K_x$. Note that taking the stalk at $y$ of (\ref{stalk at y vanishing cycles triangle}) and using (\ref{stalk at y section iso}), we get
\begin{equation}\label{vanishing of Phi}
(\sigma_x^*R\Phi_fK)_y = 0.
\end{equation}

Thus, locally at $y$, $R\Psi_fK$ defines a family of complexes 
$$
\sigma_x^*R\Psi_fK \iso Rf_{(x,y)*}(K|X_{(x)})
$$
of $D^+(Y_{(y)},\Lambda)$, parametrized by the geometric points $x$ of $X$ above $y$, generalizing the classical $(R\Psi_fK)_x$ when $Y$ is a trait, with geometric closed point $y$ (\cite{SGA7}, XIII, 2.1.1).  

\begin{proposition}\label{transitivity}
Let $f : X \to Y$, $g : Y \to Z$, and $\fl g : X \ori_Y Y \to X \ori_Z Z$ be as in \ref{composition}. Let $x$ be a geometric point of $X$ with images $y = f(x)$ in $Y$ and $z = g(y)$ in $Z$. For $L \in D^+(X \ori_Y Y,\Lambda)$, the commutative diagram
\begin{equation}\label{transitivity diagram}
\xymatrix{x \ori_Y Y = Y_{(y)} \ar[d]_{g_{(y,z)}} \ar[r]^{\sigma_{x,y}} & X \ori_Y Y \ar[d]^{\fl g} \\
x \ori_Z Z = Z_{(z)} \ar[r]^{\sigma_{x,z}} & X \ori_Z Z}
\end{equation}
gives an isomorphism
\begin{equation}\label{stalk of Rg*}
(R\fl g_*L) | x \ori_Z Z \iso Rg_{(y,z)*}(L | x \ori_Y Y)
\end{equation}
\end{proposition}
\begin{proof} The map (\ref{stalk of Rg*}) is the base change map associated with the square (\ref{transitivity diagram}). This square is the composite of the following squares
\begin{equation}\label{square1}
\xymatrix{Y_{(y)} \ar[d] \ar[r]^{\widetilde{\sigma}_{x,y}} & X_{(x)} \ori_{Y_{(y)}} Y_{(y)} \ar[d] \ar[r] & X_{(x)} \ori_Y Y \ar[d] \ar[r] & X \ori_Y Y \ar[d] \\
Z_{(z)} \ar[r]^{\widetilde{\sigma}_{x,z}} & X_{(x)} \ori_{Z_{(z)}} Z_{(z)} \ar[r] & X_{(x)} \ori_Z Z \ar[r] & X \ori_Z Z}. 
\end{equation}
Base change in the right square is clear, as $X_{(x)} \to X$ is a limit of étale neighborhoods of $x$. Base change in the left square follows from Gabber's formula $p_{2*} = \widetilde{\sigma_x}^*$ (\ref{Gabber's isom}) applied to $\widetilde{\sigma}_{x,y}$ and $\widetilde{\sigma}_{x,z}$. It remains to show base change in the middle square. Define $Y_{(z)}$ by the cartesian square
\begin{equation*}
\xymatrix{Y_{(z)} \ar[r] \ar[d]  & Y \ar[d] \\
Z_{(z)} \ar[r] & Z}.
\end{equation*}
The middle square is thus decomposed into 
\begin{equation*}
\xymatrix{X_{(x)} \ori_{Y_{(y)}} Y_{(y)} \ar[dr] \ar[r]^{u} & X_{(x)} \ori_{Y_{(z)}} Y_{(z)} \ar[d] \ar[r] & X_{(x)} \ori_Y Y \ar[d] \\
{} & X_{(x)} \ori_{Z_{(z)}} Z_{(z)} \ar[r] & X_{(x)} \ori_Z Z},
\end{equation*}
where $u$ is the localization map defined by $Y_{(y)} \to Y_{(z)}$. Base change in the square is clear, as $Y_{(z)} \to Z_{(z)}$ is a limit of pull-backs of $Y \to Z$ by étale neighborhoods of $z$ in $Z$. It then suffices to show base change in the triangle.  Let $Y_1 := Y_{(z)}$. We have $Y_{(y)} = (Y_1)_{(y)}$. By (\cite{IllusieLaszloOrgogozo14}, XI 1.11) applied to the case the map $g$ of \textit{loc. cit.} is the identity (and $Y$ of \textit{loc. cit.} our $Y_1$), and passing to the limit on the étale neighborhoods of $y$ in $Y_1$, $u$ is an equivalence, hence base change in the triangle is trivial, which finishes the proof.
\end{proof}

In particular, via the isomorphisms $\gamma$ (\ref{section iso}) and the commutative square
\begin{equation*}
\xymatrix{X_{(x)} \ori_{Y_{(y)}}\ar[d]_{\fl g} Y_{(y)} \ar[r]^{p_2} & Y_{(y)}\ar[d]_{g_{(y,z)}}\\
X_{(x)} \ori_{Z_{(z)}} Z_{(z)} \ar[r]^{p_2} & Z_{(z)}},
\end{equation*}
the restriction to the stalk of $X \ori_Z Z$ at $x$ of the transitivity isomorphism (\ref{psi of composite}),
\begin{equation*}
\sigma_{x,z}^*R\Psi_{gf}K \iso \sigma_{x,z}^*(R\fl g_* R\Psi_fK),
\end{equation*}
 translates into the isomorphism
\begin{equation}\label{transitivity 1}
R(gf)_{(x,z)*}K \iso Rg_{(y,z)*}Rf_{(x,y)*}K ,
\end{equation}
coming from $g_{(y,z)}f_{(x,y)} = (gf)_{(x,z)}$.

\medskip
Recall the following result of Gabber (\cite{Orgogozo06}, 3.1):
\begin{theorem}\label{cd of psi}
Assume that $\Lambda$ is a $\mathbf{Z}/n\mathbf{Z}$-algebra for some integer $n \ge 1$. Let $f : X \to S$ be a morphism locally of finite type such that the dimension of its fibers is bounded by an integer $N$. Then, for any sheaf $\mathcal{F}$ of $\Lambda$-modules on $X$, we have
\begin{equation*}
R^q\Psi_f(\mathcal{F}) = 0
\end{equation*}
for all $q > 2N$.
\end{theorem}

\begin{corollary}\label{cd of localized map} 
Under the assumptions of \ref{cd of psi}, for any geometric point $x$ of $X$, with image $s$ in $S$, and any sheaf of $\Lambda$-modules $\mathcal{\mathcal{G}}$ on $X_{(x)}$, we have
\begin{equation*}
R^qf_{(x,s)*}\mathcal{G} = 0
\end{equation*}
for all $q > 2N$.
\end{corollary}
\begin{proof} By standard limit arguments (\cite{SGA4} VII 5.11, IX 2.7.2, 2.7.4) we may assume that $G = \mathcal{F}|X_{(x)}$ for some sheaf of $\Lambda$-modules $\mathcal{F}$ on $X$. Then the conclusion follows from \ref{cd of psi} by (\ref{section iso}).
\end{proof}

\begin{proposition}\label{phi of composite-cone} Let $\Lambda$ be as in \ref{cd of psi}, and let $f : X \to Y$, $g : Y \to Z$ be morphisms of schemes. Assume that $g$ is locally of finite type. Let $K \in D^+(X,\Lambda)$. Then, with the notation of (\ref{composite psi}), the canonical map (\ref{phi of composite}) fits in a distinguished triangle
\begin{equation}\label{cone}
R\Phi_{gf}K \to R\fl{g}_*R\Phi_fK \to q_1^*K \otimes^L \fl{(f \times \mathrm{Id}_Z)}^*R\Phi_g(\Lambda)[1] \to.
\end{equation}
\end{proposition}
\begin{proof}
Consider the composition given by (\ref{composite psi}):
\begin{equation*}
q_1^*K \to R\fl{g}_*p_1^*K \to R\fl{g}_*R\Psi_fK. 
\end{equation*}
By the usual method, we can view it as making $R\fl{g}_*R\Psi_fK$ (which, by (\ref{psi of composite}) is identified with $R\Psi_{gf}K$) into an object of $DF^{[0,1]}(X \ori_Z Z,\Lambda)$, the filtered derived category of complexes with filtration of length one.  By (\cite{Illusie03}, 4.1), it gives rise to a cross (i.e., a refined special 9-diagram):
\begin{equation}\label{cross}
\xymatrix{{} & R\fl{g}_*p_1^*K \ar[d] \ar[r] &C \ar[d] \\
q_1^*K \ar[r] \ar[ur] &R\Psi_{gf}K \ar[r] \ar[d] & R\Phi_{gf}K \ar[dl]  \\
{} &R\fl{g}_*R\Phi_fK } ,
\end{equation}
where
\begin{equation}\label{cone1}
C \simeq \mathrm{Cone}(q_1^*K \to R\fl{g}_*p_1^*K) \simeq \mathrm{Cone}(R\Phi_{gf}K \to R\fl{g}_*R\Phi_fK)[-1].
\end{equation}
By \ref{cd of localized map}, $g$ satisfies the condition (*) of A.1.  Therefore, by A.7 for $K' = \Lambda$, the projection formula map
\begin{equation*}
R\fl{g}_*\Lambda \otimes^L q_1^*K \to R\fl{g}_*p_1^*K
\end{equation*}
is an isomorphism, and the upper oblique map in (\ref{cross}) is identified with $q_1^*K \otimes^L \alpha$, where $\alpha : \Lambda \to R\fl{g}_*\Lambda$ is the adjunction map. Therefore, we get
\begin{equation}\label{cone2}
C \simeq q_1^*K \otimes^L \mathrm{Cone}(\alpha : \Lambda \to R\fl{g}_*\Lambda). 
\end{equation}
Consider the commutative diagram (cf. (\cite{TSaito16}, proof of 2.1):
\begin{equation*}
\xymatrix{{} & Y \\
X \ori_Y Y \ar[ur]^{p_{2}} \ar[d]_{\fl{\mathrm{Id}_X \times g}} \ar[r]^{\fl{f \times \mathrm{Id}_Y}} & Y \ori_Y Y\ar[u]^{p_{2}} \ar[d]_{\fl{\mathrm{Id}_Y \times g}} \\
X \ori_Z Z \ar[r]^{\fl{f \times \mathrm{Id}_Z}} & Y \ori_Z Z}.
\end{equation*}
For $M \in D^+(Y,\Lambda)$, the induced base change map
\begin{equation*}
(\fl{f \times \mathrm{Id}_Z})^* R(\fl{\mathrm{Id}_Y \times g})_*p_2^*M \to R(\fl{\mathrm{Id}_X \times g})_*p_2^*M
\end{equation*}
is an isomorphism, as one sees by applying $\sigma_{x,z}$ and using \ref{transitivity}, both sides being identified with $Rg_{(y,z)*}(M|Y_{(y)})$. Note that
\begin{equation*}
R(\fl{\mathrm{Id}_Y \times g})_*p_2^*M \simeq R\Psi_g(M)
\end{equation*}
by (\ref{psi of composite}) and $R\Psi_{\mathrm{Id}_Y} = p_2^*$. Finally, for $M = \Lambda$, we get
\begin{equation*}
R\fl{g}_*\Lambda \simeq (\fl{f \times \mathrm{Id}_Z})^*R\Psi_g\Lambda,
\end{equation*}
and $\alpha : \Lambda \to R\fl{g}_*\Lambda$ is identified with the canonical map $(\fl{f \times \mathrm{Id}_Z})^*(\Lambda \to R\Psi_g\Lambda)$, so we get (\ref{cone}) from (\ref{cone1}) and (\ref{cone2}).
\end{proof}

\begin{corollary}\label{locally acyclic} Under the assumptions of  \ref{phi of composite-cone}, let $x$ be a geometric point of $X$, with images $y$ in $Y$ and $z$ in $Z$. The following conditions are equivalent :

(i) The morphism 
\begin{equation}\label{stalk at x of RPhi(gf)}
\sigma_{x,z}^*R\Phi_{gf}K \to \sigma_{x,z}^*(R\fl g)_*R\Phi_fK 
\end{equation}
induced by (\ref{phi of composite}) is an isomorphism.

(ii) The cone of $K_x \otimes^L (\alpha : \Lambda_{Z_{(z)}} \to Rg_{(y,z)*}\Lambda_{Y_{(y)}})$, where $\alpha$ is the adjunction map, is zero, $K_x$ denoting the constant complex of value $K_{x}$ on $Z_{(z)}$.

In particular, if $g$ is locally acyclic at $y$ (\cite{SGA41/2}, Th. finitude, 2.12), e. g. smooth at $y$, these conditions are satisfied.
\end{corollary}
\begin{proof} It suffices to apply $\sigma_{(x,z)}^*$ to (\ref{cone}), as 
\begin{equation*}
\sigma_{(x,z)}^*(q_1^*K \otimes^L \fl{(f \times \mathrm{Id}_Z)}^*R\Phi_g(\Lambda)) \iso \mathrm{Cone}(K_x \otimes^L (\alpha : \Lambda_{Z_{(z)}} \to Rg_{(y,z)*}\Lambda_{Y_{(y)}})),
\end{equation*}
as we have seen in the proof of \ref{phi of composite-cone}. The last assertion follows from (ii)  by the definition of local acyclicity (or from (\ref{cone}) by (\ref{psi-good-examples} (b))). 
\end{proof}

\medskip
By (\ref{section iso}) (applied to $g$ and $\sigma_{y,z}$), for $L = \sigma_{x,y}^*R\Phi_fK = R\Phi_f(K)|Y_{(y)}$, we have an exact triangle
\begin{equation}\label{local phi}
(L_y)_{Z_{(z)}} \to Rg_{(y,z)*}L \to \sigma_{x,z}^*(R\Phi_{g_{(y,z)}}L) \to,
\end{equation}
where as before $(L_y)_{Z_{(z)}}$ denotes the constant complex of value $L_y = (Rg_{(y,z)*}L)_z = R\Gamma(Y_{(y)},L)$, and the first map is the specialization map. Taking into account that $(\sigma_{x,y}^*R\Phi_fK)_y = 0$ (\ref{vanishing of Phi}), it gives an isomorphism:
$$
Rg_{(y,z)*}\sigma_{x,y}^*R\Phi_fK \iso R\Phi_{g_{(y,z)}}\sigma_{x,y}^*R\Phi_fK,
$$ 
Thus, we get:
\begin{corollary} If, in \ref{locally acyclic}, $g$ is locally acyclic at $y$, then (i) yields an isomorphism
\begin{equation}\label{locally acyclic bis}
\sigma_{x,z}^*R\Phi_{gf}K \iso R\Phi_{g_{(y,z)}}\sigma_{x,y}^*R\Phi_fK.
\end{equation}
\end{corollary}

\begin{remark}\label{locally acyclic ter} In particular, in view of (\ref{psi-good-examples} (b)), if $(f,K)$ is locally acyclic at $x$ and $g$ locally acyclic at $y$, then $(gf,K)$ is locally acyclic at $x$ : this is (\cite{SGA41/2}, Th. finitude, Appendice, 2.7). 
\end{remark}

\begin{remark}  Let $f : X \to Y$, $g : Y \to Z$ be as \ref{composition}, $h = gf$, and let $\Lambda$ be as in \ref{cd of psi}. Let $K \in D^+(Y,\Lambda)$. Consider the corresponding map (\ref{base change morphism})
\begin{equation}\label{inverse image}
\fl{f}^*R\Psi_gK \to R\Psi_h(f^*K),
\end{equation}
where $\fl{f} = X \ori_{Z} Z \to Y \ori_Z Z$ is the functoriality map. Assume that $f$ is locally of finite type, and locally acyclic, i.e., $(f,\Lambda)$ is locally acyclic (cf. (\ref{psi-good-examples}, (b)). Then (\ref{inverse image}) is an isomorphism (cf. (\cite{Laumon:1983ad}, 3.2.3), where this result is stated without proof). 

Indeed, by (\ref{section iso}), taking a geometric point $x \to X$, and replacing $X$, $Y$, $Z$ by their strict localizations at $x$ and its images $y$, $z$ in $Y$ and $Z$ respectively, and $f$, $g$, $h$ by the corresponding localized morphisms, we are reduced to showing that the canonical map
\begin{equation}\label{inverse image1}
Rg_{(y,z)*}K \to Rh_{(x,z)*}(f_{(x,y)}^*K) 
\end{equation}
is an isomorphism, where $K$ denotes $K|Y_{(y)}$ by abuse.  As $f$ is locally of finite type, by \ref{cd of localized map}, $f_{(x,y)*}$ is of finite cohomological dimension. As $f$ is locally acyclic, the formation of $Rf_{(x,y)*}$ commutes with finite base change (\ref{psi-good-examples}, (b)). Therefore, by (\cite{SGA41/2}, Th. finitude, Appendice, 1.2 (a)) the projection formula map
\begin{equation}\label{inverse image2}
Rf_{(x,y)*}\Lambda \otimes^L K \to Rf_{(x,y)*}f_{(x,y)}^*K 
\end{equation}
is an isomorphism\footnote{For $K \in D^b(Y,\Lambda)$. For $K \in D^+(Y,\Lambda)$, see A.5 (1).}. By (\ref{psi-good-examples}, (b)), we have $\Lambda \iso Rf_{(x,y)*}\Lambda$. As  $Rh_{(x,z)*}(f_{(x,y)}^*K) = Rg_{(y,z)*}Rf_{(x,y)*}f_{(x,y)}^*K$, (\ref{inverse image2}) implies that (\ref{inverse image1}) is an isomorphism, which finishes the proof\footnote{As W. Zheng observes, instead of using the projection formula of (\cite{SGA41/2}, Th. finitude, Appendice, 1.2 (a)), one can apply Artin's results in (\cite{SGA4}, XV). This argument does not use the assumption that $f$ is locally of finite type. However, as Artin works with abelian torsion sheaves rather than with sheaves of $\Lambda$-modules, to deduce from (\cite{SGA4}, XV 1.17) that $K \to Rf_{(x,y)*}f_{(x,y)}^*K$ is an isomorphism requires some technical preliminary reductions.}.
\end{remark}

\section{External tensor products}\label{external tensor products}

\subsection{}\label{the kunneth map}
Fix a noetherian scheme $S$, and a noetherian $\mathbf{Z}/n\mathbf{Z}$-algebra $\Lambda$, with $n$ invertible on $S$. For $i = 1, \, 2$, let $f_i : X_i \to Y_i$ be a morphism of finite type between noetherian $S$-schemes, and let $ f = f_1 \times_S f_2 : X = X_1 \times_S X_2 \to Y = Y_1 \times_S Y_2$. Let $K_i \in D_{tf}(X_i,\Lambda)$, where the subscript \textit{tf} means \textit{finite tor-dimension}, and 
$$
K = K_1 \boxtimes^L_S K_2 := \mathrm{pr}_1^*K_1 \otimes^L \mathrm{pr}_2^*K_2 \in D_{tf}(X,\Lambda)
$$
their external tensor product ($\otimes^L$ being taken over $\Lambda$). Consider the oriented topoi $X_i \ori_{Y_i} Y_i$, $X \ori_Y Y$. The morphisms $\mathrm{pr}_i : X \to X_i$, $Y \to Y_i$ define morphisms $\fl {\mathrm{pr}_i} : X \ori_Y Y \to X_i \ori_{Y_i} Y_i$. Let 
\begin{equation}\label{external product}
R\Psi_{f_1}K_1 \boxtimes^L_S R\Psi_{f_2}K_2 := \fl {\mathrm{pr}_1}^*R\Psi_{f_1}K_1  \otimes^L \fl {\mathrm{pr}_2}^*R\Psi_{f_2}K_2 \in D_{tf}(X \ori_Y Y, \Lambda),
\end{equation}
(note that $R\Psi_{f_i}K_i$ is in $D_{tf}(X_i \ori_{Y_i} Y_i,\Lambda)$ by (\cite{Orgogozo06}, Proposition 3.1)).
We define a natural map
\begin{equation}\label{kunneth}
c : R\Psi_{f_1}K_1 \boxtimes^L_S R\Psi_{f_2}K_2  \to R\Psi_f K 
\end{equation}
as follows. Consider the diagram with commutative squares
\begin{equation}\label{kunneth diagram}
\xymatrix{X_1 \ar[d]^{\Psi_{f_1}} & X \ar[l]_{\mathrm{pr_1}} \ar[d]^{\Psi_f} \ar[r]^{\mathrm{pr_2}} &X_2 \ar[d]^{\Psi_{f_2}} \\
X_1 \ori_{Y_1} Y_1 & X \ori_Y Y \ar[l]_{\fl{\mathrm{pr}_1}} \ar[r]^{\fl{\mathrm{pr}_2}} & X_2 \ori_{Y_2} Y_2}. 
\end{equation}
It produces base change maps (\ref{base change morphism})
\begin{equation}\label{base change morphisms}
\fl{\mathrm{pr}_i}^*R\Psi_{f_i}K_i \to R\Psi_f(\mathrm{pr}_i^*K_i),
\end{equation}
hence a tensor product map
$$
\fl{\mathrm{pr}_1}^*R\Psi_{f_1}K_1 \otimes^L \fl{\mathrm{pr}_2}^*R\Psi_{f_2}K_2 \to R\Psi_f(\mathrm{pr}_1^*K_1) \otimes^L R\Psi_f(\mathrm{pr}_2^*K_2),
$$
which, composed with the canonical map
$$
R\Psi_f(\mathrm{pr}_1^*K_1) \otimes^L R\Psi_f(\mathrm{pr}_2^*K_2) \to R\Psi_f K,
$$
yields (\ref{kunneth}). 

\begin{remark}\label{kunneth and base change} In the situation of \ref{the kunneth map}, take $X_2 = Y_2$, $f_2 = \mathrm{Id}_{Y_2}$, so that $f : X \to Y$ is base changed from $f_1 : X_1 \to Y_1$ by $\mathrm{pr}_1 : Y = Y_1 \times_S Y_2 \to Y_1$. The base change map (\ref{base change morphisms}) for $i = 1$
\begin{equation*}
\fl{\mathrm{pr}}_1^*R\Psi_{f_1}K_1 \to R\Psi_f(\mathrm{pr}_1^*K_1)
\end{equation*}
is a particular case of the Künneth map. Indeed, take $K_2 = \Lambda_{Y_2}$. then, by (\cite{Illusie12}, 4.7) applied to $\Psi = \Psi_{\mathrm{Id}_{Y2}} : Y_2 \to Y_2 \ori_{Y_2} Y_2$, $R\Psi_\Lambda = \Lambda$, hence $\fl{\mathrm{pr}}_2^*R\Psi_{f_2}\Lambda = \Lambda$, the base change map $\fl{\mathrm{pr}}_2^*R\Psi_{f_2}\Lambda \to R\Psi_f(\mathrm{pr}_2^*\Lambda)$ is the adjunction map $\Lambda \to R\Psi_f\Lambda$, and the composition given above is just (\ref{base change morphisms}) for $i = 1$.

In particular, one cannot expect (\ref{kunneth}) to be an isomorphism in general (cf. (\ref{psi-good-examples}, (d))).
\end{remark}

\begin{theorem}\label{kunneth theorem} With the notation of \ref{the kunneth map}, assume that, for $i = 1, 2$, $Y_i$ is of finite type over $S$, $K_i$ is in $D_{ctf}(X_i,\Lambda)$, and $(f_i,K_i)$ is $\Psi$-good (\ref{psi-good}). Then $c$ (\ref{kunneth}) is an isomorphism.
\end{theorem}
As mentioned in the introduction, a more general statement, with a much simpler proof, is given in A.3.

\medskip
\begin{proof} We proceed in several steps.

\textit{Step 1}. \textit{We may assume $Y_1 = Y_2 = S$.}

Consider the commutative diagram with cartesian squares 
\begin{equation}
\xymatrix{X \ar[r] \ar[d] & Y \times_{Y_2} X_2 \ar[r]^{p_2} \ar[d] & X_2 \ar[d]_{f_2}\\
X_1 \times_{Y_1} Y \ar[r] \ar[d]_{p_1} & Y \ar[r] \ar[d] &Y_2 \ar[d]\\
X_1 \ar[r]^{f_1} & Y_1 \ar[r] & S }. 
\end{equation}
By the assumptions on $(f_i,K_i)$ the base change maps
$$
\fl{p_i}^*R\Psi_{f_i}K_i \to R\Psi_{X_i \times_{Y_i} Y/Y}(p_i^*K_i),
$$
are isomorphisms, hence we have an identification
$$
R\Psi_{f_1}K_1 \boxtimes^L_S R\Psi_{f_2}K_2  = R\Psi_{(X_1 \times_{Y_1} Y)/Y}(p_1^*K_1) \boxtimes^L_Y R\Psi_{(Y \times_{Y_2} X_2)/Y}(p_i^*K_2),
$$
which reduces the proof of \ref{kunneth theorem} to the case where $Y_1 = Y_2 = S$.

\textit{Step 2}. \textit{(\ref{kunneth}) is an isomorphism if we assume moreover that $Y_1 = Y_2 = S$ and $(f,K)$ is $\Psi$-good.} 

We check that (\ref{kunneth}) is an isomorphism on the slices (\ref{slices}). Let $t \to s$ be a specialization of geometric points of $S$. With the notations of (\ref{slices}), we have morphisms $(i_{\alpha})_s : (X_{\alpha})_s \to (X_{\alpha})_{(s)}$,  $i_s : X_s \to X_{(s)}$, $(j_{\alpha})_{(t)} : (X_{\alpha})_{(t)} \to (X_{\alpha})_{(s)}$, $j_{(t)} : X_{(t)} \to X_{(s)}$, ($\alpha = 1, 2)$, where the subscript $(s)$ (resp. $(t)$) means base change by the strict localization $S_{(s)} \to S$ (resp. $S_{(t)} \to S$). We need to show that the morphism
\begin{equation}\label{sliced kunneth}
(i_1)_s^*R(j_1)_{(t)*}K_1 \boxtimes^L (i_2)_s^*R(j_2)_{(t)*}K_2 \to i_s^*Rj_{(t)*}K
\end{equation}
induced by (\ref{kunneth}) is an isomorphism, where $K_{\alpha}$, $K$ still denote the inverse images of $K_{\alpha}$, $K$ on $(X_{\alpha})_{(t)}$, $X_{(t)}$. Choose a strictly local trait $S'$ with closed point $s'$ and generic geometric point $t'$, and a morphism $S' \to S$ sending $s'$ to $s$ together with a morphism $t' \to t$, compatible with the specialization map $t \to s$.  As $(f_1,K_1)$, $(f_2,K_2)$, $(f,K)$ are $\Psi$-good, the morphism deduced from (2.2.2) by base change by $S' \to S$ is the morphism
\begin{equation}\label{sliced kunneth'}
(i_1)_{s'}^*R(j_1)_{(t')*}K'_1 \boxtimes^L (i_2)_{s'}^*R(j_2)_{(t')*}K'_2 \to i_s^*Rj_{(t')*}K' 
\end{equation}
similar to (\ref{sliced kunneth}), with $X \to S$ replaced by the base changed $X' \to S'$, and $K'_{\alpha}$ (resp. $K'$) induced by $K_{\alpha}$ (resp. $K$). Therefore it suffices to check that (\ref{sliced kunneth'}) is an isomorphism, in other words, we may assume that $S$ is a strictly local trait, with closed point $s$ and generic geometric point $t$. By the comparison recalled in (\ref{classical psi eta}), (\ref{classical psi}), the conclusion then follows from Gabber's theorem on the compatibility of classical nearby cycles with external products (\cite{Illusie94}, 4.7).

\textit{Step 3}. \textit{End of proof : (\ref{kunneth}) is an isomorphism if $Y_1 = Y_2 = S$.}

Let $a : T_{\bullet} \to S$ be a hypercovering for the h-topology which is admissible for $(f : X = X_1 \times_S X_2 \to S, K = K_1 \boxtimes^L_S K_2)$ (\ref{cohomological descent}). Denote by $(K_i)_{\bullet}$ (resp. $K_{\bullet}$) the inverse image of $K_{i}$ (resp. $K$) on $(X_i)_{\bullet} = X_{i} \times_S T_{\bullet}$ (resp. $X_{\bullet} = X \times_S T_{\bullet}$). We denote by
$$
\fl{a} : (X_i)_{\bullet} \ori_{T_{\bullet}} T_{\bullet} \to X_i \ori_S S
$$
($i = 1,2)$ and
$$
\fl{a} : X _{\bullet} \ori_{T_{\bullet}} T_{\bullet} \to X \ori_S S
$$
the morphisms of topoi induced by $a$ (\cite{Orgogozo12b}, 2.2.2), where by abuse the left hand sides denote the total topoi defined by the simplicial topoi. We have commutative squares
\begin{equation}\label{a pr}
\xymatrix{X_{\bullet} \ori_{T_{\bullet}} T_{\bullet} \ar[d]_{\fl{\mathrm{pr}_i}} \ar[r]^{\fl{a}} & X \ori_S S \ar[d]_{\fl{\mathrm{pr}_i}} \\
(X_i)_{\bullet} \ori_{T_{\bullet}} T_{\bullet} \ar[r]^{\fl{a}} & X_i \ori_S S}.
\end{equation}
As $(f_i,K_i)$ is $\Psi$-good, the base change maps
\begin{equation}\label{a* iso1}
\fl{a}^*R\Psi_{f_i}K_i \to R\Psi_{(f_i)_{\bullet}}(K_i)_{\bullet}
\end{equation}
are isomorphisms. Hence the same holds for their external tensor product, which, thanks to (\ref{a pr}), can be re-written
\begin{equation}\label{a* iso2}
\fl{a}^*(R\Psi_{f_1}K_1 \boxtimes^L_S R\Psi_{f_2}K_2) \to R\Psi_{(f_1)_{\bullet}}(K_1)_{\bullet} \boxtimes^L_S R\Psi_{(f_2)_{\bullet}}(K_2)_{\bullet}.
\end{equation}
Consider the morphism of $D(X \ori_S S,\Lambda)$
\begin{equation}\label{basic map}
c_{\bullet} : R\Psi_{(f_1)_{\bullet}}(K_1)_{\bullet} \boxtimes^L_S R\Psi_{(f_2)_{\bullet}}(K_2)_{\bullet} \to R\Psi_{f_{\bullet}}K_{\bullet} 
\end{equation}
defined similarly to (\ref{kunneth}), using the variant of diagram (\ref{kunneth diagram}) with $X_{i}$ (resp. $X$) replaced by $(X_i)_{\bullet}$ (resp. $X_{\bullet}$). (Note that because of (\ref{a* iso1}) and (\ref{a* iso2}) the left hand side is in $D_{ctf}$, while \textit{a priori} the right hand side is only in $D^+$.) We will show that (\ref{basic map}) is an isomorphism. As the family of restriction maps $D^+(X_{\bullet} \ori_{S_{\bullet}} S_{\bullet},\Lambda) \to D^+(X_n \ori_{S_n} S_n,\Lambda)$ is conservative, it is enough to check that, for each $n \in \mathbf{N}$, 
$$
c_{n} : R\Psi_{(f_1)_{n}}(K_1)_{n} \boxtimes^L_S R\Psi_{(f_2)_{n}}(K_2)_{n} \to R\Psi_{f_{n}}K_{n} 
$$
is an isomorphism. This is indeed the case by Step 2, as $(f_{\alpha}, K_{\alpha})$ is $\Psi$-good and $a$ is admissible for $(f,K)$, hence $(f_n,K_n)$ is $\Psi$-good. Now, consider the commutative square deduced from the functoriality of the construction of the Künneth maps,
\begin{equation}\label{basic diagram}
\xymatrix{R\Psi_{f_1}K_1 \boxtimes^L_S R\Psi_{f_2}K_2 \ar[d] \ar[r]^c & R\Psi_fK \ar[d] \\
R\fl a_*(R\Psi_{(f_1)_{\bullet}}(K_1)_{\bullet} \boxtimes^L_S R\Psi_{(f_2)_{\bullet}}(K_2)_{\bullet}) \ar[r] & R\fl a_*R\Psi_{f_{\bullet}}K_{\bullet}},
\end{equation}
where the bottom horizontal map is $R\fl a_*c_{\bullet}$,  and the vertical maps are the canonical maps deduced from the isomorphisms $a^*K_{i} = (K_{i})_{\bullet}$, $a^*K = K_{\bullet}$. The bottom horizontal map is an isomorphism as we have just shown that $c_{\bullet}$ is an isomorphism. As (\ref{a* iso2}) is an isomorphism, the left vertical map is an isomorphism by oriented descent (\ref{oriented descent iso}). The right vertical arrow is the isomorphism (\ref{recovering psi}). Therefore $c$ is an isomorphism, which concludes the proof.
\end{proof}
\begin{corollary}\label{psi-goodness of (f,K)}
Under the assumptions of \ref{kunneth theorem}, $(f,K)$ is $\Psi$-good. 
\end{corollary}
\begin{proof} Indeed, as $(f_i,K_i)$ is $\Psi$-good, the formation of $R\Psi_{f_1}K_1 \boxtimes^L_S R\Psi_{f_2}K_2$ commutes with arbitrary base change. More explicitly, for $g : S' \to S$, let $Y'_i := Y_i \times_S S'$, $X'_i := X_i \times_S S'$, $f'_i = f_i \times_S S' : X'_i \to Y'_i$, $X' := X \times_S S'$, $K'_i := K_i | X_i$. We then have a commutative diagram
\begin{equation*}
\xymatrix{X'_1 \ori_{Y'_1} Y'_1 \ar[d]_{\fl{g}} & X' \ori_{Y'} Y' \ar[l]_{\fl{\mathrm{pr}_1}} \ar[r]^{\fl{\mathrm{pr}_2}} \ar[d]_{\fl{g}} & X'_2 \ori_{Y'_2} Y'_2  \ar[d]_{\fl{g}} \\
X_1 \ori_{Y_1}Y_1 & X \ori_Y Y \ar[l]_{\fl{\mathrm{pr}_1}} \ar[r]^{\fl{\mathrm{pr}_2}} & X_2 \ori_{Y_2} Y_2}
\end{equation*}
Write $\Psi_i$ (resp. $\Psi'_i$) for short for $R\Psi_{f_i}K_i$ (resp. $R\Psi_{f'_i}K'_i$). The above diagram gives a (trivial) isomorphism
\begin{equation}\label{trivial iso}
\fl{g}^*(\fl{\mathrm{pr}_1}^*\Psi_1 \otimes^L \fl{\mathrm{pr}_2}^*\Psi_2) \iso \fl{\mathrm{pr}_1}^*(\fl{g}^*\Psi_1) \otimes^L \fl{\mathrm{pr}_2}^*(\fl{g}^*\Psi_2).
\end{equation}
By $\Psi$-goodness of $(f_i,K_i)$, the base change maps $\fl{g}^*\Psi_i \to \Psi'_i$ are isomorphisms. By composing them with (\ref{trivial iso}) we get an isomorphism
\begin{equation}\label{nontrivial iso}
\fl{g}^*(\Psi_1 \boxtimes^L_S \Psi_2) \iso \Psi'_1 \boxtimes^L_{S'} \Psi'_2,
\end{equation}
with the notation of (\ref{external product}). This map fits into a commutative diagram
\begin{equation*}
\xymatrix{\fl{g}^*(\Psi_1 \boxtimes^L_S \Psi_2) \ar[d]_{(\ref{nontrivial iso})} \ar[r]^{\fl{g}^*c} &\fl{g}^*\Psi \ar[d] \\
\Psi'_1 \boxtimes^L_{S'} \Psi'_2 \ar[r]^{c} & \Psi'} 
\end{equation*}
where $\Psi := R\Psi_fK$, $\Psi' := R\Psi_{f'}K'$, and the right vertical map is the base change map, which is therefore an isomorphism, as $c$ is an isomorphism.
\end{proof}

In particular, combining with \ref{kunneth theorem}, we get:

\begin{corollary}\label{kunneth acyclic} In the situation of \ref{kunneth theorem}, assume that, for $i = 1, 2$, $(f_i,K_i)$ is universally locally acyclic (\ref{psi-good-examples} (b)). Then $(f,K)$ is universally locally acyclic.
\end{corollary}

One can make \ref{kunneth theorem} explicit on the stalks:

\begin{corollary}\label{punctual kunneth}
Under the assumptions of \ref{kunneth theorem}, let $s$ be a geometric point of $S$. For $i = 1, 2$, let $y_i$ be a geometric point of $Y_i$ above $s$, $x_i$ a geometric point of $X_i$ above $y_i$, $x = (x_1,x_2)$ and $y = (y_1,y_2)$ the corresponding geometric points of $X = X_1 \times_S X_2$ and $Y = Y_1 \times_X Y_2$. Then, with the notation of (\ref{canonical section}), the isomorphism $c$ (\ref{kunneth}) and the commutative diagram 
\begin{equation}\label{punctual kunneth diagram}
\xymatrix{X_1 \ori_{Y_1} Y_1 & X \ori_Y Y \ar[l]_{\fl{pr_1}} \ar[r]^{\fl{pr_2}} & X_2 \ori_{Y_2} Y_2 \\
(Y_1)_{(y_1)} \ar[u]_{\sigma_{x_1,y_1}} & Y_{(y)} \ar[u]_{\sigma_{x,y}} \ar[l]_{pr_1} \ar[r]^{pr_2} & (Y_2)_{(y_2)} \ar[u]_{\sigma_{x_2,y_2}}},
\end{equation}
induce an isomorphism in $D_{ctf}(Y_{(y)},\Lambda)$ :
\begin{equation}\label{punctual kunneth1}
\sigma_{x_1,y_1}^*R\Psi_{f_1}K_1 \boxtimes^L \sigma_{x_2,y_2}^*R\Psi_{f_2}K_2 \iso \sigma_{x,y}^*(R\Psi_fK), 
\end{equation}
where $\boxtimes^L$ in the left hand side means $pr_1^* \otimes^L pr_2^*$, with $\mathrm{pr}_i$ as in (\ref{punctual kunneth diagram}).
\end{corollary}

\begin{corollary} In the situation of \ref{the kunneth map}, assume that $Y_i$ ($i = 1,2$) is regular of dimension $\le 1$. Then the map $c$ (\ref{kunneth}) is an isomorphism, and $(f,K)$ is $\Psi$-good. 
\end{corollary}
\begin{proof}
Indeed, for $(i = 1,2)$, $(f_i,K_i)$ is $\Psi$-good (\ref{psi-good-examples} (a)).
\end{proof}

\subsection{}\label{application to phi}
Let's come back to the situation of \ref{the kunneth map}. In order to analyze the behavior of $R\Phi$ under the Künneth map (\ref{kunneth}), we need to use the refined objects $R\underline{\Psi}_{f_i}K_i$ (\ref{filtered Rpsi}) of $DF^{[0,1]}_{tf}(X_i \ori_{Y_i}Y_i,\Lambda)$. Their external tensor product is a 2-step filtered object :
\begin{equation}\label{filtered boxtimes}
R\underline{\Psi}_{f_1}K_1 \boxtimes^L_S R\underline{\Psi}_{f_2}K_2 \in DF^{[0,2]}_{tf}(X \ori_Y Y,\Lambda),
\end{equation}
with associated graded
\begin{equation}\label{gr boxtimes}
\mathrm{gr}(R\underline{\Psi}_{f_1}K_1 \boxtimes^L_S R\underline{\Psi}_{f_2}K_2 ) = \mathrm{gr} R\underline{\Psi}_{f_1}K_1 \boxtimes^L_S \mathrm{gr}R\underline{\Psi}_{f_2}K_2,
\end{equation}
i. e.
\begin{align*}
\mathrm{gr}^0(R\underline{\Psi}_{f_1}K_1 \boxtimes^L_S R\underline{\Psi}_{f_2}K_2 ) &=& R\Phi_{f_1}K_1 \boxtimes^L_S R\Phi_{f_2}K_2 \\
\mathrm{gr}^1(R\underline{\Psi}_{f_1}K_1 \boxtimes^L_S R\underline{\Psi}_{f_2}K_2 ) &=& (R\Phi_{f_1}K_1 \boxtimes^L_S p_1^*K_2) \oplus (p_1^*K_1 \boxtimes^L_S R\Phi_{f_2}K_2) \\
\mathrm{gr}^2(R\underline{\Psi}_{f_1}K_1 \boxtimes^L_S R\underline{\Psi}_{f_2}K_2 ) &=& p_1^*K_1 \boxtimes^L_S p_2^*K_2.
\end{align*}
It follows that under the assumptions of \ref{kunneth theorem}, the isomorphism $c$ (\ref{kunneth}) defines a filtered object $(R\Psi_fK,F_2) \in DF^{[0,2]}(X \ori_Y Y,\Lambda)$, with associated graded given by (\ref{gr boxtimes}), where $F_2$ refines the filtration $F$ of  (\ref{filtered Rpsi}) :
\begin{align*}
F^2_2 = F^1 = p_1^*K_1 \boxtimes^L_S p_2^*K_2 \\ 
R\Phi_fK = \mathrm{gr}^0_F = F^0_2/F^2_2.
\end{align*}
In particular, we have a distinguished triangle
\begin{equation}\label{Rphi triangle}
(R\Phi_{f_1}K_1 \boxtimes^L_S p_1^*K_2) \oplus (p_1^*K_1 \boxtimes^L_S R\Phi_{f_2}K_2) \to R\Phi_fK \to R\Phi_{f_1}K_1 \boxtimes^L_S R\Phi_{f_2}K_2 \to .
\end{equation}
Thus, in the situation of \ref{punctual kunneth}, as  by (\ref{stalk at y vanishing cycles triangle}) $\sigma_{x_i,y_i}^*p_1^*K_i$  is the constant complex on $Y_{(y_i)}$ of value $(K_i)_{x_i}$, we get a distinguished triangle
\begin{equation}\label{punctual Rphi triangle}
(\sigma_{x_1,y_1}^*R\Phi_{f_1}K_1 \boxtimes^L (K_2)_{x_2}) \oplus ((K_1)_{x_1} \boxtimes^L \sigma_{x_2,y_2}^*R\Phi_{f_2}K_2) \to \sigma_{x,y}^*R\Phi_fK 
\end{equation}
\begin{equation*}
\to \sigma_{x_1,y_1}^*R\Phi_{f_1}K_1 \boxtimes^L \sigma_{x_2,y_2}^*R\Phi_{f_2}K_2 \to .
\end{equation*}
The following is a generalization of \ref{kunneth acyclic}:
\begin{corollary}\label{isolated1} In the situation of \ref{the kunneth map}, assume that, for $i = 1, 2$, there is an open subset $U_{i}$ of $X_{i}$ with $\Sigma_{i} = X_{i} -U_{i}$ quasi-finite over $Y_{i}$ such that $(f_{i},K_{i})|U_{i}$ is universally locally acyclic (\cite{SGA41/2}, Th. finitude, 2.12). Let $U := U_1 \times_S U_2$, $\Sigma := X - U = (\Sigma_1 \times_S X_2) \cup (X_2 \times_S \Sigma_2)$. Then $(f_{i},K_{i})$ ($i = 1, 2$) and $(f,K)$ are $\Psi$-good,  and $R\Phi_f(K)$ is concentrated on $\Sigma \ori_Y Y$.
\end{corollary}

\begin{proof} As $(f_i,K_i)$ is $\Psi$-good (\ref{psi-good-examples} (c)), $(f,K)$ is also $\Psi$-good (\ref{psi-goodness of (f,K)}).  As $(f_i,K_i)|U_i$ is universally locally acyclic, $R\Phi_{(f_i|U_i)}(K_i|U_i) = 0$ (\ref{psi-good-examples} (b)). By (\ref{Rphi triangle}) we therefore have $R\Phi_fK | U \ori_Y Y = 0$, which implies the last assertion, as the complement of the closed subtopos $\Sigma \ori_Y Y$ of $X \ori_Y Y$ is the open subtopos $U \ori_{Y} Y$.
\end{proof}

\section{Interlude: additive convolution}\label{interlude}
\subsection{}\label{notation} We fix a perfect field $k$ of characteristic exponent $p$, an algebraic closure $\overline{k}$ of $k$, and a \textit{finite} ring $\Lambda$ annihilated by an integer invertible in $k$. We denote by
$\pi_1^t(\mathbf{G}_{m,k},\{\overline{1}\})$ the tame quotient of the fundamental group of $\mathbf{G}_{m,k}$, which is an extension
\begin{equation}\label{pi1tame}
1 \to I_t \to \pi_1^t(\mathbf{G}_{m,k},\{\overline{1}\}) \to \mathrm{Gal}(\overline{k}/k) \to 1,
\end{equation}
where 
\begin{equation}\label{It}
I_t \, \, (\iso \widehat{\mathbf{Z}}'(1) = \varinjlim_{(n,p) = 1} \mu_n(\overline{k})) 
\end{equation}
is the tame quotient of  the geometric fundamental group $I = \pi_1(\mathbf{G}_{m,\overline{k}},\{\overline{1}\})$ of $\mathbf{G}_{m,k} = \mathbf{A}^1_k - \{0\}$ (here $\overline{1}$ means the unit of $\mathbf{G}_{m,\overline{k}}$).  If $\overline{\eta}_{0}$ (resp. $\overline{\eta}_{\infty}$) is a geometric point above the generic point $\eta_0$ (resp. $\eta_{\infty}$) of the henselization of $\mathbf{A}^1_k$ at $\{0\}$ (resp. $\{\infty\}$), we can consider the tame quotient $(I_t)_0$ (resp. $(I_t)_{\infty}$) of the inertia subgroup $I_0 \subset \mathrm{Gal}(\overline{\eta}_0/\eta_0)$ (resp. $(I_t)_{\infty} \subset \mathrm{Gal}(\overline{\eta}_{\infty}/\eta_{\infty})$, which maps isomorphically to $I_t$ (once a path is chosen from $\overline{\eta}_0$ (resp. $\overline{\eta}_{\infty}$) to $\{\overline{1}\}$). 

Recall that a lisse sheaf $L$ of $\Lambda$-modules on $\mathbf{G}_{m,k}$ is tamely ramified at $0$ (resp. $\infty$) if the action of $I_0$ (resp. $I_{\infty}$) on $L_{\overline{\eta}_0}$ (resp. $L_{\overline{\eta}_{\infty}}$) factors through $(I_t)_0$ (resp. $(I_t)_{\infty}$) (or, equivalently, through a finite quotient of it, as $\Lambda$ has been assumed finite). The sheaf $L$ is tamely ramified at $0$ and $\infty$ if and only if the action of $I$ on $L_{\{\overline{1}\}}$ factors through $I_t$. 

More generally, given $L \in D^b_c(\mathbf{A}^1_k,\Lambda)$, we say that $L$ is \textit{tamely ramified} at $0$ (resp. $\infty$) if, for all $q$, $(\mathcal{H}^qL)_{\eta_0}$ (resp. $\mathcal{H}^qL)_{\eta_{\infty}}$) is tame (i.e., the action of $I_0$ (resp. $I_{\infty}$) on $\mathcal{H}^qL_{\overline{\eta}_0}$ (resp. $\mathcal{H}^qL_{\overline{\eta}_{\infty}}$) factors through $(I_t)_0$ (resp. $(I_t)_{\infty}$)). 

We denote by 
\begin{equation}\label{global sum map}
a : \mathbf{A}^1_k \times \mathbf{A}^1_k \to \mathbf{A}^1_k
\end{equation}
the sum map.  

\medskip
\textbf{A. Global additive convolution}

\medskip
Recall the following definition, due to Deligne.
\begin{definition}\label{global convolution} The (additive) convolution functor
\begin{equation*}
*^L : D_{ctf}(\mathbf{A}^1_k, \Lambda) \times D_{ctf}(\mathbf{A}^1_k,\Lambda) \to D_{ctf}(\mathbf{A}^1_k,\Lambda)
\end{equation*}
is defined by
\begin{equation*}
K_1 *^L K_2 := Ra_*(K_1 \boxtimes^L_{\Lambda} K_2).
\end{equation*}
\end{definition}
\begin{remark}\label{convolution !} (a) There is a variant $*^L_!$ of $*^L$, with $Ra_*$ replaced by $Ra_!$. These two constructions are exchanged by the dualizing functor on $\mathbf{A}^1_k$, $D := R\mathcal{H}om(-,\Lambda(1)[2])$:
\begin{equation}\label{D*}
D(K_1 *_!^L K_2) = DK_1 *^L DK_2.
\end{equation}
The $*^L_!$ definition has the advantage of corresponding to the usual convolution of trace functions, when $k$ is a finite field $\mathbf{F}_q$, i. e. if, for $K \in D_{ctf}(\mathbf{A}^1_k,\Lambda)$, and $n \ge 1$, $f_K$ denotes the function on $\mathbf{A}^1_k(\mathbf{F}_{q^n})$ defined by $f_K(x) = \mathrm{Tr}(F^n,K_x)$, then
\begin{equation*}
f_{K_1 *^L_! K_2} = f_{K_1} * f_{K_2}. 
\end{equation*}
(b) Assume $p>1$. If $\psi : \mathbf{F}_p \to \Lambda^*$ is a non trivial additive character, we can consider the Fourier transform $\mathcal{F} : \mathcal{F}_{\psi} : D_{ctf}(\mathbf{A}^1_k,\Lambda) \to D_{ctf}(\mathbf{A}^1_k,\Lambda)$, defined by the formula (1.2.1.1) of \cite{Laumon87}, with $\overline{\mathbf{Q}}_{\ell}$ replaced by $\Lambda$. By (\cite{Laumon87}, 1.2.2.7), we have
\begin{equation}\label{Fourier !}
\mathcal{F}(K_1 *^L_! K_2) \iso (\mathcal{F}K_1 \otimes^L \mathcal{F}K_2)[-1].
\end{equation}
By (\cite{Laumon87}, 1.3.2.2), $D\mathcal{F}_{\psi} = (\mathcal{F}_{\psi^{-1}}D)(1)$, so that applying (\ref{D*}), we get
\begin{equation}\label{Fourier *}
\mathcal{F}(K_1 *^L K_2) \iso (\mathcal{F}K_1 \otimes^{!L} \mathcal{F}K_2)(1)[1],
\end{equation}
where $L_1 \otimes^{!L} L_2 := D(DL_1 \otimes^L DL_2)$. As $L_1 \otimes^{!L} L_2 = \Delta^!(L_1 \boxtimes^L L_2)$ where $\Delta : \mathbf{A}^1_k \to \mathbf{A}^2_k$ is the diagonal, we have $L_1 \otimes^{!L} L_2 = (L_1 \otimes^L L_2)[-2](-1)$, so, from (\ref{Fourier *}) we get
\begin{equation}\label{Fourier *bis}
\mathcal{F}(K_1 *^L K_2) \iso (\mathcal{F}K_1 \otimes^{L} \mathcal{F}K_2)[-1],
\end{equation}
which is actually the same formula as (\ref{Fourier !}). 
\end{remark}
The following results are standard (cf.  \cite{Deligne80}, and \cite{Laumon83}, \cite{Laumon87} for the case of $*!$ and $\overline{\mathbf{Q}}_{\ell}$-coefficients):

\begin{proposition}\label{global * formulas} Let $e : \mathbf{A}^1_k \to \SP \, k$ be the projection, and $i : \{0\} \to \mathbf{A}^1_k$, $j : \mathbf{A}^1_k - \{0\} \to \mathbf{A}^1_k$ be the inclusions. Let $K_i \in D_{ctf}(\mathbf{A}^1_k,\Lambda)$ ($i = 1,2$).

(1) For $K_1$ concentrated at $\{0\}$, i. e. $K_1 = i_*i^*K_1$, we have a canonical isomorphism
\begin{equation*}
K_1 *^L K_2 \iso e^*(i^*K_1) \otimes^L K_2.
\end{equation*}

(2) Assume that $K_2$ is geometrically constant, i. e. is of the form $e^*L_2$ for $L_2 \in D_{ctf}(\SP \, k,\Lambda)$. Then we have a canonical isomorphism
\begin{equation}\label{global * vanishing1}
K_1 *^L K_2 \iso e^*Re_*K_1 \otimes^L K_2,
\end{equation}
and
\begin{equation}\label{global vanishing phi}
R\Phi_a(K_1 \boxtimes^L K_2) = 0.
\end{equation}
If in addition $K_1$ is concentrated in degree $0$, $(K_1)_0 = 0$, lisse on $\mathbf{A}^1_k - \{0\}$, and tamely ramified at $\{0\}$ and $\{\infty\}$, 
then 
\begin{equation}\label{global * vanishing2}
K_1 *^L K_2 = 0.
\end{equation}
Similar statement with $K_1$ and $K_2$ interchanged.

(3) Assume that, for $i = 1, 2$, $K_i$ is concentrated in degree $0$, lisse on $\mathbf{A}^1_k - \{0\}$, tamely ramified at $\infty$, and $(K_i)_0 = 0$.  Then
\begin{equation}\label{concentration in degree 1}
\mathcal{H}^q(K_1 *^L K_2) = 0
\end{equation}
for $q \ne 1$, and $\mathcal{H}^1(K_1 *^L K_2)$ is a constructible sheaf of finitely generated and projective $\Lambda$-modules, lisse on $\mathbf{A}^1_k - \{0\}$.  Its generic rank $r$ is 
\begin{equation}\label{generic rank formula}
r = r_1r_2 + r_1s_2 + r_2s_1,
\end{equation}
where $r_i = \mathrm{rk}(K_i)$, $s_i = \mathrm{sw}_0(K_i)$, $\mathrm{sw}_0$ denoting the Swan conductor at $0$.  Its rank at zero is 
\begin{equation}\label{rank at zero}
\mathrm{rk}(\mathcal{H}^1(K_1 *^L K_2)_0) = \mathrm{sw}_0(K_1 \otimes [-1]^*K_2),
\end{equation}
where $[-1]$ is the automorphism of $\mathbf{A}^1_k = \SP \,k[t]$ given by $t \mapsto -t$.

For any homomorphism $\Lambda \to \Lambda'$ of rings satisfying the assumptions of \ref{notation}, if $K'_i := K_i \otimes_{\Lambda} \Lambda'$, the natural map
\begin{equation*}
\mathcal{H}^1(K_1 *^L_{\Lambda} K_2) \otimes_{\Lambda} \Lambda' \to \mathcal{H}^1(K'_1 *^L_{\Lambda'} K'_2)
\end{equation*}
is an isomorphism. 

For any geometric point $t$ of $\mathbf{A}^1_k$, the canonical map
\begin{equation}\label{specialization}
\mathcal{H}^1(K_1 *^L K_2)_t \to H^1(a^{-1}(t),K_1 \boxtimes^L K_2) 
\end{equation}
is an isomorphism. 

(4) Under the assumptions of (3), let $a_{(0)} : \mathbf{A}^2_k \times_{\mathbf{A}^1_k} S \to S$ be the map deduced from $a$ by base change to the henselization $S$ at $\{0\}$ of $\mathbf{A}^1_k$. Then $R\Phi_{a_{(0)}}(K_1 \boxtimes^L K_2) \in D_{ctf}(a^{-1}(0)\ori_S \eta,\Lambda)$ is concentrated at $(0,0)$ and in degree 1, and we have
\begin{equation}\label{global *vanishing3}
R\Psi_{a_{(0)}}(K_1 \boxtimes^L K_2)_{(0,0)}[1] = R\Phi_{a_{(0)}}(K_1 \boxtimes^L K_2)_{(0,0)}[1] = R^0\Phi_{\mathrm{Id}}(\mathcal{H}^1(K_1 *^L K_2))_0 
\end{equation}
where $\mathrm{Id}$ means the identity of $\mathbf{A}^1_{(0)}$, with
\begin{equation}\label{rank formula}
\mathrm{rk}(R^0\Phi_{\mathrm{Id}}(\mathcal{H}^1(K_1 *^L K_2)))_0 = r_1r_2 + r_1s_2 + r_2s_1 - \mathrm{sw}_0(K_1 \otimes [-1]^*K_2).
\end{equation}
\end{proposition}

\begin{proof}
We may (and will) assume $k$ algebraically closed. We write $A$ for $\mathbf{A}^1_k$.

(1) Put $L_1 := i^*K_1$, and let $i_1 := i \times \mathrm{Id} : \{0\} \times A \to A^2$. We have
\begin{equation*}
K_1 \boxtimes^L K_2 = i_{1*}(e^*L_1 \otimes^L K_2),
\end{equation*}
and, as $ai_1 : \{0\} \times A \to A$ is the identity, we get
\begin{equation*}
K_1 *^L K_2 = Ra_*(K_1 \boxtimes^L K_2) = R(ai_1)_*(e^*L_1 \otimes^L K_2) \iso e^*L_1 \otimes^L K_2.
\end{equation*}

(2) As $\mathrm{pr}_2^*e^* = \mathrm{pr}_1^*e^*$, we have
\begin{equation*}
K_1 *^L K_2 = Ra_*\mathrm{pr}_1^*(K_1 \otimes^L e^*L_2).
\end{equation*}
Let 
\begin{equation*}
\varphi : A^2 \iso A^2\end{equation*}
be the isomorphism given by $\varphi(x,y) = (x+y,x)$. We have $\mathrm{pr}_1\varphi = a$, $\mathrm{pr}_2\varphi = \mathrm{pr}_1$, hence
\begin{equation*}
K_1 *^L K_2  = R\mathrm{pr}_{1*}\varphi_* \mathrm{pr}_1^*(K_1 \otimes^L e^*L_2),
\end{equation*}
but the base change  map $\mathrm{pr}_2^* \to \varphi_*\mathrm{pr}_1^*$ is (trivially) an isomorphism, so we get
\begin{equation*}
K_1 *^L K_2 \iso R\mathrm{pr}_{1*}\mathrm{pr}_2^*(K_1 \otimes^L e^*L_2).
\end{equation*}
Finally, by smooth base change and the projection formula, we get
\begin{equation*}
K_1 *^L K_2 \iso e^*Re_*(K_1 \otimes^L e^*L_2) \iso e^*Re_*K_1 \otimes^L K_2.
\end{equation*}
Similarly, we have
\begin{equation*}
R\Phi_a(K_1 \boxtimes^L K_2) \iso R\Phi_{\mathrm{pr}_1}\mathrm{pr}_2^*(K_1 \otimes^L e^*L_2),
\end{equation*}
and for any $M \in D_{ctf}(A,\Lambda)$, $R\Phi_{\mathrm{pr}_1}\mathrm{pr}_2^*M = 0$ by universal local acyclicity for schemes of finite type over a field (\cite{SGA41/2}, Th. finitude, 2.16) (cf. \ref{psi-good-examples} (b)).
For the last assertion, it suffices to show $Re_*K_1 = 0$. If $L_1$ is the largest constant subsheaf of $j^*K_1$, we have $e_*K_1 = e_*(j_!L_1) = 0$. It remains to show $R^1e_*K_1 = 0$. We have 
\begin{equation}\label{Laumon's theorem}
\chi(A,K_1) = \chi(A,j_!j^*K_1) = \chi(A,Rj_*j^*K_1)
\end{equation}
as $\chi(\{0\},i^*Rj_*j^*K_1) = 0$ by Laumon's theorem \cite{Laumon81}\footnote{As the referee points out, in the case of a strictly local trait $(S, i : s \to S, j : \eta \to S)$, Laumon's theorem to the effect that for a sheaf $\mathcal{F}$ on $\eta$, $\chi(s,i^*Rj_*\mathcal{F}) = 0$ is an immediate consequence of the standard formulas giving $i^*R^qj_*\mathcal{F}$ (cf. \cite{SGA41/2}, Dualité, proof of 1.3).}. As $K_1$ is tamely ramified at $\{0\}$ and $\{\infty\}$, we have 
\begin{equation*}
\chi(A - \{0\},j^*K_1) = \chi(A -\{0\},\Lambda)\mathrm{rk}(j^*K_1) = 0,
\end{equation*}
by the (Ogg-Shafarevitch case of the) Grothendieck-Ogg-Shafarevitch formula, which finishes the proof of (\ref{global * vanishing2}).

(3) Put $K = K_1 \boxtimes K_2$. We use the partial compactification of (\cite{Laumon87}, 2.7.1.1 (iii)),
$$
A \times_k A  \overset{(\mathrm{pr}_1,a)}{\to} A \times_k A \overset{\alpha \times \mathrm{Id}}{\to} D \times A,
$$
where $D = \mathbf{P}^1_k = A \cup \{\infty\}$, $\alpha : A \inj D$ is the inclusion, and let $\overline{pr}_2 : D \times_k A \to A$ denote the projection. If $t$ is a geometric point of $A$, we have
$$
(Ra_*K)_t = R\Gamma(D \times \{t\},(R(\alpha \times \mathrm{Id})_*(\mathrm{pr}_1, a)_*K)|D \times \{t\})
$$
by properness of $\overline{\mathrm{pr}}_2$. As $(\mathrm{pr}_1,a)_*K$ is tamely ramified along $\{\infty\} \times_k A$, by (\cite{SGA41/2}, Th. finitude, Appendice, 1.3.3 (i))  we have
$$
(R(\alpha \times \mathrm{Id})_*(\mathrm{pr}_1,a)_*K)|D \times \{t\} \iso R\alpha_{t*}((\mathrm{pr}_1, a)_*K | A \times \{t\}),
$$
where $\alpha_t : A \times \{t\} \inj D \times \{t\}$ is the fiber of $\alpha \times \mathrm{Id}$ at $t$. Finally, $(\mathrm{pr}_1, a)_*K | A \times \{t\} = a_{t*}(K|a^{-1}(t))$, where $a_t : a^{-1}(t) \iso A \times \{t\}$ is induced by $\mathrm{pr}_1$, and we get that the specialization map (\ref{specialization})
$$
(Ra_*K)_t \to R\Gamma(a^{-1}(t),K).
$$
is an isomorphism. This also implies (\ref{concentration in degree 1}), since $a^{-1}(t) \iso A$ is affine, and $H^0(a^{-1}(t),K) = 0$ as $K | a^{-1}(t)$ is lisse outside $(t,0)$ and $(0,t)$ and vanishes at these points.  The compatibility of the formation of $\mathcal{H}^1(K_1 *^L K_2)$ with extension of scalars $\Lambda \to \Lambda'$ follows, as $K_1 *^L K_2$ is in $D_{ctf}(A,\Lambda)$.  Let us show that $\mathcal{H}^1(K_1 *^L K_2)$ is lisse on $A - \{0\}$. Let $E = R(\alpha \times \mathrm{Id})_*(\mathrm{pr}_1,a)_*K$. For a geometric point $t$ of $A - \{0\}$, let $a_{(t)}$ (resp. $\overline{pr}_{2(t)}$) denote the map deduced from $a$ (resp. $\overline{pr}_2$) by henselization at $t$. As $(\mathrm{pr}_1,a)_*K$ is tamely ramified along $\{\infty\} \times_k A$,  by a classical result (cf. (\cite{TSaito1301.4632}, 3.14)) $(\overline{pr}_2,E)$ is universally locally acyclic along $\{\infty\} \times_k A$, hence $R\Phi_{\overline{pr}_{2(t)}}E \in D^b_c(D \times_k \{t\},\Lambda)$ is zero at $(\infty,t)$. On the other hand, at $(0,t)$ (resp. $(t,t)$), $R\Phi_{\overline{pr}_{2(t)}}E$, which is isomorphic to $R\Phi_{a_{(t)}}K$ at $(0,t)$ (resp. $(t,0)$), is also zero. Indeed, on the strict localization of $A^2$ at (a geometric point above) $(0,t)$, $\mathrm{pr}_2^*K_2$ is constant, so we may assume $K_2 = \Lambda$. Then, as above, by the isomorphism $\varphi$ in the proof of (2), $R\Phi_{a_{(t)}}K$ is identified with $R\Phi_{\mathrm{pr}_1}\mathrm{pr}_2^*K_1$ at $(t,0)$, which is zero. The proof at $(t,0)$ is similar. It follows that $(\overline{\mathrm{pr}}_2,E)$ is universally locally acyclic, and as $\overline{pr}_2$ is proper, this implies that the specialization and cospecialization maps at $t$ are isomorphisms (\cite{SGA41/2}, Th. finitude, Appendice, 2.4), in other words, that $\mathcal{H}^1(K_1 *^L K_2)$ is lisse on $\mathbf{G}_{m,k}$. 

Let us prove (\ref{generic rank formula}). By (\ref{specialization}), 
\begin{equation*}
r = \mathrm{rk}(H^1(a^{-1}(t),K))
\end{equation*}
for $t \in A - \{0\}$. As $K|a^{-1}(t)$ is tamely ramified at $\infty$, and of rank $r_1r_2$, and taking into account that $\chi(a^{-1}(t),K) = \chi_c(a^{-1}(t),K)$ by \cite{Laumon81}, the Grothendieck-Ogg-Shafarevich formula gives
\begin{equation*}
\chi(a^{-1}(t),K) = r_1r_2 - (r_1r_2 + \mathrm{sw}_{(0,t)}(K)) - (r_1r_2 + \mathrm{sw}_{(t,0)}(K)).
\end{equation*}
As $K_2$ (resp. $K_1$) is lisse at $(0,t)$ (resp. $(t,0)$), we have $\mathrm{sw}_{(0,t)}(K) = s_1r_2$ (resp. $\mathrm{sw}_{(t,0)}(K) = s_2r_1$), which gives
\begin{equation*}
\chi(a^{-1}(t),K) = -r_1r_2 -s_1r_2 - s_2r_1,
\end{equation*}
hence (\ref{generic rank formula}) by (\ref{concentration in degree 1}). Finally, on $a^{-1}(0)$, $K$ is isomorphic to 
$K_1 \otimes [-1]^*K_2$ on $A$, hence $\mathrm{sw}_0(K |a^{-1}(0)) = \mathrm{sw}_0(K_1 \otimes [-1]^*K_2)$. Using again (\ref{specialization}) (at $t=0$) and the Grothendieck-Ogg-Shafarevich formula, we get (\ref{rank at zero}).

(4) This is similar to (\cite{Laumon87}, 2.7.1.1 (iii)), whose proof works for $\Lambda$ as a coefficients ring.  However, as here we consider $*^L$ and not $*_!^L$, a justification is needed. We have just seen that, for $E = R(\alpha \times \mathrm{Id})_*(\mathrm{pr}_1, a)_*K$, where $K = K_1 \boxtimes^L K_2$, $(\overline{\mathrm{pr}}_2,E)$ is universally locally acyclic outside $(0,0)$, i.e., $R\Phi_{\overline{pr}_2}E$ is concentrated at $(0,0)$. The triangle deduced from $E|\overline{pr}_2^{-1}(0) \to R\Psi_{\overline{pr}_2}(E) \to R\Phi_{\overline{pr}_2}(E) \to$ by applying $R\Gamma(\overline{pr}_2^{-1}(0),-)$ reads
\begin{equation}\label{triangle}
Ra_*(K)_{0} \to Ra_*(K)_{\overline{\eta}} \to R\Phi_{a_{(0)}}(K)_{(0,0)} \to.
\end{equation}
Indeed, $R\Gamma(\overline{pr}_2^{-1}(0),E|\overline{pr}_2^{-1}(0)) = (R\overline{pr}_{2*}E)_0$ by proper base change, \begin{equation*}
R\Gamma(\overline{pr}_2^{-1}(0),R\Psi_{\overline{pr}_2}(E)) = R\Gamma(\overline{pr}_2^{-1}(\overline{\eta}),E) = (R\overline{pr}_{2*}E)_{\overline{\eta}}
\end{equation*}
by properness of $\overline{pr}_2$, and $R\Gamma(\overline{pr}_2^{-1}(0),R\Phi_{\overline{pr}_2}(E)) = R\Phi_{\overline{pr}_2}E_{(0,0)} = \Phi_a(K)_{(0,0)}$ by concentration of $R\Phi_{\overline{pr}_2}(E)$ at $(0,0)$. As $R\overline{pr}_{2*}E = Ra_*K = K_1 *^L K_2$ by definition of $E$, and taking account (\ref{concentration in degree 1}), we get (\ref{global *vanishing3}), and (\ref{rank formula} by (\ref{triangle}), (\ref{generic rank formula}),  and (\ref{rank at zero}).
\end{proof}
\begin{remark}\label{at infinity} In \cite{Laumon87}, 2.7.1.1 (i)) it is proved that $K_1 *_!^L K_2$ is tamely ramified at $\{\infty\}$, using a Fourier transform for $\overline{\mathbf{Q}}_{\ell}$-sheaves. We give a proof of the similar results for $K_1 *^L K_2$ (and torsion coefficients) in \ref{tameness at infinity}.
\end{remark}

If, for $i = 1,2$, $M_i$ is a lisse sheaf of finitely generated and projective $\Lambda$-modules on $\mathbf{G}_{m,k}$, tamely ramified at $\infty$, we will write
\begin{equation}\label{global *2}
M_1 *_1 M_2 := j^*(j_!M_1 *_1 j_!M_2).
\end{equation}
This is a sheaf on $\mathbf{G}_{m,k}$.

\begin{corollary}\label{global unit} In the situation of \ref{global * formulas}, assume that $K_1$ is concentrated in degree $0$, $(K_1)_0 = 0$, lisse on $\mathbf{G}_{m,k}$, and tamely ramified at $\{0\}$ and $\{\infty\}$, and that $K_2$ is geometrically constant (so that we have (\ref{global * vanishing2})). Then the triangle
\begin{equation}
j_!j^*K_2 \to K_2 \to i_*(K_2)_0  \to,
\end{equation}
where $i : \{0\} \to \mathbf{A}^1_k$ is the inclusion, induces an isomorphism
\begin{equation}\label{global unit2}
K_1 *^L j_!j^*K_2 \iso K_1 \otimes^L e^*((K_2)_0)[-1].
\end{equation}

If $M_i$ ($i = 1,2$) is a lisse sheaf of finitely generated and projective $\Lambda$-modules on $\mathbf{G}_{m,k}$, if $M_1$ is tamely ramified at $\{0\}$ and $\{\infty\}$ and  $M_2$ is geometrically constant, the above isomorphism yields an isomorphism
\begin{equation}\label{global unit4}
M_1 *_1 M_2 \iso M_1 \otimes M_2.
\end{equation}

Similar statements with $K_1$ and $K_2$ (resp. $M_1$ and $M_2$) interchanged.
In particular, we have natural isomorphisms
\begin{equation}\label{global unit3}
\Lambda *_1 M \iso M *_1 \Lambda \iso M,
\end{equation}
for $M$ satisfying the assumption of $M_1$ above.
\end{corollary}

\begin{remark}\label{associativity-global} Let $\mathcal{C}$ be the category of lisse, finitely generated and projective $\Lambda$-modules on $\mathbf{G}_{m,k}$, which are tamely ramified at $\infty$. We leave it to the reader to check that, in addition to (\ref{global unit3}), the functor 
\begin{equation*}
\mathcal{C} \times \mathcal{C} \to \mathcal{C}, (M_1,M_2) \mapsto M_1 *_1 M_2
\end{equation*}
verifies associativity and commutativity isomorphisms (for the associativity, use the Künneth formula for $a_m \times_k a_n : \mathbf{A}^m_k \times_k \mathbf{A}^n_k \to \mathbf{A}^1_k \times_k \mathbf{A}^1_k$ (\cite{SGA5}, III, (1.6.4)), $a_r : \mathbf{A}^r_k \to \mathbf{A}^1_k$ denoting the sum map), satisfying the usual constraints.
\end{remark}

\medskip
\textbf{B. Local additive convolution}

\medskip
The following construction is a slight variant of a construction of Laumon (\cite{Laumon87}, 2.7.2). Denote by $A_h = \SP \,k\{t\}$ (resp. $A_{sh} = \SP \, \overline{k}\{t\}$) the henselization at $0$ (resp. the strict henselization at $\{\overline{0}\}$ of the affine line $A = \mathbf{A}^1_k = \SP \,k[t]$. Denote again by
\begin{equation}\label{sum}
a : (A^2)_h \to A_h 
\end{equation}
(resp. $a : (A^2)_{sh} \to A_{sh}$)
the morphism induced by the sum map $a : A^2 \to A$.  Let $\mathrm{pr}_i : (A^2)_h \to A_h$ (resp. $\mathrm{pr}_i : (A^2)_{sh} \to A_{sh}$) be the morphism induced by the $i$-th projection. We denote by $\overline{\eta}$ a geometric point of $A_{sh}$ above the generic point $\eta$ of $A_{h}$.
 
\begin{definition}\label{local convolution} The \textit{local (additive) convolution functor} 
\begin{equation*}
*^L : D_{ctf}(A_h,\Lambda)\times D_{ctf}(A_h,\Lambda) \to D_{ctf}(A_h,\Lambda)
\end{equation*}
is defined by
\begin{equation}\label{derived local additive convolution}
K_1 *^L K_2 := R\Psi_a(K_1 \boxtimes^L K_2)_{(0,0)},
\end{equation}
where $K_1 \boxtimes^L K_2 := \mathrm{pr}_1^*K_1 \otimes^L \mathrm{pr}_2^*K_2 \in D_{ctf}((A^2)_h,\Lambda)$, and the subscript $(0,0)$ denotes the restriction to the closed point $(0,0)$ of the usual nearby cycles complex for $a : (A^2)_h \to A_h$ [SGA 7 XIII 2.1.1]; in the notation of \ref{vanishing topos}, this is $R\Psi_a(K_1 \boxtimes^L K_2) |(0,0) \ori_A A$ (cf. (\ref{classical psi})), where $(0,0) \ori_A A$ is identified with $(0) \ori_A A \iso A_h$ by 
(\ref{oriented topos of trait}). 

For $K_i \in D_{ctf}(A_{sh},\Lambda)$, we define $K_1 *^L K_2 \in D_{ctf}(A_{sh},\Lambda)$ similarly.  
\end{definition}

For $K_i \in D_{ctf}(A_h,\Lambda)$, we have a distinguished triangle
\begin{equation}
(K_1)_0 \otimes^L (K_2)_0  \to K_1 *^L K_2 \to R\Phi_a(K_1 \boxtimes^L K_2)_{(0,0)} \to,
\end{equation}
so, if $(K_1)_0$ or $(K_2)_0$ is zero,
\begin{equation}\label{* and RPhi}
K_1 *^L K_2 \iso R\Phi_a(K_1 \boxtimes^L K_2)_{(0,0)}.
\end{equation}
On the other hand, by (\ref{section iso}), we have
\begin{equation}\label{stalk of *}
(K_1 *^L K_2) | A_{sh} \iso \sigma_{(\overline{0},\overline{0})}^*R\Psi_a(K_1 \boxtimes^L K_2) \iso (K_1 |A_{sh}) *^L (K_2 |A_{sh}).
\end{equation}

The following result is a variant of (\cite{Laumon87}, 2.7.1.3) :
\begin{proposition}\label{vanishing}
 For $V_i \in D_{ctf}(\eta,\Lambda)$ concentrated in degree zero (i. e. given by a representation of $\mathrm{Gal}(\overline{\eta}/\eta)$ in a finitely generated projective $\Lambda$-module), $j_!V_1 *^L j_!V_2$ is concentrated in degree 1, i. e.
\begin{equation*}
 \mathcal{H}^i(j_!V_1 *^L j_!V_2) = 0
\end{equation*}
for $i \ne 1$, and $\mathcal{H}^1(j_!V_1 *^L j_!V_2)_{\overline{\eta}}$ is projective and of finite type over $\Lambda$.
\end{proposition}
\begin{proof} The proof is the same as that of (\cite{Laumon87}, 2.7.1.3). By the Gabber-Katz extension theorem (\cite{Katz:1986ab}, 1.4) (cf. (\cite{Laumon87}, 2.2.2.2)), we can choose a lisse sheaf $K_i$ of finitely generated and projective $\Lambda$-modules on $\mathbf{G}_{m,k}$ extending $V_i$, i. e. such that $K_i | \eta = V_i$, and $K_i$ is tamely ramified at $\{\infty\}$. Then the result follows from (\ref{global *vanishing3}) and (\ref{* and RPhi}).
\end{proof}

A purely local proof of \ref{vanishing}, in a more general setting, will be given in \ref{Deligne's variant}.

\medskip
In the situation of \ref{vanishing}, similarly to (\ref{global *2}), we will put
\begin{equation}\label{*1}
V_1 *_1 V_2 = j^*\mathcal{H}^1(j_!V_1 *^L j_!V_2)
\end{equation}
(in (\cite{Laumon87}, 2.7.2) the notation $V_1 *  V_2$ is used instead of $V_1 *_1 V_2$). If we denote by 
\begin{equation}\label{G}
\mathcal{G} = \mathcal{G}(\eta,\Lambda)
\end{equation}
the category of sheaves of finitely generated and projective $\Lambda$-modules over $\eta$, we thus have a functor
\begin{equation}\label{*2}
*_1 : \mathcal{G} \times \mathcal{G} \to \mathcal{G}.
\end{equation}
Similarly to \ref{associativity-global}, one checks that it verifies associativity and commutativity constraints (for the associativity, the Künneth formula of (\cite{SGA5}, III, (1.6.4)) has to be replaced by the Künneth formula of \ref{kunneth theorem}, in its punctual form (A.3.1). It also has a two-sided unit (see \ref{local unit} below). 

\begin{proposition}\label{local * vanishing} Let $K_i \in D_{ctf}(A_h,\Lambda)$ ($i = 1,2)$. 

(1) If $K_1$ is concentrated at $\{0\}$, we have a canonical isomorphism 
\begin{equation}\label{local * vanishing0}
i_*(K_1)_0 *^L K_2 \iso e^*((K_1)_0) \otimes^L K_2 
\end{equation}
(resp. the similar one with $K_1$ and $K_2$ interchanged), where $e : A_h \to \SP \,k$ is the projection and $i : \{0\} \to A_h$ the inclusion.

(2) Assume that $K_1$  (resp. $K_2$) has lisse cohomology sheaves and that $(K_2)_0 =0$ (resp. $(K_1)_0 = 0$). Then
\begin{equation}\label{local * vanishing 1}
K_1 *^L K_2 = 0.
\end{equation}
\end{proposition}
\begin{proof} (1) The proof is the same as that of (\ref{global * formulas} (1)).

(2) It suffices to check the resp. assertion. We may assume that $A_h$ is strictly local, hence $K_2$ is constant, so we are reduced to the case where $K_2$ is the constant sheaf $\Lambda$. By (\cite{SGA41/2}, Rapport, 4.6), we may assume that $K_1 = j_!V_1$, where $V_1$ is a bounded complex of sheaves of finitely generated and projective $\Lambda$-modules over $\eta$.  By dévissage, we may assume that $K_1 = j_!V$ with $V$ concentrated in degree zero. By the Gabber-Katz extension theorem, choose a lisse sheaf $L$ on $\mathbf{G}_{m,k}$, of finitely generated and projective $\Lambda$-modules, tamely ramified at $\{\infty\}$, such that $L_{\eta} \iso V$. By (\ref{* and RPhi}),
\begin{equation*}
j_!V *^L \Lambda \iso R\Phi_a(j_!V \boxtimes \Lambda)_{(0,0)} \iso R\Phi_a(j_!L \boxtimes \Lambda)_{(0,0)},
\end{equation*}
which is zero by (\ref{global vanishing phi}). 

Note that, for the resp. assertion of (2), it suffices to show $R\Phi_a(\mathrm{pr}_1^*K_1)_{(0,0)} = 0$. A purely local argument for this, in a more general setting, will be given in \ref{Deligne's variant}.
\end{proof}

In view of (\ref{local * vanishing} (1)), (\ref{local * vanishing} (2)) implies the following local analogue of \ref{global unit} :

\begin{corollary}\label{local unit} Under the resp. assumption of \ref{local * vanishing} (2) the triangle
\begin{equation}
j_!j^*K_2 \to K_2 \to i_*(K_2)_0  \to,
\end{equation}
where $i : \{0\} \to A_h$ is the inclusion, induces an isomorphism
\begin{equation}\label{local unit2}
K_1 *^L j_!j^*K_2 \iso K_1 \otimes^L e^*((K_2)_0)[-1]).
\end{equation}
Similar statement with $K_1$ and $K_2$ interchanged.
For $V_i \in \mathcal{G}(\eta,\Lambda)$, if $V_1$ or $V_2$ is unramified, i. e. extends to a lisse sheaf on $A_h$, the above isomorphism yields an isomorphism
\begin{equation}\label{local unit4}
V_1 *_1 V_2 \iso V_1 \otimes V_2.
\end{equation}
In particular, we have natural isomorphisms
\begin{equation}\label{local unit3}
\Lambda *_1 V \iso V *_1 \Lambda \iso V.
\end{equation}
\end{corollary}

\section{Applications to Thom-Sebastiani type theorems}\label{Thom-Sebastiani type theorems}

In order to formulate an analogue of the Thom-Sebastiani for the sum map, it is convenient to introduce the following generalization of the local additive convolution product \ref{local convolution}. We keep the notation of \ref{interlude}. In particular, $A$ (resp. $A_h$, resp. $A_{sh}$) denotes the affine line $\mathbf{A}^1_k$ (resp. its henselization, resp. strict henselisation) at the origin, $a  : A^2 \to A$ (resp. $a : (A^2)_h \to A_h$, resp. $a : (A^2)_{sh} \to A_{sh}$) denotes the sum map (resp. the map induced by it), and $\overline{\eta}$ a geometric point of $A_{sh}$ above the generic point $\eta$ of $A_h$.

\begin{definition}\label{X-convolution} For $i = 1,2$, let $f_i : X_i \to A_h$ be a morphism of finite type. The \textit{local (additive) convolution along $(X_1,X_2)$} is the functor
\begin{equation}\label{X-convolution1}
*^L : D_{tf}(X_1 \ori_{A_h} A_h,\Lambda) \times D_{tf}(X_2 \ori_{A_h} A_h,\Lambda) \to D^+(X_h \ori_{A_h} A_h,\Lambda)
\end{equation}
defined by 
\begin{equation}\label{X-convolution2}
K_1 *^L K_2 = R\fl a_*K,
\end{equation}
where, in the right hand side of (\ref{X-convolution1}), 
$$
X_h := (X_1 \times_k X_2) \times_{A_h \times_k A_h} (A^2)_h
$$
 is viewed as a scheme over $A_h$ by $a : (A^2)_h \to A_h$, 
$$
K =  \fl{\mathrm{pr}_1}^*K_1 \otimes^L \fl{\mathrm{pr}_2}^*K_2,
$$
$\fl{\mathrm{pr}_i} : X_h \ori_{(A^2)_h} (A^2)_h \to X_i \ori_{A_h} A_h$ is defined by the $i$-th projection, and 
$$
\fl a : X_h \ori_{(A^2)_h} (A^2)_h \to X_h \ori_{A_h} A_h
$$
is induced by $a: (A^2)_h \to A_h$. The subscript tf means finite tor-dimension. 
\end{definition}

\begin{remark} This definition generalizes \ref{local convolution}. Indeed, take $X_i$ to be the closed point $\{0\}$ of $A_h$. Then $\{0\} \ori_{A_h} A_h \iso A_h$, $\{0,0\} \ori_{(A^2)_h} (A^2)_h \iso (A^2)_h$, and for $K_i \in D_{ctf}(A_h,\Lambda)$, by (\ref{stalk of Rg*}) we have $R\fl a_*(K_1  \boxtimes^L K_2) \iso K_1 *^L K_2$.

There is an obvious variant of \ref{X-convolution} with $A_h$ replaced by $A_{sh}$. In the situation of \ref{X-convolution},  let $x_i$ be a geometric point of $X_i$ above $\{0\}$, $x$ the geometric point of $X_h$ defined by $(x_1,x_2)$. We have $x_i \ori_{A_h} A_h \iso A_{sh}$, $x \ori_{A_h} A_h \iso A_{sh}$. Assume that $K_i \in D_{tf}(X_i \ori_{A_h} A_h,\Lambda)$. Let $(K_i)_{|x_i} := K_i |x_i \ori_{A_h} A_h \in D_{tf}(A_{sh},\Lambda)$. Then (by (\ref{stalk of Rg*}) again)
\begin{equation}\label{punctual X-convolution}
(K_1 *^L K_2) | x \ori_{A_h} A_h \iso (K_1)_{|x_1} *^L (K_2)_{|x_2} \in D_{ctf}(A_{sh},\Lambda).
\end{equation}
\end{remark}

\begin{proposition}\label{finiteness of X-convolution} For $K_i \in D_{ctf}(X_i \ori_{A_h} A_h,\Lambda)$, we have 
$$
K_1 *^L K_2  \in D_{ctf}(X \ori_{A_h} A_h,\Lambda).
$$ 
\end{proposition}
\begin{proof} This is proved by a standard dévissage, based on the following facts, whose proof is straightforward from the definitions. Let $Z$ be a scheme separated and of finite type over a trait $S = (s,\eta,\overline{s} \leftarrow \overline{\eta})$. 

(1) The oriented product $Z \ori_S S$ is the union of the closed subtopos $Z_s \ori_S S$ and the complementary open subtopos $Z_{\eta} \times_S S$. 

(a) The topos $Z_{\eta} \times_S S$ is equivalent to $Z_{\eta} \times_{\eta} \eta$, in turn equivalent to $Z_{\eta}$ (as continuous $\mathrm{Gal}(\overline{\eta}/\eta)$-equivariant sheaves on $Z_{\overline{\eta}}$ are just sheaves on $Z_{\eta}$). 

(b) As noted in \ref{vanishing topos}, the topos $Z_s \ori_S S$ is the union of the closed subtopos $Z_s \ori_{S} s$ and the complementary open subtopos $Z_s \ori_S \eta$. The topos $Z_s \ori_{S} s$ is equivalent to $Z_s \ori_s s = Z_s$. The topos $Z_s \ori_S \eta$ is the topos of sheaves $F_{\overline{\eta}}$ on $X_{\overline{s}}$ with a continuous action of $\mathrm{Gal}(\overline{\eta}/\eta)$ compatible with the action of $\mathrm{Gal}(\overline{s}/s)$ on $X_{\overline{s}}$.

(2) A sheaf of $\Lambda$-modules $F$ on $Z \ori_S S$ is constructible if and only if its inverse images on $Z_{\eta}$, $Z_s$, and $Z_s \ori_S \eta$ are.

(3) Let $\widetilde{S}$ be the strict localization of $S$ at $\overline{s}$, and let $F$ be as in (2). Then $F$ is constructible if and only if its inverse image on $X \ori_{\widetilde{S}} \widetilde{S}$ is. 

(4) Let $K$ be an object of $D_{ctf}(Z \ori_S S,\Lambda)$. Then $K$ is isomorphic, in $D(Z \ori_S S,\Lambda)$, to a bounded complex of constructible sheaves of $\Lambda$-modules, projective over $\Lambda$\footnote{Flat and of finite type is equivalent to projective and of finite type.}. This follows from the argument of (\cite{SGA41/2}, Rapport, 4.7, 4.8).

(5) Assume $S$ strictly local, and let $F$ be a constructible sheaf of $\Lambda$-modules on $Z_s \ori_S \eta$, projective over $\Lambda$. Then there exists a finite stratification $Z_s = \coprod Z_{\alpha}$ into locally closed subschemas, and, for each $\alpha$, a finite quotient $G_{\alpha}$ of $G = \mathrm{Gal}(\overline{\eta}/\eta)$, a $\Lambda$-$G_{\alpha}$-module $M_{\alpha}$ projective and of finite type over $\Lambda$, such that the restriction of $F$ to $Z_{\alpha}$ is a locally constant sheaf of  $\Lambda[G_{\alpha}]$-modules of value $M_{\alpha}$. 

Let us now prove \ref{finiteness of X-convolution}. By (3) we may assume $k$ algebraically closed, i. e. replace $A_h$ by $A_{sh}$, which we will denote by $S = (S,s,\eta, \overline{\eta})$. By (4) we may assume that $K_i$ is a constructible sheaf of projective $\Lambda$-modules on $X_i \ori_S S$.  By (1) and (2) we may assume that each $K_i$ is of the form $u_*E_i$ or $v_!F_i$, for $E_i$ (resp. $F_i$) constructible on $(X_i)_s\ori_S S$ (resp. $(X_i)_{\eta}$, where $u : (X_i)_s \ori_S S \inj X_i \ori_S S$, $v : (X_i)_{\eta} \inj X_i \ori_S S$ are the inclusions. As the external tensor product of $u_*E_1$ and $v_!F_2$ (resp. $v_!F_1$ and $u_*E_2$) is zero, we only have to treat the two cases: (i) $K_i = u_*E_i$ ($i = 1,2$), (ii) $K_i = v_!F_i$ ($i = 1,2$). Case (ii) follows from (\cite{SGA41/2}, Th. finitude): we may assume that $a : A^2_h \to A_h$ is induced by a map $b : V_1 \times_k V_2 \to W$, where $V_1$, $V_2$ and $W$ are étale neighborhoods of $\{0\}$ in $\mathbf{A}^1_k$, $K_i$ comes from a constructible sheaf $L_i$ of projective $\Lambda$-modules on a scheme $Z_i$ separated and of finite type over $V_i-\{0\}$ ; then, if $g_i : Z_i \to V_i -\{0\}$, $g = g_1 \times_k g_2 : Z_1 \times_k Z_2 \to (V_1 -\{0\})\times_k (V_2-\{0\})$, $K_1 *^L K_2$ is induced by $Rb_*Rg_*(L_1 \boxtimes L_2)$, which is in $D^b_{ctf}(W-\{0\},\Lambda)$. 

Let's treat case (i), i. e. $K_i = u_*E_i$. As $u$ is a closed subtopos, $u_*E_1 *^L u_*E_2 = u_*(E_1 *^L E_2)$, with $E_1 *^L E_2 \in D(X_s \ori_S S,\Lambda)$ defined in a similar way to (\ref{X-convolution2}). By (1) (b), we need only treat the two cases (a) $E_i = u_*u^*E_i$ ($i = 1,2$), (b) $E_i = v_!v^*E_i$ ($i = 1,2$), where this time $u : (X_i)_s \inj (X_i)_s \ori_S S$ and $v : (X_i)_s \ori_S \eta \inj (X_i)_s \ori_S S$ denote respectively the closed and open subtopoi. We may again assume that $E_1$ and $E_2$ are both of type (a) or both of type (b). In case of type (a), the conclusion follows from (\cite{SGA41/2}, Th. finitude). So let us assume we are in case (b), i. e. $E_i = v_!L_i$, with $L_i$ a sheaf of constructible and projective $\Lambda$-modules on $(X_i)_s \ori_S \eta$.  If $Z_i$ is a (locally closed) subscheme of $(X_i)_s$ and $Z := Z_1 \times_k Z_2$, then (by (\ref{punctual X-convolution}))
\begin{equation}
 (K_1 | Z_1 \ori_S S) *^L (K_2 | Z_2 \ori_S S) \to (K_1 *^L K_2)|Z \ori_S S
\end{equation}
is an isomorphism. Therefore, by (5) we may assume that $L_i$ is a locally constant, constructible sheaf of $\Lambda[G_i]$-modules of value $M_i$, projective over $\Lambda$, for $G_i$ a finite quotient of $\mathrm{Gal}(\overline{\eta}/\eta)$. Let $t_i : T_i \to (X_i)_s$ be a finite étale cover such that $t_i^*L_i$ becomes constant. If $T := T_1 \times_k T_2$, then, by (\ref{punctual X-convolution}) again, the natural map
\begin{equation}
(v_!L_1 |T_1 \ori_S S) *^L (v_!L_2|T_2 \ori_S S) \to (v_!L_1 *^L v_!L_2) | T \ori_S S
\end{equation}
is an isomorphism, so we may assume $L_i$ constant. As similarly the natural map 
$$
(v_!M_1 |(X_1)_s \ori_S S) *^L (v_!M_2|(X_2)_s \ori_S S) \to v_!L_1 *^L v_!L_2 
$$
is an isomorphism, $v_!L_1 *^L v_!L_2$ is then constant, of value $v_!M_1 *^L v_!M_2$ , which is in $D^b_{ctf}(s \ori_S S, \Lambda)$ (= $D^b_{ctf}(S,\Lambda)$) by \ref{vanishing}.
\end{proof}

\begin{remark} I don't know if, more generally, for morphisms of the form $\fl g : X \ori_Y Y \to X \ori_S S$, for $g : Y \to S$ of finite type, with $S$ regular of dimension 1, and $X \to Y$ of finite type, $R\fl g_*$ preserves $D_{ctf}$.
\end{remark}

The following is an analogue of the classical Thom-Sebastiani theorem for the sum map.

\begin{theorem}\label{Thom-Sebastiani} With the notation of \ref{X-convolution}, consider the composite morphism
\begin{equation*}
af_h : X_h \overset{f_h}{\to} (A^2)_h \overset{a}{\to} A_h,
\end{equation*}
where $f_h$ is deduced from $f = f_1 \times_k f_2$ by base change by $(A^2)_h \to A_h \times_k A_h$. Let $K_i \in D_{ctf}(X_i,\Lambda)$, $K = (K_1 \boxtimes K_2) | X_h$.
Then  (\ref{psi of composite}), (\ref{phi of composite}), (\ref{kunneth}) induce isomorphisms in $D_{ctf}(X_h \ori_{A_h} A_h,\Lambda)$
\begin{equation}\label{Thom-Sebastiani1}
(R\Psi_{f_1}K_1) *^L (R\Psi_{f_2}K_2) \iso R\Psi_{af_h}(K) 
\end{equation}
\begin{equation}\label{Thom-Sebastiani2}
(R\Phi_{f_1}K_1) *^L (R\Phi_{f_2}K_2)  \iso R\Phi_{af_h}(K) .
\end{equation}
\end{theorem}

\begin{proof} Up to replacing $A_h$ by an étale neighborhood $Y_i$ of $\{0\}$ in $A$ and extending $(f_i,K_i)$ over $Y_i$, the hypotheses of \ref{kunneth theorem} are satisfied by (\ref{psi-good-examples} (a)). By passing to the limit under such neighborhoods, we deduce that the morphism
\begin{equation}\label{TS3}
c : R\Psi_{f_1}K_1 \boxtimes^L R\Psi_{f_2}K_2 \iso R\Psi_f(K_1 \boxtimes^L K_2).
\end{equation}
of (\ref{kunneth}) is an isomorphism. By (\ref{Milnor tube1}), the base change map
\begin{equation}\label{TS4}
R\Psi_f(K_1 \boxtimes^L K_2) | X_h \ori_{(A^2)_h} (A^2)_h \to R\Psi_{f_h}K
\end{equation}
is an isomorphism. By (\ref{psi of composite}), applying $R\fl a_*$ to the restriction of (\ref{TS3}) to $X_h \ori_{(A^2)_h} (A^2)_h$ yields (\ref{Thom-Sebastiani1}).

As $a$ is universally locally acyclic, by (\ref{locally acyclic}) the map (\ref{phi of composite})
\begin{equation}\label{TS5}
R\Phi_{af_h}K \to R\fl a_* R\Phi_{f_h}K
\end{equation}
is an isomorphism. By (\ref{Rphi triangle}), to deduce (\ref{Thom-Sebastiani2}) it thus suffices to show
\begin{equation}\label{TS6}
R\fl a_*((R\Phi_{f_1}K_1 \boxtimes^L p_1^*K_2) |X_h \ori_{(A^2)_h} (A^2)_h)) = 0.
\end{equation} 
and similarly with $R\Phi_{f_1}K_1 \boxtimes^L p_1^*K_2$ replaced by $p_1^*K_1 \boxtimes^L R\Phi_{f_2}K_2$.
We check (\ref{TS6}) on local sections. Let $L := (R\Phi_{f_1}K_1 \boxtimes^L p_1^*K_2) |X_h \ori_{(A^2)_h} (A^2)_h$. Let $x_i$ be a geometric point of $X_i$ above $\{0\}$, $x$ the corresponding point of $X_h$ (above the closed point $(0,0)$ of $(A^2)_h$). By (\ref{punctual X-convolution}), we have
\begin{equation}\label{TS7}
(R\fl a_*L) | x \ori_{A_h} A_h = (R\Phi_{f_1}K_1)_{|x_1} *^L (p_1^*K_2)_{|x_2}.
\end{equation}
As $(p_1^*K_2)_{|x_2}$ is the \textit{constant} complex on $A_{sh}$ of value $(K_2)_{x_2}$ and $(R\Phi_{f_1}K_1)_{|x_1}$ vanishes at $\{0\}$, \ref{local * vanishing} implies $(R\fl a_*L) | x \ori_{A_h} A_h = 0$, hence (\ref{TS6}). The proof for the similar one is the same. 

\end{proof}

\subsection{}\label{isolated singularities} With the notation of \ref{X-convolution}, let $x_i$ be a rational point of $(X_i)_0$. Assume that $f_i$ is smooth along $(X_i)_0$ outside $\{x_i\}$, flat, and locally of complete intersection of relative dimension $n_i$ at $x_i$ (this last condition is satisfied for example if $X_i$ is regular at $x_i$ of dimension $n_i +1$). Recall (\cite{Illusie03}, 2.10) that under these assumptions $R\Phi_{f_i}(\Lambda)$ is concentrated at $x_i$, and we have
\begin{equation}
R^q\Phi_{f_i}(\Lambda)_{x_i} = 0
\end{equation}  
for $q \ne n_i$, and that $R^{n_i}\Phi_{f_i}(\Lambda)_{x_i}$ is a finitely generated and projective $\Lambda$-module. From (\ref{Thom-Sebastiani2}) we deduce:

\begin{corollary}\label{isolated2} Under the assumptions of \ref{isolated singularities}, the restriction of $R\Phi_{af_h}(\Lambda)$ to $(X_1)_0 \times_k (X_2)_0$ is concentrated at the rational point $x = (x_1,x_2)$,
\begin{equation*}
R^q\Phi_{af_h}(\Lambda)_x = 0
\end{equation*}
for $q \ne n+1$, and, with the notation of  (\ref{*1}), (\ref{Thom-Sebastiani2}) induces an isomorphism of sheaves of $\Lambda$-modules over $\eta$, with finitely generated and projective stalks :
\begin{equation}\label{isolated3}
R^{n_1}\Phi_{f_1}(\Lambda)_{x_1} *_1 R^{n_2}\Phi_{f_2}(\Lambda)_{x_2} \iso R^{n+1}\Phi_{af_h}(\Lambda)_x. 
\end{equation}
\end{corollary}
\begin{remark}\label{isolated4} If in \ref{isolated2} we assume furthermore that $X_i$ is essentially smooth over $k$ and $f_i$ is smooth outside $x_i$, then $af_h$ is essentially smooth outside $x$, and $R\Phi_{af_h}(\Lambda)$ is concentrated at $x$.
\end{remark}
\begin{remark}\label{convolution vs tensor product}
Under the assumptions of \ref{isolated4}, let $\mu(f_i)$ (resp. $\mu(af_h)$) be the Milnor number of $f_i$ (resp. $af_h$) at $x_i$ (resp. $x$) (\cite{SGA7}, XVI 1.2). Assume $\Lambda$ local. By Deligne's theorem (\cite{SGA7} XVI, 2.4) we have
\begin{equation*}
\mu(f_i) = \mathrm{dimtot} \, R^{n_i}\Phi_{f_i}(\Lambda)_{x_i}, \, \, \mu(af_h) = \mathrm{dimtot} \, R^{n+1}\Phi_{af_h}(\Lambda)_x,
\end{equation*}
Here, for $V \in \mathcal{G}(\eta,\Lambda)$ (\ref{*2}), $\mathrm{dim tot} \, V = \mathrm{rk}(V) + \mathrm{sw}(V)$, where $\mathrm{rk}(V)$ is the rank of $V$ over $\Lambda$, and $\mathrm{sw}(V)$ the Swan conductor of $V$. By the definition of the Milnor number, in this situation, we have
\begin{equation*}
\mu(af_h) = \mu(f_1)\mu(f_2).
\end{equation*}
For coefficients $\overline{\mathbf{Q}}_{\ell}$, this formula agrees with (\ref{isolated3}), as by a result of Laumon (\cite{Laumon87}, (2.7.2.1)) we have
\begin{equation*}
\mathrm{dimtot}(V_1 *_1 V_2) = \mathrm{dimtot}(V_1)\mathrm{dimtot}(V_2)
\end{equation*}
for $V_1$, $V_2$ in $\mathcal{G}(\eta,\overline{\mathbf{Q}}_{\ell})$. It is not clear that Laumon's arguments extend to finite coefficients. We give an independent proof in \ref{dimtot}. 
\end{remark}

\subsection{}\label{multiple convolution} If $r$ is an integer $\ge 1$, one defines a multiple convolution product
$$
*^L_{1 \le i \le r} : \prod_{1 \le i \le r} D_{ctf}(A_h,\Lambda) \to D_{ctf}(A_h,\Lambda)
$$
by a formula similar to (\ref{derived local additive convolution})
\begin{equation}\label{generalized convolution}
*^L_{1 \le i \le r}M_i :=R\fl a_*(\boxtimes_{1 \le i \le r}^L M_i)_{(0,\cdots,0)},
\end{equation}
where $\fl a : \{0\} \ori_{(A^r)_h}(A^r)_h \to \{0\} \ori_{A_h} A_h$ is induced by the sum map $a : \mathbf{A}^r_k \to \mathbf{A}^1_k$, $(A^r)_h$ denoting the henselization of $\mathbf{A}^r_k$ at $(0,\cdots,0)$. We also have a variant of \ref{X-convolution} for a family $(X_i)_{1 \le i \le r}$ of schemes of finite type over $A_h$, and $M_i \in D_{ctf}(X_i \ori_{A_h} A_h,\Lambda)$ :
\begin{equation}\label{multiple X-convolution}
*^L M_i := R\fl a_*M \in D_{ctf}(X_h \ori_{A_h} A_h,\Lambda)
\end{equation}
where $X$ is the pull-back to $(A^r)_h$ of $\prod_k X_i$, $M$ the inverse image of $\boxtimes^L M_i$ on $X_h$ and $\fl a : X_h \ori_{(A^r)_h} (A^r)_h \to X_h \ori_{A_h} A_h$ is induced by the sum map $A^r \to A$. 

For $f_i : X_i \to A_h$ of finite type, $K_i \in D_{ctf}(X,\Lambda)$ ($1 \le i \le r$), one gets isomorphisms in $D_{ctf}(X\ori_{A_h} A_h,\Lambda)$ :
\begin{equation}\label{multiple Thom-Sebastiani1}
*^L_{1 \le i \le r} (R\Psi_{f_i}K_i) | X \ori_{A_h} A_h \iso R\Psi_{af}(K) | X \ori_{A_h} A_h 
\end{equation}
\begin{equation}\label{multiple Thom-Sebastiani2}
*^L_{1 \le i \le r} (R\Phi_{f_i}K_i) | X \ori_{A_h} A_h \iso R\Phi_{af_h}(K)| X \ori_{A_h} A_h
\end{equation}
where $f_h$ is deduced from $\prod f_i$ by base change to $A^r_h$.

\begin{remark}\label{remark on Thom-Sebastiani}
(a) Let $\ell$ be a prime not equal to the characteristic of $k$.  The results in \ref{Thom-Sebastiani} imply variants with $\Lambda$ replaced by $\mathbf{Z}_{\ell}$, $\mathbf{Q}_{\ell}$, a finite extension $E_{\lambda}$ of $\mathbf{Q}_{\ell}$, $\mathcal{O}_{E_{\lambda}}$, and $\overline{\mathbf{Q}}_{\ell}$. 

(b) If in \ref{multiple convolution}, $X_i = A_h$ and $f_i$ induced by $t \mapsto t^2$, $R\Phi_{af_h}(\Lambda)$ is concentrated at $\{0\} := \{0,\cdots, 0\}$ and in degree $r-1$, and given by
\begin{equation*}
R^{r-1}\Phi_{af_h}(\Lambda)_{\{0\}} = (*_1)_{1 \le i \le r} R^0\Phi_{f_i}(\Lambda)_{\{0\}}
\end{equation*} 

For $\mathrm{char}(k) \ne 2$,  $R^{r-1}\Phi_{af_h}(\Lambda)_{\{0\}}$ is tame, and of rank 1, but differs from $\otimes_{1 \le i \le r}R^0\Phi_{f_i}(\Lambda)_{\{0\}}$ by arithmetic characters ((\cite{SGA7} XV, 2.2.5, D, E), (\cite{Fu13}, 1.3)).

(c) If $G$ is a smooth commutative group (or monoid) scheme over $k$, one can define a local convolution $*_G$ associated with the product map $m : G \times G \to G$, and we have results similar to \ref{Thom-Sebastiani}, with $A_h$ replaced by the henselization of $G$ at $\{1\}$. The case $G = \mathbf{G}_m$ (resp. $G$ an elliptic curve) seems to be of special interest, in view of \cite{Rojas-Leon:2013aa} (resp. \cite{Sawin14}).  Note that in the case of the multiplicative convolution for $\mathbf{G}_m$ near $\{0\}$ or infinity, studied in \cite{Katz88}, which is equivalent to the case of the local convolution for the multiplicative monoid $A$ near $\{0\}$, while the analogue of (\ref{Thom-Sebastiani1}) will hold, (\ref{Thom-Sebastiani2}) will have to be replaced by formulas involving the nearby cycles of the product map near zero, through (\ref{cone})) and the defect of vanishing of the analogue of  (\ref{TS6}). 

\end{remark}

\subsection{}\label{Deligne's variant}
In \cite{Deligne11} Deligne proposes the following variant of these local convolutions. For $i = 1, 2$, let $C_i$ be a germ of smooth curve over $k$, by which we mean the henselization at a rational point of a smooth curve over $k$. Let $D$ be another germ of smooth curve over $k$. Denote by $o$ the closed points of $C_i$, $D$, and by $C$ the henselization of $C_1 \times_k C_2$ at the point $(o,o)$. Let again denote by $o$ the closed point of $D$. Consider a local $k$-morphism
\begin{equation}
a : C \to D
\end{equation}
satisfying the condition

\medskip
(D) \textit{The restriction of $a$ to the closed subscheme $C_1 \times \{o\}$ (resp. $\{o\} \times C_2$) of $C$ is an isomorphism onto $D$}.

\medskip
For $K_i \in D_{ctf}(C_i,\Lambda)$, let $K_1 \boxtimes^L K_2$ denote the object of $D_{ctf}(C,\Lambda)$ induced by $\mathrm{pr}_1^*K_1 \otimes^L \mathrm{pr}_2^*K_2$. As in \ref{local convolution} we define the $a$-local convolution
$$
*^L_a : D_{ctf}(C_1,\Lambda) \times D_{ctf}(C_2,\Lambda) \to D_{ctf}(D,\Lambda)
$$
by
\begin{equation}
K_1 *^L_a K_2 := R\Psi_a(K_1 \boxtimes^L K_2)_o.
\end{equation}
If $V_i$ is a sheaf of projective, finitely generated $\Lambda$-modules over the generic point $\eta_i$ of $C_i$, and $j : \eta_i \to C_i$ denotes the inclusion, then $\mathcal{H}^q(j_!V_1 *_aj_!V_2) = 0$
 for $q \ne 1$, and $H^1(j_!V_1 *^L_aj_!V_2)$ is projective, finitely generated over $\Lambda$. Indeed, condition (D) implies that $a$ is essentially smooth at $o$, as the map it induces on the Zariski tangent spaces at $o$ is of the form $(x,y) \mapsto \lambda x + \mu y$, with $\lambda$ and $\mu$ in $k^*$. The local acyclicity of smooth maps implies that the Milnor fiber $F = a^{-1}(\overline{\eta})$ (where $\overline{\eta}$ is a geometric point of $D$ over the generic point $\eta$) is irreducible. Therefore, if $z_1$ (resp. $z_2$) is the closed point of $F$ cut out by $C_1 \times \{o\}$ (resp. $\{o\} \times C_2$, and $u : U := F - \{z_1\} - \{z_2\} \inj F$ is the inclusion, then for any lisse sheaf $\mathcal{F}$ of $\Lambda$-modules on $U$, and in particular for $\mathcal{F} = (j_!V_1 \boxtimes j_! V_2)|U$, $H^0(F,u_!\mathcal{F}) = 0$.  As $a$ is of relative dimension 1, $H^q(j_!V_1 *_aj_!V_2) = 0$ for $q > 1$.  
 
As in (\ref{*1}), we define
 \begin{equation}
 V_1 *_{a,1} V_2 := j^*\mathcal{H}^1(j_!V_1 *_aj_!V_2),
 \end{equation}
where $j : \eta \inj D$ is the inclusion.

Furthermore, for any $K_i \in D_{ctf}(C_i,\Lambda)$, one has
\begin{equation}
R\Phi_a(p_i^*K_i)_o = 0,
\end{equation}
where $p_1$ (resp. $p_2$) is induced by $\mathrm{pr}_1$ (resp. $\mathrm{pr}_2$), composed with $a(-,o)$ (resp. $a(o,-)$). This generalizes (\ref{local * vanishing 1}). One uses the isomorphism $\varphi : C \iso C$ induced by $\varphi(x,y) = (a(x,y),x)$ ($C_1$, $C_2$ being identified with $D$ by (D)), which verifies $p_1\varphi = a$, $p_2\varphi = p_1$. One finds that $\varphi_*R\Phi_a(p_1^*K_1)_o = R\Phi_{p_1}(p_2^*K_1)_o$, which is zero as follows (by a limit argument) from the universal local acyclicity for schemes of finite type over a field.

Now, if $f_i : X_i \to C_i$ is a morphism of finite type, and $X := (X_1 \times_k X_2) \times_{C_1 \times_k C_2}C$, one defines $*^L_a$ as in (\ref{X-convolution2}) :
\begin{equation}
K_1 *^L_a K_2 = R\fl a_*K,
\end{equation}
where $K = (K_1 \boxtimes^L K_2)| X \ori_D D$. Then, for $f : X \to C$ induced by $f_1 \times_k f_2$, by the same argument as for Theorem \ref{Thom-Sebastiani}, one gets formulas similar to (\ref{Thom-Sebastiani1}) and (\ref{Thom-Sebastiani2}), namely isomorphisms in $D_{ctf}(X_o \ori_D D,\Lambda)$ : 
\begin{equation}\label{Thom-Sebastiani3}
(R\Psi_{f_1}K_1) *^L_a (R\Psi_{f_2}K_2)  \iso R\Psi_{af}(K) 
\end{equation}
\begin{equation}\label{Thom-Sebastiani4}
(R\Phi_{f_1}K_1) *^L_a (R\Phi_{f_2}K_2)  \iso R\Phi_{af}(K).
\end{equation}

\section{The tame case}\label{tame case}

\subsection{}\label{kummer torsor} 
We keep the notation and hypotheses of \ref{notation}. Let $n \ge 1$ be an integer invertible in $k$. As in (\cite{Laumon83}, 2.2) we denote by 
\begin{equation}\label{kummer torsor1}
\pi_n : \mathcal{K}_n \to \mathbf{G}_{m,k}
\end{equation}
the $\mu_n$-torsor on $\mathbf{G}_{m,k} = \mathbf{A}^1_k - \{0\}$ defined by the exact sequence 
\begin{equation*}
1 \to \mu_n \to \mathbf{G}_{m,k} \overset{\pi_n}{\to} \mathbf{G}_{m,k}  \to 1.
\end{equation*}
(the \textit{Kummer torsor}) ($\mathcal{K}_n = \mathbf{G}_{m,k}$ and $\pi_n$ is the map $a \mapsto a^n$). A lisse sheaf $L$ of $\Lambda$-modules on $\mathbf{G}_{m,k}$ is tamely ramified at $\{0\}$ and $\{\infty\}$ if and only if, for some $n \ge 1$, after a finite extension of $k$ containing $\mu_n(\overline{k})$, $L$ is trivialized by $\pi_n$. 

Assume that $k$ contains $\mu_n(\overline{k})$. If $L$ is trivialized by $\pi_n$, $L$ is recovered from the constant sheaf $\pi_n^*L$ on $\mathcal{K}_n$ of value $L_{\{\overline{1}\}}$ by the canonical isomorphism
\begin{equation}\label{contracted product1}
L \iso (\pi_{n*}\pi_n^*L)^{\mu_n} \iso (\pi_{n*}\Lambda \otimes_{\Lambda} L_{\{\overline{1}\}})^{\mu_n},
\end{equation}
where  $\mu_n$ $(= \mu_n(\overline{k}))$ acts diagonally (on $\pi_{n*}\Lambda$ via its action on $\mathcal{K}_n$, and on $L_{\{\overline{1}\}}$ via $(\mu_n)_{\{\overline{1}\}}$, a quotient of $\pi_1^t(\mathbf{G}_{m,k},\{\overline{1}\})$). Note that, as $\pi_{n*}\Lambda$ is locally free of rank one over $\Lambda[\mu_n]$, $\pi_{n*}\Lambda \otimes_{\Lambda} L_{\{\overline{1}\}}$ is a co-induced module, hence the norm map $N : x \mapsto  \sum_{g \in \mu_n} gx$ on $\pi_{n*}\Lambda \otimes_{\Lambda} L_{\{\overline{1}\}}$ induces an isomorphism
\begin{equation*}
N : \pi_{n*}\Lambda \otimes_{\Lambda[\mu_n]} L_{\{\overline{1}\}} \to (\pi_{n*}\Lambda \otimes_{\Lambda} L_{\{\overline{1}\}})^{\mu_n},
\end{equation*}
so that (\ref{contracted product1}) yields an isomorphism
\begin{equation}\label{contracted product1bis}
L \iso \pi_{n*}\Lambda \otimes_{\Lambda[\mu_n]} L_{\{\overline{1}\}},
\end{equation}
whose stalk at $\overline{1}$ is the identity, as the stalk of $\pi_{n*}\Lambda$ at $\overline{1}$ is naturally identified with $\Lambda[\mu_n(\overline{k})]$.  (In the identification of $\pi_{n*}\Lambda \otimes_{\Lambda[\mu_n]} L_{\{\overline{1}\}}$ with the co-invariants $(\pi_{n*}\Lambda \otimes_{\Lambda} L_{\{\overline{1}\}})_{\pi_n}$ we let $\mu_n$ act on the right on $\pi_{n*}\Lambda$ by $ag = g^{-1}a$.)

If $L$ is only assumed to be \textit{geometrically} trivialized by $\pi_n$, i.e., after extension to $\overline{k}$, we get an isomorphism (\ref{contracted product1bis}) over $\mathbf{G}_{m,\overline{k}}$,
\begin{equation}\label{contracted product2bis}
L_{\overline{k}} \iso \pi_{n*}\Lambda_{\overline{k}} \otimes_{\Lambda[\mu_n]} L_{\{\overline{1}\}},
\end{equation}
where the subscript $\overline{k}$ means inverse image on $\mathbf{G}_{m,\overline{k}}$.

\begin{remark}\label{invariants} The above argument is a particular case of the following standard facts, that we recall for later use. Let $T$ be a scheme, $A$ a commutative ring, $G$ a finite group.

(1) Let $E$ be a lisse sheaf of finitely generated and projective $A$-modules on $T$, equipped with an action of $G$, and let $E_0$ denote the underlying sheaf of $A$-modules, with trivial action of $G$. Then the map
\begin{equation*}
A[G] \otimes_A E_0 \to A[G] \otimes_A E
\end{equation*}
sending $g \otimes x$ to $g \otimes gx$,  is a $G$-equivariant isomorphism, where the right hand side is equipped with the diagonal action of $G$.

(2) Let $B$ be a locally free sheaf of left $A[G]$-modules of rank $1$. For $E$ as in (1), we have $\mathcal{H}^q(G, B \otimes_A E) = 0$ for $q >0$, and the norm map $N : B \otimes_A E \to B \otimes_A E$, $b \otimes x \mapsto \sum g(b \otimes x)$ induces an isomorphism
\begin{equation*}
(B \otimes_A E)_G \, (= B \otimes_{A[G]} E) \to \mathcal{H}^0(G, B \otimes_A E) = (B \otimes_A E)^{G}
\end{equation*}
(where in $B \otimes_{A[G]} E$, $G$ acts on the right on $B$ by $b.g = g^{-1}b$).

(3) Let $t$ be a geometric point of $T$. Let $\pi : P \to T$ be a $G$-torsor on $T$.  Then the natural map $E \to (\pi_*\pi^*E)^G$ is an isomorphism, and, if we assume moreover that $E$ is trivialized by $P$, then, identifying $\pi^*E$ with the constant sheaf of value $E_t$ on $P$ (with its $G$-action), the projection formula and (2) (with $B = \pi_*A$) yield a $G$-equivariant isomorphism
\begin{equation*}
\pi_*A \otimes_{A[G]} E_t \iso E,
\end{equation*}
where $G$ acts diagonally on the left hand side.
\end{remark}

\begin{remark}\label{action of tame pi1}
When $k$ is algebraically closed, the functor $L \mapsto L_{\{1\}}$ is an equivalence from the category of constructible lisse sheaves of $\Lambda$-modules on $\mathbf{G}_{m,k}$ which are tamely ramified at $\{0\}$ and $\{\infty\}$, to the category of representations of $I_t$ on finitely generated $\Lambda$-modules. As $I_t$ is commutative, for such a sheaf $L$, the action of $I_t$ on $L_{\{1\}}$ is $I_t$-equivariant, hence extends uniquely to an action on $L$, compatible with its action on the stalk.

Let $g \in I_t$, with image $g_n$ in $\mu_n$.  Denote by $g^*_L$ the corresponding automorphism of $L$. As (\ref{contracted product1bis}) is functorial and compatible with passing to the stalk at $1$, $g^*_L$ corresponds to the automorphism of the right hand side of (\ref{contracted product1bis}) given by $a \otimes b \mapsto a \otimes g_nb$. Equivalently, $g^*_L$ is given by $a \otimes b \mapsto g^*_na \otimes b$, where now $g^*_n$ denotes the automorphism of $\pi_{n*}\Lambda$ induced by $g$, which is also that given by the action of $g_n \in \mu_n$.   In other words, $g^*_L$ is induced by the automorphism $g^*_{\pi_{n*}\Lambda}$ of the universal object $\pi_{n*}\Lambda$ via $I_t \surj \mu_n$.
\end{remark}

\begin{remark}\label{profinite coefficients} The preceding definitions and results have variants for $\Lambda$ replaced by profinite rings $R$ like $\mathbf{Z}_{\ell}$ ($\ell \ne p$), the ring of integers $O_{\lambda}$ of a finite extension $E_{\lambda}$ of $\mathbf{Q}_{\ell}$, $E_{\lambda}$, or $\overline{\mathbf{Q}}_{\ell}$, and $D_{ctf}(-,\Lambda)$ replaced by $D^b_c(-,R)$, taken in the sense of Deligne \cite{Deligne80a} when $k$ satisfies the cohomological finiteness condition of (\textit{loc. cit.}, 1.1.2 (d)), or in general in the sense of Ekedahl \cite{Ekedahl90}. We will use them freely.  
\end{remark}

\begin{theorem}\label{B} Let $n_i \ge 1$ ($i = 1,2$) integers invertible in $k$, and $r = a_1n_1 = a_2n_2 = \mathrm{lcm}(n_1,n_2)$. Assume that $k$ contains $\mu_r(\overline{k})$. With the notation of (\ref{global *2}), consider
\begin{equation}\label{B(n1,n2)}
\mathcal{B}_{\Lambda,k}(\underline{n}) = \mathcal{B}_{\Lambda,k}(n_1,n_2) = \pi_{n_1*}\Lambda *_1 \pi_{n_2*}\Lambda,
\end{equation}
which is a lisse sheaf of finitely generated and free $\Lambda$-modules on $\mathbf{G}_{m,k}$, equipped with a natural action of $\mu_{\underline{n}} = \mu_{n_1} \times \mu_{n_2}$. Then :

(1) $\mathcal{B}_{\Lambda,k}(\underline{n})$ is geometrically trivialized by $\pi_r$ (cf. \ref{kummer torsor}), in particular, tamely ramified at $\{0\}$ and $\{\infty\}$. Moreover, we have
\begin{equation}\label{B(n1,n2)0}
(j_!\pi_{n_1*}\Lambda *^L j_!\pi_{n_2*}\Lambda)_0 = 0,
\end{equation}
where $j : \mathbf{G}_{m,k} \inj \mathbf{A}^1_k$ is the inclusion.

(2) As a sheaf of $\Lambda[\mu_{\underline{n}}]$-modules, $\mathcal{B}_{\Lambda,k}(\underline{n})$ is locally free of finite type and of rank 1.

(3) The tame inertia group $I_t$ acts on $\mathcal{B}_{\Lambda}(\underline{n})_{\{\overline{1}\}}$ through its quotient   $\mu_r$ via $\mu_r \to \mu_{n_1} \times \mu_{n_2}$, $\varepsilon \mapsto (\varepsilon^{a_1},\varepsilon^{a_2})$.
\end{theorem}

\begin{proof} We may assume $k$ algebraically closed. 

(1) The following argument is borrowed from Deligne's seminar \cite{Deligne80}. Write $G$ for $\mathbf{G}_{m,k}$, $A$ for $\mathbf{A}^1_k$, and drop the subscript $k$ in the products for short. Let $N$ be an integer $\ge 1$. As $k$ is algebraically closed, to prove that a sheaf $L$ on $G$ (for the étale topology) is trivialized by $\pi_N$, it suffices to show that it is equipped with an isomorphism $\alpha : m_N^*L \iso \mathrm{pr}_2^*L$, where $m_N : G \times G \to G$ is the map $(\lambda,t) \mapsto \lambda^Nt$. Indeed, restricting $\alpha$ to the subscheme $G \times \{1\}$ of $G$ gives an isomorphism $\pi_N^*L \iso e^*(L_{\{1\}})$, where $e : G \to \mathrm{Spec} \,k$ is the projection. We'll call such an $\alpha$ a weak action of $G$ on $L$ above $m_N$ (no associativity or unity constraints are required on $\alpha$).

We have
\begin{equation*}
Ra_*((j_!\pi_{n_1*}\Lambda \boxtimes j_!\pi_{n_2*}\Lambda) | A^2 - a^{-1}(0)) = Ra_{\underline{n}*}(u_!\Lambda),
\end{equation*}
where $a_{\underline{n}} : Z \to G$ is the map $(x_1,x_2) \mapsto x_1^{n_1} + x_2^{n_2}$, with 
\begin{equation*}
X = A \times A - a_{\underline{n}}^{-1}(0) - \{0\} \times A - A \times \{0\}
\end{equation*}
and $u$ the open inclusion
\begin{equation*}
u : X \inj Z := A \times A - a_{\underline{n}}^{-1}(0)
\end{equation*}
with complement $D = \{0\} \times A \cup A \times \{0\} - \{(0,0)\}$. Therefore it suffices to show that the pair $(Z,D)$ is trivialized by $\pi_r$, which is equivalent to constructing a weak $G$-action on $a_{\underline{n}} : (Z,D) \to G$ above $m_r$. Such an action is given by
\begin{equation*}
h : G \times Z\to Z, \, \, (\lambda,x_1,x_2) \mapsto (\lambda^{a_1}x_1,\lambda^{a_2}x_2).
\end{equation*} 
Indeed, one readily checks that the diagram
\begin{equation*}
\xymatrix{G \times Z \ar[d]_f \ar[r]^h & Z \ar[d]_{a_{\underline{n}}} \\
G \times G \ar[r]^{m_r} & G}
\end{equation*}
is cartesian, where $f(\lambda,x_1,x_2) = (\lambda,x_1^{n_1} + x_2^{n_2})$. The other verifications are straightforward.

Let us prove (\ref{B(n1,n2)0}). Recall that by (\ref{global * formulas} (3)), we have
\begin{equation*}
(j_!\pi_{n_1*}\Lambda *^L j_!\pi_{n_2*}\Lambda)_0 \iso R\Gamma(a^{-1}(0), j_!\pi_{n_1*}\Lambda \boxtimes j_!\pi_{n_2*}\Lambda | a^{-1}(0)).
\end{equation*}
Let $K$ denote the complex on the right hand side. Let $D := a^{-1}(0)$, $v : D- \{0\} \inj D$ the inclusion,  $L := j_!\pi_{n_1*}\Lambda \boxtimes j_!\pi_{n_2*}\Lambda | D - \{0\}$, so that $K = R\Gamma(D,v_!L)$. We have $H^iK = 0$ for $i = 0$ and $i > 1$. So it suffices to show that $\chi(D, v_!L) = 0$. We have $\chi(D,v_!L) = \chi(D,Rv_*L) = \chi(D -\{0\},L)$ (cf. end of proof of (\ref{global * formulas}, (2))). By the isomorphism $G \iso D-\{0\}, \, t \mapsto (t,-t)$, we have $L \iso v_!(\pi_{n_1*}\Lambda \otimes [-1]^*\pi_{n_2*}\Lambda)$, where $[-1] : G \iso G, \, t \mapsto -t$, so $L$ is tamely ramified at $\{0\}$ and $\{\infty\}$. By the Ogg-Shafarevich formula, we have $\chi(D -\{0\},L) = \chi(D-\{0\})\mathrm{rk}(L) = 0$, which finishes the proof of (\ref{B(n1,n2)0}).

(2) By (\ref{global * formulas} (3)) we know that $\mathcal{H}^q(Ra_*((j_!\pi_{n_1*}\Lambda \boxtimes j_!\pi_{n_2*}\Lambda) | A^2 - a^{-1}(0))$ is zero for $q \ne 1$, and $\mathcal{H}^1 = \pi_{n_1*}\Lambda *_1 \pi_{n_2*}\Lambda$ is a lisse sheaf of finitely generated and projective $\Lambda$-modules on $G$. Let $m$ be its rank. We first show that $m = n_1n_2$. By the second formula of (\ref{global * formulas} (3)), we have
\begin{equation*}
m = \mathrm{rk} \,H^1(a^{-1}(1),M | a^{-1}(1))
\end{equation*}
where 
\begin{equation*}
M = (j_!\pi_{n_1*}\Lambda \boxtimes j_!\pi_{n_2*}\Lambda) | a^{-1}(1) = u_{1!}((\pi_{n_1*}\Lambda \boxtimes \pi_{n_2*}\Lambda) | a^{-1}(1) - \{(0,1)\} - \{(1,0)\}),
\end{equation*}
and $u_1 : a^{-1}(1) - \{(0,1)\} - \{(1,0)\} \inj a^{-1}(1)$ is the open inclusion. Then
\begin{equation*}
m = -\chi(a^{-1}(1),u_{1!}((\pi_{n_1*}\Lambda \boxtimes \pi_{n_2*}\Lambda) | a^{-1}(1) - \{(0,1)\} - \{(1,0)\})).
\end{equation*}
By (\ref{Laumon's theorem}),
\begin{equation*}
m = -\chi(a^{-1}(1),Ru_{1*}(M | a^{-1}(1) - \{(0,1)\} - \{(1,0)\})),
\end{equation*}
so
\begin{equation}\label{B1}
m = -\chi(a^{-1}(1) - \{(0,1)\} - \{(1,0)\},N)
\end{equation}
where $N$ is the (lisse) sheaf $M | a^{-1}(1) - \{(0,1)\} - \{(1,0)\}$. This sheaf is tamely ramified at $\infty$, $\{(0,1)\}$, and $\{(1,0)\}$. Indeed, if $\varphi$ denotes the isomorphism $A \iso a^{-1}(1), \, \, \, t \mapsto (t,1-t)$, we have
\begin{equation*}
\varphi^*N \iso \pi_{n_1*}\Lambda \otimes [t \mapsto 1-t]^*(\pi_{n_2*}\Lambda),
\end{equation*}
and each factor is tamely ramified at $0$, $1$, and $\infty$. 
As $N$ is of rank $n_1n_2$, by the Ogg-Shafarevitch formula, (\ref{B1}) gives
\begin{equation*}
m = -n_1n_2\chi(a^{-1}(1) - \{(0,1)\} - \{(1,0)\},\Lambda) = n_1n_2.
\end{equation*}
As $K = j_!\pi_{n_1*}\Lambda \boxtimes j_!\pi_{n_2*}\Lambda$ is a constructible sheaf of $\Lambda[\mu_{\underline{n}}]$-modules on $A^2$ whose stalks are zero or free of rank 1, $Ra_*K$ is in $D_{ctf}(A,\Lambda)$. As it is concentrated in one degree, it follows that $\mathcal{B}_{\Lambda,k}(\underline{n}) = R^1a_*(K | A^2 - a^{-1}(0))$ is finitely generated and projective over $\Lambda[\mu_{\underline{n}}]$. As its rank over $\Lambda$ is $n_1n_2$, its rank over $\Lambda[\mu_{\underline{n}}]$ is 1.

(3) With the notation of the proof of (1), as $(Z,D)$ is trivialized by $\pi_r$, the action of $I_t$ on $\mathcal{B}_{\Lambda,k}(\underline{n})$ via $\mu_r$ is given by transportation of structure via the action of $\mu_r$ on the second factor of $(Z_r,D_r) := (Z,D) \times_{G} (G,\pi_r)$, the pull-back of $(Z,D)$ by $\pi_r : G \to G, t \mapsto t^r$. We have:
$$
Z_r = \{(x_1,x_2,t) | x_1^{n_1} + x_2^{n_2} = t^r\}.
$$  
Consider the curve
$$
C_{\underline{n}} = (\{(x_1,x_2) | x_1^{n_1} + x_2^{n_2} = 1\} - \{(\zeta_1,0) |\zeta_1^{n_1} = 1 \} - \{(0,\zeta_2) |\zeta_1^{n_2}=1\},
$$
fiber (by $a_{\underline{n}}$) of $Z-D$ over 1. The isomorphism 
$$ 
\alpha_t : (x_1,x_2) \mapsto  (t^{a_1}x_1,t^{a_2}x_2)
$$
from $C_{\underline{n}}$ to $(Z_r-D_r)_t$ defines a trivialization 
$$
\alpha :  C_{\underline{n}} \times G \iso (Z_r-D_r) | G
$$
by which the action of $\mu_r$ on the right hand side corresponds to that on the left hand side given on the first factor by $\varepsilon(x_1,x_2) \to (\varepsilon^{-a_1}x_1,\varepsilon^{-a_2}x_2)$. 
Thus the action of $\varepsilon$ on the stalk at $\{1\}$ of $R^1a_{n*}(u_!\Lambda)$, i.e., $H^1(C_{\underline{n}},\Lambda)$, is the action deduced from $(x_1 \mapsto \varepsilon^{-a_1}x_1, x_2 \mapsto \varepsilon^{-a_2}x_2)$. As the action on the direct image sheaf corresponds to the inverse of the action on the space, the conclusion follows. 
\end{proof}

\begin{definition}\label{universal global} We will call $\mathcal{B}_{\Lambda,k}(\underline{n})$ the \textit{universal (global) convolution sheaf} (relative to $\underline{n} = (n_1,n_2)$). We will write it $\mathcal{B}_{\Lambda}(\underline{n})$ (or even $\mathcal{B}(\underline{n})$) when no confusion may arise.
\end{definition}

\begin{remark}\label{Jacobi sums}
(a) If $k'$ is a perfect extension of $k$, $\mathcal{B}_{\Lambda,k'}(\underline{n})$ is deduced from $\mathcal{B}_{\Lambda,k}(\underline{n})$ by base change from $\mathbf{G}_{m,k}$ to $\mathbf{G}_{m,k'}$. Thus, $\mathcal{B}_{\Lambda,k}(\underline{n})$ is deduced by base change from $\mathcal{B}_{\Lambda,\mathbf{F}_p}(\underline{n})$. 

If, for $i = 1,2$, $n_i$ divides $m_i$, and $s = \mathrm{lcm}(m_1,m_2)$, then $\mathcal{K}_{n_i}$ is deduced from $\mathcal{K}_{m_i}$ by $d_i : \mu_{m_i} \to \mu_{n_i}$ ($d_i = m_i/n_i$), hence $\mathcal{B}(\underline{n})$, as a sheaf of $\Lambda[\mu_{\underline{n}}]$-modules is deduced from $\mathcal{B}(\underline{m})$ by $\underline{d} = (d_1,d_2) : \mu_{\underline{m}} \to \mu_{\underline{n}}$ (and its action of $\mu_r$ from that of $\mu_s$ by $s/r : \mu_s \to \mu_r$).

(b) If $\Lambda$ is a $\mathbf{Z}/\ell^{\nu}\mathbf{Z}$-algebra, with $\ell \ne p$, then, by (\ref{global * formulas}, (3)), $\mathcal{B}_{\Lambda,k}(\underline{n}) = \mathcal{B}_{\mathbf{Z}/\ell^{\nu}\mathbf{Z},k}(\underline{n}) \otimes_{\mathbf{Z}/\ell^{\nu}\mathbf{Z}} \Lambda$. For $\nu \ge 1$ the $\mathcal{B}_{\mathbf{Z}/\ell^{\nu}\mathbf{Z},k}(\underline{n})$ define a $\mathbf{Z}_{\ell}$-sheaf, locally free of finite type and rank one over $\mathbf{Z}_{\ell}[\mu_{\underline{n}}]$,
\begin{equation}
\mathcal{B}_{\mathbf{Z}_{\ell},k}(\underline{n}) := \varprojlim_{\nu}\mathcal{B}_{\mathbf{Z}/\ell^{\nu}\mathbf{Z},k}(\underline{n}) = \pi_{n_1*}\mathbf{Z}_{\ell} *_1 \pi_{n_2*}\mathbf{Z}_{\ell}. 
\end{equation}
By extension of scalars to $R$, where $R$ is $O_{\lambda}$, $E_{\lambda}$, or $\overline{\mathbf{Q}}_{\ell}$, as in \ref{profinite coefficients}, we get an $R$-sheaf, locally free of finite type and rank 1 over $R[\mu_{\underline{n}}]$,
\begin{equation}\label{B_R}
\mathcal{B}_{R,k}(\underline{n}) := R \otimes_{\mathbf{Z}_{\ell}} \mathcal{B}_{\mathbf{Z}_{\ell},k}(\underline{n}) = \pi_{n_1*}R *_1 \pi_{n_2*}R. 
\end{equation}

Assume that $m \ge 1$ is invertible in $k$, $k$ contains $\mu_m = \mu_m(\overline{k})$, and the finite extension $E_{\lambda}$ of $\mathbf{Q}_{\ell}$ contains $\mu_m$.  We have a canonical decomposition (cf. (\cite{Laumon83}, (2.0.6))) 
\begin{equation}\label{characters-bare}
\oplus_{\chi : \mu_m \to E_{\lambda}^*} \mathcal{K}_{m}(\chi^{-1}) \iso \pi_{m*}E_{\lambda},
\end{equation}
where, on the left hand side, $\chi$ runs through the characters $\mu_m \to E_{\lambda}^*$, and $\mathcal{K}_m(\chi^{-1})$ is the lisse sheaf of rank one $\chi(\mathcal{K}_m) = \mathcal{H}om_{\mu_m}(V_{\chi},\pi_{m*}E_{\lambda})$ in the notation of (\cite{Laumon83}, (2.2.1)). This decomposition is in fact $\mu_m$-equivariant, for the action of $\mu_m$ on the left hand side described at the end of \ref{action of tame pi1}. This can be reformulated in the following way. For each character $\chi : \mu_m \to E_{\lambda}^*$, let $V_{\chi}$ denote a 1-dimensional $E_{\lambda}$-vector space, endowed with the action of $\mu_m$ via $\chi$. Then (\ref{characters-bare}) can be rewritten as a $\mu_m$-equivariant decomposition
\begin{equation}\label{characters} 
\oplus_{\chi : \mu_m \to E_{\lambda}^*} V_{\chi} \otimes \mathcal{K}_{m}(\chi^{-1}) \iso \pi_{m*}E_{\lambda},
\end{equation}
given by composition, where $\mu_m$ acts on the left hand side by its action on the first factor. From this and (\ref{B_R}) we deduce (for $\mu_r(\overline{k}) \subset k)$ a canonical $\mu_{\underline{n}}$-equivariant decomposition
\begin{equation}\label{decomposition of B}
\oplus_{(\chi_1,\chi_2)} V_{\chi_1^{-1}} \otimes V_{\chi_2^{-1}} \otimes \mathcal{J}_{E_{\lambda},k}(\chi_1,\chi_2) \iso \mathcal{B}_{E_{\lambda},k} ,
\end{equation}
where
\begin{equation}\label{J sheaf}
\mathcal{J}_{E_{\lambda},k}(\chi_1,\chi_2) := \mathcal{K}_{n_1}(\chi_1) *_1 \mathcal{K}_{n_2}(\chi_2),
\end{equation}
and $\chi_i$ runs through the characters of $\mu_{n_i}$ with values in $E_{\lambda}$. In the sequel we will write $\mathcal{J}(\chi_1,\chi_2)$ for $\mathcal{J}_{E_{\lambda},k}(\chi_1,\chi_2)$. This is a  lisse $E_{\lambda}$-sheaf of rank 1 on $\mathbf{G}_{m,k}$, geometrically trivialized by $\pi_r$ (\ref{B}). The reason for the notation $\mathcal{J}$ is that, if $k$ is a finite field, then this sheaf is related to a Jacobi sum. Indeed, we have the following result, due to Deligne (see (\cite{Laumon83}, 7.3.1, 7.3.4)):
\end{remark}
\begin{proposition}\label{J and Jacobi sum} Assume that $k = \mathbf{F}_q$ with $n_i$ (hence $r$) dividing $q-1$, and that $\chi_1$ and $\chi_2$ are nontrivial. Then, we have a natural isomorphism (of lisse sheaves of rank 1) on $\mathbf{G}_{m,k}$ :
\begin{equation}\label{J and Jacobi sum1}
\mathcal{J}(\chi_1,\chi_2) \iso \mathcal{K}_{n_1}(\chi_1) \otimes \mathcal{K}_{n_2}(\chi_2) \otimes \varepsilon(\chi_1,\chi_2)
\end{equation}
where $\varepsilon(\chi_1,\chi_2)$ is the geometrically constant lisse sheaf inverse image on $\mathbf{G}_{m,k}$ of the sheaf on $\SP \, k$ defined by the restriction of $\mathcal{J}(\chi_1,\chi_2)$ to the rational point $1 \in \mathbf{G}_{m}(k)$, corresponding to the representation of $\mathrm{Gal}(\overline{k}/k)$ sending the geometric Frobenius $F$ relative to $\mathbf{F}_q$ to $q/J(\chi_1,\chi_2)$, where $J(\chi_1,\chi_2)$ is the Jacobi sum\footnote{Our $J(\chi_1,\chi_2)$ differs from Deligne's $J(\chi_1,\chi_2)$ in (\cite{SGA41/2}, Sommes trigonométriques, (4.14.2)) by the sign $(\chi_1\chi_2)(-1)$.}
\begin{equation*}
J(\chi_1,\chi_2) = - \sum_{t_i \in \mathbf{F}_q^*, t_1 + t_2 = 1}\chi_1^{-1}(t_1^{(q-1)/n_1})\chi_2^{-1}(t_2^{(q-1)/n_2})
\end{equation*} 
(a $q$-Weil number of weight 1 (resp. 2) if $\chi_1 \otimes \chi_2$ is not trivial (resp. trivial)). In particular, if $\overline{t}$ is a geometric point over a rational point $t \in \mathbf{G}_m(k)= \mathbf{F}_q^*$, we have
\begin{equation*}
\mathrm{Tr}(F^*,\mathcal{J}(\chi_1,\chi_2)_{\overline{t}}) = \chi_1(t)\chi_2(t)(q/J(\chi_1,\chi_2)).
\end{equation*}
\end{proposition}

\begin{remark}\label{Jacobi-Gauss} (a) For $R$ equal to $O_{\lambda}$ or $O_{\lambda}/\mathfrak{m}^{\nu}$ no decomposition of $\mathcal{B}_{R,k}(\underline{n})$ analogous to (\ref{decomposition of B}) exists, as in general $\mathcal{B}_{R,k}(\underline{n})$ contains unipotent components. 

(b) As the referee observes, instead of quoting \cite{Laumon83}, one can deduce \ref{J and Jacobi sum} from \cite{Laumon87}, as follows. For simplicity, we'll treat only the case where $\chi_1\chi_2$ is nontrivial. We will write $\mathcal{K}$ for $\mathcal{K}_{n_i}$. Let $\mathcal{F}$ be the Fourier transform as in (\ref{convolution !} (b)), with here $\Lambda = E_{\lambda}$. By (\ref{Fourier *bis}) we have
$$
\mathcal{F}(j_!\mathcal{K}(\chi_1) *^L j_!\mathcal{K}(\chi_2)) = \mathcal{F}(j_!\mathcal{K}(\chi_1)) \otimes \mathcal{F}(j_!\mathcal{K}(\chi_2))[-1],
$$
hence, by (\cite{Laumon87}, 1.4.3.1, 1.4.3.2),
\begin{equation}\label{Jacobi-Gauss1}
\mathcal{F}(j_!\mathcal{K}(\chi_1) *^L j_!\mathcal{K}(\chi_2)) = j_!\mathcal{K}((\chi_1\chi_2)^{-1}) \otimes G(\chi_1,\psi)G(\chi_2,\psi)[-1].
\end{equation}
On the other hand, as $\chi_1\chi_2$ is nontrivial, by (\cite{Laumon87}, 1.4.3.1, 1.4.3.2) again, we have
\begin{equation}\label{Jacobi-Gauss2}
\mathcal{F}(j_!\mathcal{K}(\chi_1\chi_2)) = j_!\mathcal{K}((\chi_1\chi_2)^{-1}) \otimes G(\chi_1\chi_2,\psi).
\end{equation}
As 
$$
j_!\mathcal{K}(\chi_1) *^L j_!\mathcal{K}(\chi_2) = j_!(\mathcal{K}(\chi_1) *_1 \mathcal{K}(\chi_2))[-1],
$$
comparing (\ref{Jacobi-Gauss1}) and (\ref{Jacobi-Gauss2}), we get
\begin{equation}\label{Jacobi-Gauss3}
\mathcal{K}(\chi_1) *_1 \mathcal{K}(\chi_2) = \mathcal{K}(\chi_1\chi_2) \otimes \frac{G(\chi_1,\psi)G(\chi_2,\psi)}{G(\chi_1\chi_2,\psi)}.
\end{equation}
This isomorphism is equivalent to (\ref{J and Jacobi sum1}). Indeed, by the classical relation between Jacobi sums and Gauss sums (cf. (\cite{SGA41/2}, Sommes trigonométriques, 4.15.1)), 
\begin{align*}
J(\chi_1,\chi_2) &= \frac{\tau(\chi_1,\psi)\tau(\chi_2,\psi)}{\tau(\chi_1\chi_2,\psi)} \\
{} &= \frac{g(\chi_1^{-1},\psi)g(\chi_2^{-1},\psi)}{g((\chi_1\chi_2)^{-1},\psi)} \\
{} &= \frac{qg(\chi_1\chi_2,\psi^{-1})}{g(\chi_1,\psi^{-1})g(\chi_2,\psi^{-1})},
\end{align*}
hence
$$
\frac{q}{J(\chi_1,\chi_2)} = \frac{g(\chi_1,\psi^{-1})g(\chi_2,\psi^{-1})}{g(\chi_1\chi_2,\psi^{-1})} = \frac{g(\chi_1,\psi)g(\chi_2,\psi)}{g(\chi_1\chi_2,\psi)}.
$$
\end{remark} 
\medskip
Theorem \ref{B} has the following consequence:
\begin{corollary}\label{contracted product2} Under the assumptions of \ref{B}, for $i = 1,2$, let $M_i$ be a lisse sheaf of finitely generated and projective $\Lambda$-modules on $\mathbf{G}_{m,k}$ geometrically trivialized by $\pi_{n_i}$ (cf. \ref{kummer torsor}). Let $V_i$ be the stalk of $M_i$ at $\{\overline{1}\}$ (with its action of $\mu_{n_i}$), and $V = V_1 \otimes V_2$ (with its action of $\mu_{\underline{n}}$). Then the isomorphisms (\ref{contracted product2bis}) for $M_1$ and $M_2$ (over $\mathbf{G}_{m,\overline{k}}$) induce an isomorphism
\begin{equation}\label{contracted product2*}
\mathcal{B}(\underline{n})_{\overline{k}} \otimes_{\Lambda[\mu_{\underline{n}}]} V \iso (M_1 *_1 M_2)_{\overline{k}}
\end{equation}
(with the notation of (\ref{global *2})), where $\mu_{\underline{n}}$ acts on $\mathcal{B}(\underline{n})$ on the right by $bg = g^{-1}b$. In particular, $M_1 *_1 M_2$ is geometrically trivialized by $\pi_r$, and we have
\begin{equation*}
\mathrm{rk}(M_1 *_1 M_2) = \mathrm{rk}(M_1)\mathrm{rk}(M_2).
\end{equation*}
\end{corollary}
\begin{proof} We may assume $k$ algebraically closed. As in the proof of (\ref{B} (1)), consider the open immersion 
\begin{equation*}
u : U = A \times A - a^{-1}(0) - \{0\} \times A - A \times \{0\} \inj A \times A - a^{-1}(0).
\end{equation*}
Let $R_i := \pi_{n_i*}\Lambda$, and write $B$ for $\mathcal{B}(\underline{n})$. By definition, $B = R^1a_*(u_!(R_1 \boxtimes R_2 | U))$, so by (\ref{invariants} (2)) and the fact that $B$ is locally free of rank 1 over $\Lambda[\mu_{\underline{n}}]$, the norm map gives an isomorphism
\begin{equation*}
B \otimes_{\Lambda[\mu_{\underline{n}}]} V  \iso \mathcal{H}^0(\mu_{\underline{n}}, R^1a_*(u_!(R_1 \boxtimes R_2) \otimes_{\Lambda} V)|U).
\end{equation*}
By (\ref{invariants} (2)) and the fact that $R_1 \boxtimes R_2$ is locally free of rank 1 over $\Lambda[\mu_{\underline{n}}]$, we have
\begin{equation*}
\mathcal{H}^0(\mu_{\underline{n}}, (R^1a_*(u_!(R_1 \boxtimes R_2) \otimes_{\Lambda} V)|U) \iso R^1a_*u_!\mathcal{H}^0(\mu_{\underline{n}},(R_1 \boxtimes R_2) \otimes_{\Lambda} V)|U)).
\end{equation*}
As the referee observes, this isomorphism could also be viewed as deduced from the canonical isomorphism
\begin{equation*}
R\underline{\Gamma}^{\mu_{\underline{n}}}Ra_* \iso Ra_*R\underline{\Gamma}^{\mu_{\underline{n}}}
\end{equation*}
applied to $(u_!(R_1 \boxtimes R_2) \otimes_{\Lambda} V)|U$, via the degeneration of the corresponding spectral sequences. By (\ref{invariants} (1)), as $R_i$ is locally free of rank 1 over $\Lambda[\mu_{n_i}]$, the natural map
\begin{equation*}
\mathcal{H}^0(\mu_{n_1},R_1 \otimes V_1) \boxtimes \mathcal{H}^0(\mu_{n_2},R_2 \otimes V_2) \to \mathcal{H}^0(\mu_{\underline{n}}, (R_1 \boxtimes R_2) \otimes_{\Lambda} V))
\end {equation*}
is an isomorphism, and by (\ref{contracted product1bis}) (or (\ref{invariants} (3)), the right hand side is identified with $M_1 \boxtimes M_2$,
so we finally get the isomorphism (\ref{contracted product2*})
\begin{equation*}
B \otimes_{\Lambda[\mu_{\underline{n}}]} V \iso R^1a_*(u_!(M_1 \boxtimes M_2|U)) = M_1 *_1 M_2.
\end{equation*}
\end{proof}

\begin{proposition}\label{tameness at infinity}  For $i = 1, 2$, let $K_i \in D_{ctf}(\mathbf{A}^1_k,\Lambda)$. Assume that  $K_1$ and $K_2$ are tamely ramified at $\{\infty\}$ (\ref{notation}). Then
$K_1*^L K_2$ (\ref{global convolution}) is also tamely ramified at $\{\infty\}$. 
\end{proposition}

As mentioned in the proof of (\ref{global * formulas} (3)), a similar result for $*^L_!$ and $\overline{\mathbf{Q}}_{\ell}$-coefficients was proved by Laumon (\cite{Laumon87}, 2.7.1.1 (i)). Note that here, unlike in (\textit{loc. cit.}), we make no assumption on the ramification of $K_1$ and $K_2$ at finite distance. This generalization was suggested to me by T. Saito. I am indebted to him and to W. Zheng for the following proof, which is simpler than the one I had originally given.

\begin{proof}  We may assume $k$ algebraically closed, and we drop it from the notation. We need a compactification of $\mathbf{A}^2$ and of the sum map $a : \mathbf{A}^2 \to \mathbf{A}^1$. First, embed $\mathbf{A}^2$ in $\mathbf{P}^2$ in the usual way, $(x,y) \mapsto (x:y:1)$, where $(x:y:z)$ are homogeneous coordinates in $\mathbf{P}^2$.  Let $L$ be the line $z =0$. It joins the points $A = (0:1:0)$ and $B = (1:0:0)$. Let $C$ be the point $(1:-1:0)$ of $L$.  Let $g : X \to \mathbf{P}^2$ be the blow-up of $\mathbf{P}^2$ at $A$ and $B$. Denote by $L'$ the strict transform of $L$, and by $E_A = g^{-1}(A)$, $E_B = g^{-1}(B)$ the exceptional divisors. Let $A'$ (resp. $B'$) be the point where $L'$ meets $E_A$ (resp. $E_B$). Let $C'$ be the point of $L'$ corresponding to $C$. The pencil of projective lines in $\mathbf{P}^2$ through  $A$ (resp. $B$) defines a projection $p_1 : X \to \mathbf{P}^1$ (resp. $p_2 : X \to \mathbf{P}^1$), which sends the strict transform of $L$ to $\infty$. The map $f = (p_1,p_2) : X \to \mathbf{P}^1 \times \mathbf{P}^1$ contracts $L'$ to the point $(\infty,\infty)$ of $\mathbf{P}^1 \times \mathbf{P}^1$. As $(L')^2 = -1$,  $f$ is the blow-up of $(\infty,\infty)$ in $\mathbf{P}^1 \times \mathbf{P}^1$ (with exceptional divisor $L'$, $E_A$ (resp. $E_B$) the strict transform of $\mathbf{P}^1 \times \infty$ (resp. $\infty \times \mathbf{P}^1$)). Let $h : Z \to X$ be the blow-up of $C'$ in $X$, and  $h' : Z' \to \mathbf{P}^2$ the blow-up of $C$ in $\mathbf{P}^2$. Let $F = h^{-1}(C')$ and $F' = h'^{-1}(C)$ be the exceptional divisors. The pencil of projective lines 
\begin{equation}\label{lambda-mu}
\lambda(x+y) -\mu z = 0
\end{equation}
through $C$ induces the pencil of parallel lines $x+y - t = 0$ on $\mathbf{A}^2 = \mathbf{P}^2 - L$. It defines a $\mathbf{P}^1$-bundle $a' : Z' \to \mathbf{P}^1$, whose restriction to $F'$ is an isomorphism (see Figure~\ref{fig1}). 
\begin{figure}[h]
\caption{} \label{fig1}
\includegraphics[scale=0.6]{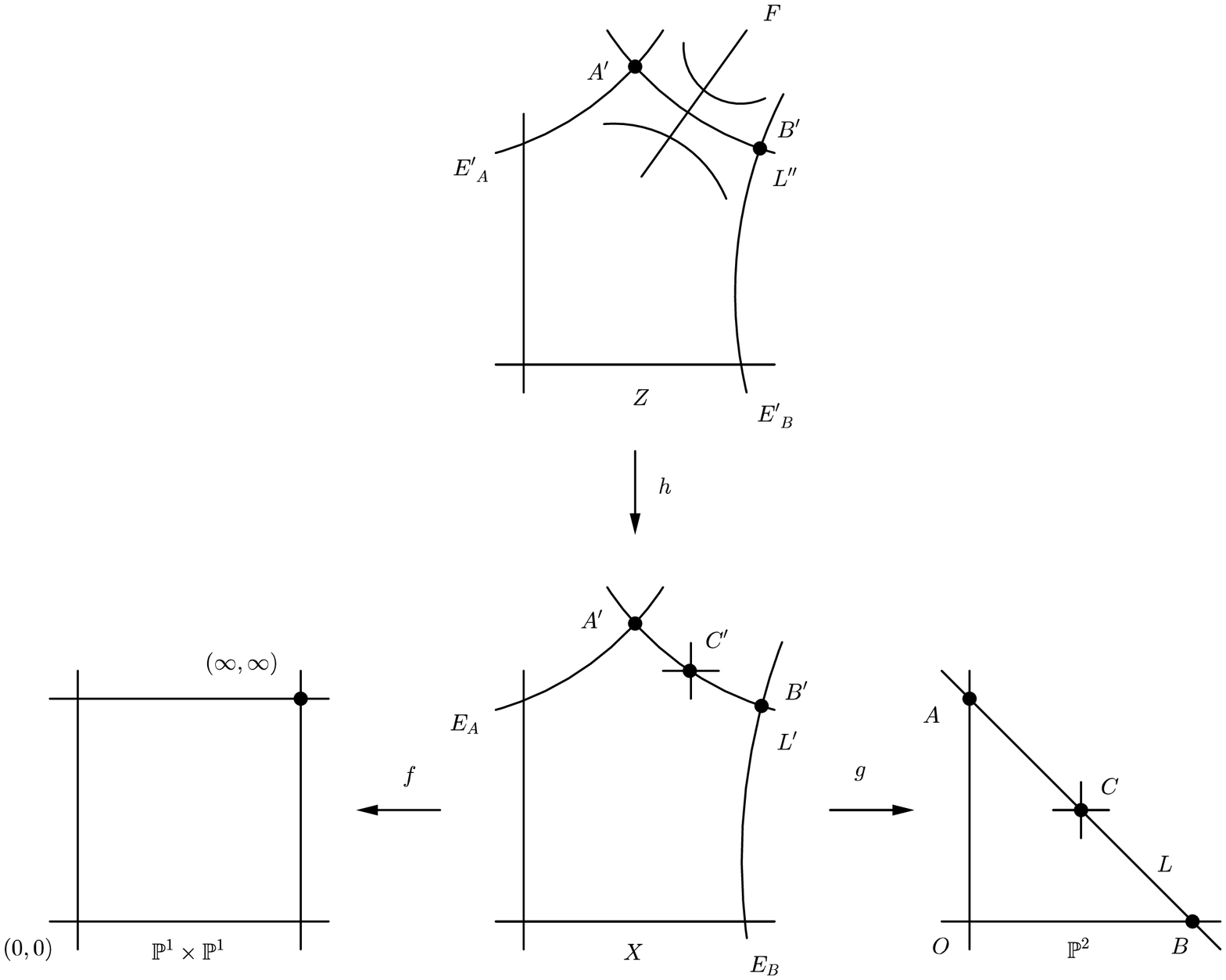}
\end{figure}
We thus get a diagram with cartesian square:
\begin{equation*}
\xymatrix{{} & Z \ar[d]_h \ar[r]^{g'}& Z' \ar[d]_{h'} \ar[r]^{a'} & \mathbf{P}^1\\
\mathbf{P}^1 \times \mathbf{P}^1 & X \ar[l]_f \ar[r]^g & \mathbf{P}^2 & {}}.
\end{equation*}
The map 
\begin{equation*}
\widetilde{a} := a'g' : Z \to \mathbf{P}^1
\end{equation*}
is a compactification of $a : \mathbf{A}^2 \to \mathbf{A}^1$. We have
\begin{equation*}
\widetilde{a}^{-1}(\infty) = E'_A \cup L'' \cup E'_B,
\end{equation*}
where $E'_A$ (resp. $E'_B$, resp. $L''$) denotes the strict transform of $E_A$ (resp. $E_B$, resp. $L'$) in $Z$, and $\infty$ is the point $(\lambda:\mu)$ with coordinate $\lambda = 0$, where $\lambda$ and $\mu$ are defined in (\ref{lambda-mu}). The map $\widetilde{a}$ is proper, and is smooth outside $A'$ and $B'$, at which points it has semistable reduction.  Indeed, $g' : Z \to Z'$ is the blow-up of $Z'$ at the points $A'$ and $B'$ above $A$ and $B$ respectively, and $\widetilde{a}$ is the composition $a'g'$, where $a'$ is smooth. The semistability at $A'$ and $B'$ in such a situation is standard and follows from an elementary calculation (see, e.g., (\cite{Saito:1993aa}, Proof of Prop. 6') or (\cite{Illusie:2004ab}, 1.3)). 

Let $j : \mathbf{A}^2 \inj Z$ be the open immersion. Let $K : = K_1 \boxtimes^L K_2$ on $\mathbf{A}^2$. We have
\begin{equation*}
Ra_*K = R\widetilde{a}_*Rj_*K | (\mathbf{P}^1 - \infty).
\end{equation*} 
We want to show that $R\widetilde{a}_*Rj_*K$ is tame at $\infty$. As $\widetilde{a}$ is proper, it suffices to show that $R\Psi_{\widetilde{a}}(Rj_*K)_{\infty}$ is tame. There are four cases:

(a) Tameness along $E'_A - A'$. The question is étale local at points of $E'_A - A'$, which we identify with $\mathbf{A}^1 \times \infty$ in $\mathbf{A}^1 \times \mathbf{P}^1$ (via $fh$). At such a point $(x_0,\infty) \in (\mathbf{A}^1 \times \mathbf{P}^1)(k)$, we have to show the tameness of $R\Psi_{\widetilde{a}}Rj_*(K)_{(x_0,\infty)}$. We may assume $K_2$ concentrated in degree zero and, by dévissage, using Abhyankar's lemma, we may further assume that $K_2 = \pi_{n*}\Lambda$, for an $n$ invertible in $k$, where $\pi_n$ is the Kummer torsor $y \mapsto y^n$. Let $T = \SP \,k\{z\}$ be the strict localization of $\mathbf{P}^1$ at $\infty$, with parameter $z$ induced by $1/y$, where $(x,y)$ are coordinates on $\mathbf{A}^2$ as above.  Let
\begin{equation*}
W := (\mathbf{A}^1 \times \mathbf{P}^1)_{(x_0,\infty)} = T\{u\} = \SP \, k\{u,z\}
\end{equation*}
be the strict localization of $\mathbf{A}^1 \times \mathbf{P}^1$ at $(x_0,\infty)$, where $u$ is the parameter induced by $x = x_0 +u$.  The morphism $\widetilde{a} : W \to T$ is defined by $\widetilde{a}^*z = z/(1+z(x_0+u))$. Let $\alpha$ be the $k\{u\}$-automorphism of $W$ defined by $\alpha^*(z) = z/(1+z(x_0+u))$, and let $\beta$ be its inverse. By definition, we have $\mathrm{pr}_2\alpha = \widetilde{a}$, where $\mathrm{pr}_2 : W = T\{u\} \to T$ is the canonical projection, hence $\widetilde{a}\beta = \mathrm{pr}_2$, and as $\beta^*R\Psi_{\widetilde{a}} = R\Psi_{\mathrm{pr}_2}\beta^*$, it suffices to show the tameness of $R\Psi_{\mathrm{pr}_2}\beta^*(Rj_*K)$ at the closed point $(0,0)$ of $W$. We have $\beta^*z = z/(1-z(x_0 +u))$, so on $W[z^{-1}]$, $\mathrm{pr}_2\beta$ is induced by $b : \mathbf{A}^2 \to \mathbf{A}^1$, $b(x,y) = y-x$. Therefore, we have
\begin{equation}\label{pull-back}
R\Psi_{\mathrm{pr}_2}\beta^*(Rj_*K) = R\Psi_{\mathrm{pr}_2}(Rj_*(\mathrm{pr}_1^*K_1 \otimes^L b^*K_2)).
\end{equation}
Let $\varepsilon : W' \to W$ be the finite étale cover of $W$ of equation $v^n = 1 -z(x_0 +u)$. The pull-backs by $\varepsilon$ of $b^*K_2$ and $\mathrm{pr}_2^*K_2$ are isomorphic. Therefore, by (\ref{pull-back}) we are reduced to showing the tameness at $(x_0,\infty)$ of $R\Psi_{\mathrm{pr}_2}(Rj_*(\mathrm{pr}_1^*K_1 \otimes^L \mathrm{pr}_2^*K_2))$, where $j = (\mathrm{Id},j_2) : \mathbf{A}^1 \times \mathbf{A}^1 \inj \mathbf{A}^1 \times \mathbf{P}^1$.
By Künneth, we have
\begin{equation*}
R\Psi_{\mathrm{pr}_2}(Rj_*(\mathrm{pr}_1^*K_1 \otimes^L \mathrm{pr}_2^*K_2)) = R\Psi_{\mathrm{pr}_2}(\mathrm{pr}_1^*K_1 \otimes^L \mathrm{pr}_2^*Rj_{2*}K_2),
\end{equation*}
and by the projection formula,
\begin{equation*}
R\Psi_{\mathrm{pr}_2}(\mathrm{pr}_1^*K_1 \otimes^L \mathrm{pr}_2^*Rj_{2*}K_2)_{(x_0,\infty)} =  R\Psi_{\mathrm{pr}_2}(\mathrm{pr}_1^*K_1)_{(x_0,\infty)} \otimes^L (K_2)_{\overline{\eta}},
\end{equation*}
where $\overline{\eta}$ is a geometric generic point of $T$. By Deligne's generic universal local acyclicity theorem (\cite{SGA41/2}, Th. finitude, 2.16), $R\Phi_{\mathrm{pr}_2}(\mathrm{pr}_1^*K_1)_{(x_0,\infty)} = 0$. Therefore, the inertia $I$ acts tamely on $R\Psi_{\widetilde{a}}(Rj_*K)_{\infty}$ along $E'_A - A'$. 

(b) Along $E'_B - B'$, the same calculation gives the desired tameness.

(c) Around $A'$ or $B'$, $\widetilde{a}$ is isomorphic to $m : \mathbf{A}^2 \to \mathbf{A}^1, m(x,y) = xy$ \footnote{The local coordinates $(x,y)$ here are not those used in (a).}, and, by the description of  $X$ as the blow-up of $(\infty,\infty)$ in $\mathbf{P}^1 \times \mathbf{P}^1$, $Rj_*K$, restricted to $(xy \ne 0)$ is an external product of tamely ramified complexes, hence tamely ramified along $m^{-1}(0)$. 
For any lisse sheaf $\mathcal{F}$ on $(xy \ne 0)$ which is tamely ramified along $(xy = 0)$, $R\Psi_m(\mathcal{F})$ is tame, as follows by Abhyankar's lemma from the case where $\mathcal{F} = \pi_*\Lambda$, for a tame cover $\pi$ of $\mathbf{A}^2$ of the form $x \mapsto x^n, y \mapsto y^n$. Therefore, $R\Psi_{\widetilde{a}}(Rj_*K)_{\infty}$ is tame around $A'$ and $B'$.

(d) Along $L'' - A' - B'$, $\widetilde{a}$ is smooth, and $Rj_*K$ is tamely ramified, hence by the last argument in (c), $R\Psi_{\widetilde{a}}(Rj_*K)_{\infty}$ is tame along $L''- A' - B'$, which finishes the proof. 
\end{proof}

The following corollary was suggested to me by T. Saito. 
\begin{corollary}\label{dimtot} Let $(A_h,0,\eta)$ be the henselization of $\mathbf{A}^1_k$ at $0$.  For $i = 1, 2$, let $V_i \in \mathcal{G}(\eta,\Lambda)$ (\ref{G}). Then, with the notation of \ref{convolution vs tensor product} and \ref{*1}, we have:
\begin{equation}\label{dimtot formula}
\mathrm{dimtot}(V_1 *_1 V_2) = \mathrm{dimtot}(V_1)\mathrm{dimtot}(V_2).
\end{equation}
Let $r_i := \mathrm{rk}(V_i)$, $s_i := \mathrm{sw}(V_i)$, so that $\mathrm{dimtot}(V_i) = r_i + s_i$, and let $r := \mathrm{rk}((V_1 *_1 V_2)$, $s := \mathrm{sw}(V_1 *_1 V_2)$. We have
\begin{equation}\label{rank formula 1}
r = r_1r_2 + r_1s_2 + r_2s_1 - \mathrm{sw}(V_1 \otimes [-1]^*V_2),
\end{equation}
\begin{equation}\label{swan formula}
s = s_1s_2 + \mathrm{sw}(V_1 \otimes [-1]^*V_2).
\end{equation}
\end{corollary}

As mentioned in \ref{convolution vs tensor product}, the analogue of (\ref{dimtot formula}) for $\overline{\mathbf{Q}}_{\ell}$-coefficients was proved by Laumon in (\cite{Laumon87}, (2.7.2.1)). He also proved (\ref{rank formula 1}) and (\ref{swan formula}) in the same setting. The proof below is similar, except that we use $*$ instead of $*_!$ for the global convolution. 

\begin{proof} We may and will assume $k$ algebraically closed. By the Gabber-Katz theorem \cite{Katz:1986ab}, extend $V_i$ to a lisse sheaf $\mathcal{F}_i$ on $(\mathbf{G}_m)_k$, tame at $\infty$. Put $A:= \mathbf{A}^1_k$. Let $j : A - \{0\} \inj A$ denote the inclusion. Recall (\ref{global *vanishing3}) that
\begin{equation}\label{local-global *}
V_1 *_1 V_2 = R^0\Phi_{\mathrm{id}}(j_!\mathcal{F}_1 *^L j_!\mathcal{F}_2[1])_0,
\end{equation}
where $j_!\mathcal{F}_1 *^L j_!\mathcal{F}_2 = Ra_*(j_!\mathcal{F}_1 \boxtimes j_!\mathcal{F}_2)$, $a : A^2 \to A$ being the global sum map. Formula (\ref{rank formula 1}) follows from (\ref{rank formula}). This does not use \ref{tameness at infinity}. We will use \ref{tameness at infinity} to prove (\ref{swan formula}). For this, we compute the Euler-Poincaré characteristic $\chi(A, Ra_*(j_!\mathcal{F}_1 \boxtimes j_!\mathcal{F}_2))$ in two ways. First,
\begin{equation*}
\chi(A, Ra_*(j_!\mathcal{F}_1 \boxtimes j_!\mathcal{F}_2)) = \chi(A^2,j_!\mathcal{F}_1 \boxtimes j_!\mathcal{F}_2),
\end{equation*}
and by Künneth (\cite{SGA5} III, (1.6.4)),
\begin{equation*}
\chi(A^2,j_!\mathcal{F}_1 \boxtimes j_!\mathcal{F}_2) = \chi(A,j_!\mathcal{F}_1)\chi(A,j_!\mathcal{F}_2).
\end{equation*}
As $\mathcal{F}_i$ is tame at $\infty$, by the Grothendieck-Ogg-Shafarevich formula, and remembering that $\chi = \chi_c$   \cite{Laumon81}, we have
\begin{equation*}
\chi(A, j_!\mathcal{F}_i) = -\mathrm{sw}_0(\mathcal{F}_i),
\end{equation*}
hence
\begin{equation*}
\chi(A,j_!\mathcal{F}_1 *^L j_!\mathcal{F}_2[1]) = -s_1s_2.
\end{equation*}
By the Grothendieck-Ogg-Shafarevich formula again, and (\ref{rank at zero}), using that $\mathcal{F}_1 *_1 \mathcal{F}_2$ is tame at $\infty$ by \ref{tameness at infinity}, we have, for $L:=j_!\mathcal{F}_1 *^L j_!\mathcal{F}_2[1]$,
\begin{equation*}
\chi(A,L) = \mathrm{rk}_0(L) + \chi_c(A - \{0\},L) = \mathrm{sw}(V_1 \otimes [-1]^*V_2) -s
\end{equation*}
which gives (\ref{swan formula}). Combining (\ref{swan formula}) and (\ref{rank formula 1}),  we get (\ref{dimtot formula}). 
\end{proof}

\begin{remark} (a) The calculation made above shows that, conversely, (\ref{dimtot formula}) (or, equivalently, (\ref{swan formula})) implies \ref{tameness at infinity}.

(b) Formula (\ref{dimtot formula}) is essentially (\cite{Laumon:1983ac}, Exemples 2.3.8 (a)), as the sum map $a$ gives a good pencil in the situation of (\textit{loc. cit.}). 

(c) T. Saito gave an independent proof of (\ref{dimtot formula}) as a corollary of his theory of the characteristic cycle, more precisely, of the compatibility of the formation of the characteristic cycle with external products (\cite{TSaito16}, 2.6, 2.7).
\end{remark}

\begin{remark}\label{associativity etc} Let $\mathcal{C} = \mathcal{C}(\mathbf{G}_{m,k},\Lambda)$ be the full subcategory of sheaves of finitely generated and projective $\Lambda$-modules $M$ on $\mathbf{G}_{m,k}$ which are tamely ramified at $\{0\}$ and $\{\infty\}$. Then, as in (\ref{*2}), we have a functor
\begin{equation*}
\mathcal{C} \times \mathcal{C} \to \mathcal{C}, \, \, (M_1,M_2) \mapsto M_1 *_1 M_2,
\end{equation*}
satisfying associativity and commutativity constraints, and, by (\ref{global unit3}), having the constant sheaf $\Lambda$ as a two-sided unit.
\end{remark} 

\begin{remark} The isomorphism (\ref{contracted product2*}) is compatible with the decomposition (\ref{decomposition of B}), in the following sense. Let $M_i = \mathcal{K}_n(\chi_i)$, for characters $\chi_i : \mu_{n_i} \to E_{\lambda}^*$, so that $V_i = V_{\chi_i^{-1}}$. Then, together with the decomposition (\ref{decomposition of B}), the isomorphism (\ref{contracted product2*}) reads 
\begin{equation*}
\oplus_{\chi'_1,\chi'_2}  \mathcal{J}_{E_{\lambda},k}({\chi'}_1,{\chi'}_2)\otimes_{E_{\lambda}} (V_{{\chi'}_1^{-1}} \otimes V_{{\chi'}_2^{-1}})\otimes_{E_{\lambda}[\mu_{\underline{n}}]} (V_{\chi_1^{-1}} \otimes V_{\chi_2^{-1}}) \iso M_1 *_1 M_2,
\end{equation*}
and the left hand side boils down to $\mathcal{J}_{E_{\lambda},k}(\chi_1,\chi_2)$.
\end{remark}
\begin{corollary}\label{tame monodromy} Assume $k$ algebraically closed. Let $g \in I_t$. Under the assumptions of \ref{contracted product2}, let $g_i^*$ denotes the automorphism of $M_i$ induced by the automorphism of the stalk of $M_i$ at $1$ given by $g$ (cf. \ref{action of tame pi1}). Then the automorphism $g^*$ of $M_1 *_1 M_2$ given by the action of $g$ on $(M_1 *_1 M_2)_1$ is :
\begin{equation}\label{tame monodromy1}
g^* = g_1^* *_1 g_2^*.
\end{equation}
\end{corollary}
\begin{proof}
This follows from (\ref{contracted product2*}) by the same argument as the one used in \ref{action of tame pi1} to describe $g^*_L$.
\end{proof}

\subsection{}\label{local variants} The constructions and statements above have local counterparts, that we will now discuss. We keep the notation of the beginning of \ref{Thom-Sebastiani type theorems}. We will now denote by $I \subset \mathrm{Gal}(\overline{\eta}/\eta)$ the inertia group, and $I_t$ its tame quotient (which is identified with the group denoted $I_t$ in (\ref{pi1tame}) by the map $A_{sh} \to A$ and the choice of a path from $\overline{\eta}$ to $\{\overline{1}\}$). For $n \ge 1$ invertible in $k$, we denote by
\begin{equation}\label{Kummer loc}
\pi_n^{\mathrm{loc}} : \mathcal{K}_n^{\mathrm{loc}} \to A_h
\end{equation}
the inverse image on $\eta \in A_h$ of the Kummer torsor $\mathcal{K}_n$ (\ref{kummer torsor}). For $L \in \mathcal{G}(\eta,\Lambda)$ (\ref{G}), if $\mu_n(\overline{k}) \subset k$ and $L$ is geometrically trivialized by $\pi_n^{\mathrm{loc}}$, we have a canonical isomorphism similar to (\ref{contracted product2}):
\begin{equation}\label{contracted product local}
L_{\overline{k}} \iso \pi_{n*}^{\mathrm{loc}}\Lambda_{\overline{k}} \otimes_{\Lambda[\mu_n]} L_{\overline{\eta}}.
\end{equation}

For a pair of integers $\underline{n} = (n_1,n_2)$ as in \ref{B}, with $\mu_r(\overline{k}) \subset k$ ($r = \mathrm{lcm}(n_1,n_2)$), we denote by
\begin{equation}\label{B loc}
\mathcal{B}_{\Lambda,k}^{\mathrm{loc}}(\underline{n})
\end{equation}
the inverse image of $\mathcal{B}_{\Lambda,k}(\underline{n})$ (\ref{B(n1,n2)}) on $\eta$. This is a tamely ramified representation of $\mathrm{Gal}(\overline{\eta}/\eta)$, geometrically trivialized by $\pi_{r}$, in the notation of $\ref{B}$. We will call $\mathcal{B}_{\Lambda,k}^{\mathrm{loc}}(\underline{n})$ the \textit{universal (local) convolution sheaf} (relative to $\underline{n} = (n_1,n_2)$), and, as in \ref{universal global}, write it $\mathcal{B}^{\mathrm{loc}}_{\Lambda}(\underline{n})$ or $\mathcal{B}^{\mathrm{loc}}(\underline{n})$ when no confusion may arise.

\begin{proposition}\label{B loc properties} (1) The isomorphism (\ref{global *vanishing3}) for $K_i = j_!\pi_{n_i*}\Lambda$ induces an isomorphism
\begin{equation}\label{B loc 1}
\pi_{n_1*}^{\mathrm{loc}}\Lambda *_1 \pi_{n_2*}^{\mathrm{loc}}\Lambda \iso \mathcal{B}_{\Lambda,k}^{\mathrm{loc}}(\underline{n}),
\end{equation}
where the left hand site denotes a local convolution product (\ref{*1}).

(2) The stalk of $\mathcal{B}_{\Lambda,k}^{\mathrm{loc}}(\underline{n})$ at $\overline{\eta}$ is finitely generated and free of rank 1 over $\Lambda[\mu_{\underline{n}}]$, and the tame inertia group $I_t$ acts on it through its quotient $\mu_r$ via $\mu_r \to \mu_{n_1} \times \mu_{n_2}$, $\varepsilon \mapsto (\varepsilon^{a_1},\varepsilon^{a_2})$.
\end{proposition}
\begin{proof} As $K_i$ is tamely ramified at zero, we have, by (\ref{global *vanishing3}) and (\ref{rank at zero}), 
\begin{equation*}
\pi_{n_1*}^{\mathrm{loc}}\Lambda *_1 \pi_{n_2*}^{\mathrm{loc}}\Lambda \iso \mathcal{H}^1(K_1 *^L K_2)_{\eta},
\end{equation*}
which proves (1). (2) follows from (\ref{B} (2), (3)).
\end{proof}

\begin{remark}\label{tameness of *} By (\ref{swan formula}), for any tame representations $L_i \in \mathcal{G}(\eta,\Lambda)$ ($i = 1,2$), $L_1 *_1 L_2$ is tame. However, as in \ref{contracted product2}, we have a more precise result:
\end{remark}

\begin{corollary}\label{contracted product2 local} For $i = 1,2$, let $n_i$  be as in \ref{B}, and let $L_i \in \mathcal{G}(\eta,\Lambda)$ be geometrically trivialized by $\pi_{n_i}^{\mathrm{loc}}$. Let $V_i = (L_i)_{\overline{\eta}}$ (with its action of $\mu_{n_i}$), and $V = V_1 \otimes V_2$ (with its action of $\mu_{\underline{n}}$). Then the isomorphisms (\ref{contracted product local}) for $L_1$ and $L_2$ induce an isomorphism
\begin{equation}\label{contracted product2* local}
\mathcal{B}^{\mathrm{loc}}(\underline{n})_{\overline{k}} \otimes_{\Lambda[\mu_{\underline{n}}]} V \iso L_1 *_1 L_2
\end{equation}
(with the notation of (\ref{*1})), where $\mu_{\underline{n}}$ acts on $\mathcal{B}(\underline{n})$ on the right by $bg = g^{-1}b$. In particular, $L_1 *_1 L_2$ is geometrically trivialized by $\pi_r$, and we have
\begin{equation*}
\mathrm{rk}(L_1 *_1 L_2) = \mathrm{rk}(L_1)\mathrm{rk}(L_2).
\end{equation*}
\end{corollary}
\begin{proof} One can imitate the proof of \ref{contracted product2}, or deduce the result from it, using the Gabber-Katz extension theorem. 
\end{proof}

\begin{remark}\label{Gt} Let $\mathcal{G}_t$ be the full subcategory of $\mathcal{G} = \mathcal{G}(\eta,\Lambda)$ consisting of tame representations. It follows from \ref{tameness of *} or \ref{contracted product2 local} that $\mathcal{G}_t$ is stable under $*_1$. 
\end{remark}

\begin{remark}\label{action of It} Assume $k$ is algebraically closed.  Let $L \in \mathcal{G}(\eta,\Lambda)$ be a tame representation. As in \ref{action of tame pi1}, since $I_t$ is commutative, the action of $I_t$ on $L_{\{1\}}$ is $I_t$-equivariant, hence extends uniquely to an action on $L$, compatible with its action on $L_{\overline{\eta}}$. For $g \in I_t$, we will denote by $g^*$ the automorphism of $L$ induced by the automorphism $g$ of $L_{\overline{\eta}}$.
\end{remark}

\begin{corollary}\label{local tame monodromy} Assume $k$ algebraically closed. Let $g \in I_t$. Under the assumptions of \ref{contracted product2 local}, let $g_i^*$ denotes the automorphism of $L_i$ induced by the automorphism of $(L_i)_{\overline{\eta}}$ given by $g$ (cf. \ref{action of It}). Then the automorphism $g^*$ of $L_1 *_1 L_2$ given by the action of $g$ on $(L_1 *_1 L_2)_{\overline{\eta}}$ is :
\begin{equation}\label{local tame monodromy1}
g^* = g_1^* *_1 g_2^*.
\end{equation}
\end{corollary}
\begin{proof} The proof is the same as for \ref{tame monodromy}.
\end{proof}

This implies the following ``convolution" variant of the Thom-Sebastiani theorem for the monodromy operators:
\begin{corollary}\label{TS - monodromy *} Under the assumptions of \ref{isolated2}, with $k$ algebraically closed, assume that, for $i = 1,2$, $R^{n_i}\Phi_{f_i}(\Lambda)_{x_i}$ is tame. Then $R^{n+1}\Phi_{af_h}(\Lambda)_x$ is tame, and if for $g \in I_t$, $g_i^*$ denotes the monodromy automorphism of $R^{n_i}\Phi_{f_i}(\Lambda)_{x_i}$, and $g^*$ that of $R^{n+1}\Phi_{af_h}(\Lambda)_x$, with the identification (\ref{isolated3}) we have
\begin{equation}
g^* = g_1^* *_1 g_2^*.
\end{equation}
\end{corollary}

Actually, one can recover the original formulation of the Thom-Sebastiani theorem ((0.1), (0.2)), involving a tensor product rather than a convolution. This is a consequence of the following corollaries of \ref{B} and \ref{B loc properties}.

\begin{corollary}\label{* = otimes over kbar} Assume $k$ algebraically closed. 
There exists a projective system of $\mu_{\underline{n}}$-equivariant isomorphisms
\begin{equation}\label{* = otimes over kbar1} 
\beta_{\underline{n}} : \mathcal{B}_{\Lambda,k}(\underline{n}) \iso \pi_{n_1*}\Lambda \otimes \pi_{n_2*}\Lambda
\end{equation}
(resp.
\begin{equation}\label{* = otimes over kbar2}
\beta_{\underline{n}}^{\mathrm{loc}} : \mathcal{B}_{\Lambda,k}(\underline{n})^{\mathrm{loc}} \iso \pi_{n_1*}^{\mathrm{loc}}\Lambda \otimes \pi_{n_2*}^{\mathrm{loc}}\Lambda),
\end{equation}
for $\underline{n} = (n_1,n_2)$ running through the pairs of integers $\ge 1$ invertible in $k$, ordered by divisibility.
\end{corollary}
\begin{proof} By \ref{B loc properties}, it suffices to prove the assertion for the global convolution. First, fix $n_i$, let $r = \mathrm{gcd}(n_1,n_2)$, and let's prove the existence of an isomorphism (\ref{* = otimes over kbar1}). The functor $M \to M_{\{1\}}$ is an equivalence from the category of lisse $\Lambda[\mu_{\underline{n}}]$-modules on $\mathbf{G}_{m,k}$ tamely ramified at $\{0\}$ and $\{\infty\}$ to the category of (continuous) $\Lambda[\mu_{\underline{n}}][I_t]$-modules, finitely generated over $\Lambda$. The stalks at $\{1\}$ of $\mathcal{B}_{\Lambda,k}(\underline{n})$ and $\pi_{n_1*}\Lambda \otimes \pi_{n_2*}\Lambda$ are both free of rank 1 over $\Lambda[\mu_{\underline{n}}]$, and $I_t$ acts on them through the diagonal map $\mu_r \to \mu_{\underline{n}}$ (\ref{B} (3)). They are therefore isomorphic, and this proves the existence of an isomorphism $\beta_{\underline{n}}$. To prove the existence of a projective system of isomorphisms $\beta_{\underline{n}}$'s, by (\ref{Jacobi sums} (a)), it suffices to do it for the diagonal system $\underline{n} = (n,n)$. Isomorphisms $\beta_{(n,n)}$ form a torsor under $(\Lambda[\mu_n \times \mu_n])^*$. As $(\Lambda[\mu_{n+1}\times\mu_{n+1}])^* \to (\Lambda[\mu_n \times \mu_n])^*$ is surjective, the conclusion follows.
\end{proof}

By \ref{contracted product2} (resp. \ref{contracted product2 local}), this implies :
\begin{corollary}\label{* = otimes over kbar3} Under the assumptions of \ref{* = otimes over kbar}, the functors on $\mathcal{C} \times \mathcal{C}$ (\ref{associativity etc}) (resp. $\mathcal{G}_t \times \mathcal{G}_t$ (\ref{Gt})), $(M_1,M_2) \mapsto M_1 *_1 M_2$ and $(M_1,M_2) \mapsto M_1 \otimes M_2$ are isomorphic. 
\end{corollary}

        For $k = \mathbf{C}$, a transcendental proof of this was given by Deligne in \cite{Deligne80}.

\begin{remark}\label{* = otimes over kbar4} Let $k$ be as in \ref{* = otimes over kbar}, and let $\ell$ be a prime number invertible in $k$. For $\Lambda = \mathbf{Z}/\ell^{\nu}\mathbf{Z}$, with $\nu \ge 1$ fixed, projective system of isomorphisms $\beta_{\bullet} = (\beta_{\underline{n}})$ form a torsor under $(\Lambda[[\mu \times \mu]])^* = \varprojlim (\Lambda[\mu_{n_1} \times \mu_{n_2}])^*$. As $(\mathbf{Z}/\ell^{\nu +1}[[\mu \times \mu]])^* \to (\mathbf{Z}/\ell^{\nu}[[\mu \times \mu]])^*$ is surjective, there exists an isomorphism of projective systems of $\mu_{n_1} \times \mu_{n_2}$-equivariant $\mathbf{Z}_{\ell}$-sheaves
\begin{equation}
(\mathcal{B}_{\mathbf{Z}_{\ell},k}(\underline{n})) \iso (\pi_{n_1*}\mathbf{Z}_{\ell} \otimes \pi_{n_2*}\mathbf{Z}_{\ell})
\end{equation}
(resp. 
\begin{equation}
(\mathcal{B}^{\mathrm{loc}}_{\mathbf{Z}_{\ell},k}(\underline{n})) \iso (\pi^{\mathrm{loc}}_{n_1*}\mathbf{Z}_{\ell} \otimes \pi^{\mathrm{loc}}_{n_2*}\mathbf{Z}_{\ell})).
\end{equation}
This implies the analogue of \ref{* = otimes over kbar3} for $\mathbf{Z}_{\ell}$-sheaves, as well as $O_{\lambda}$-sheaves, $E_{\lambda}$-sheaves, and $\overline{\mathbf{Q}}_{\ell}$-sheaves (where $E_{\lambda}$ is a finite extension of $\mathbf{Q}_{\ell}$, and $O_{\lambda}$ its ring of integers).
\end{remark}

\subsection{}\label{obstruction} No longer assume that $k$ is algebraically closed.  Let $n_i$ be as in \ref{B}, with $\mu_r(\overline{k}) \subset k$. There is in general no isomorphism $\mathcal{B}_{\Lambda,\overline{k}}(\underline{n}) \iso (\pi_{n_1*}\Lambda \otimes \pi_{n_2*}\Lambda)_{\overline{k}}$ compatible with the action of $\mathrm{Gal}(\overline{k}/k)$, as \ref{J and Jacobi sum} shows. The obstruction is the class in $H^1(\SP \,k, (\Lambda[\mu_{\underline{n}}])^*)$ of the geometrically constant, locally free $\Lambda[\mu_{\underline{n}}]$-module of rank 1, $B_{1,2}(\underline{n},\Lambda) \otimes (T_{1,2}(\underline{n},\Lambda))^{-1}$, where $B_{1,2}(\underline{n},\Lambda) = \mathcal{B}_{\Lambda,k}(\underline{n})$, $T_{1,2}(\underline{n},\Lambda) = \pi_{n_1*}\Lambda \otimes \pi_{n_2*}\Lambda$.  Let $\ell$ be a prime number invertible in $k$, and let us now restrict to rings $\Lambda$ killed by a power of $\ell$. Then $B_{1,2}(\underline{n},\Lambda) \otimes (T_{1,2}(\underline{n},\Lambda))^{-1}$ comes by change of rings $\mathbf{Z}_{\ell} \to \Lambda$ and field extension $\mathbf{F}_p \to k$ from the geometrically constant, locally free $\mathbf{Z}_{\ell}[\mu_{\underline{n}}]$-module of rank 1 over $\mathbf{G}_{m,\mathbf{F}_p}$
\begin{equation}\label{obstruction1}
J_{1,2}(\underline{n}) :=  B_{1,2}(\underline{n},\mathbf{Z}_{\ell},\mathbf{F}_p) \otimes (T_{1,2}(\underline{n},\mathbf{Z}_{\ell}, \mathbf{F}_p))^{-1}
\end{equation}
For variable $\underline{n}$, the $J_{1,2}(\underline{n})$'s form a projective system denoted $J_{1,2}$. The next result is due to Deligne :
\begin{proposition}\label{Deligne associativity} Denote by $A_i$ ($i = 1, 2, 3$) a copy of the projective system $(\pi_{n*}\mathbf{Z}_{\ell})$ ($n$ running through the integers prime to $p$) on $\mathbf{G}_{m,\mathbf{F}_p}$. Then the associativity property of $*_1$ (\ref{associativity etc}) gives an isomorphism
\begin{equation}\label{Deligne associativity2}
J_{12,3} \otimes J_{1,2} \iso J_{1,23} \otimes J_{2,3}
\end{equation}
with the notation of (\ref{obstruction1}) and $J_{12,3} := ((A_1 \otimes A_2) *_1 A_3) \otimes ((A_1 \otimes A_2) \otimes A_3)^{-1}$, $J_{1,23} := (A_1 *_1 (A_2 \otimes A_3)) \otimes (A_1 \otimes (A_2 \otimes A_3))^{-1}$.
\end{proposition}
\begin{proof} The associativity isomorphism
\begin{equation*}
(A_1 *_1 A_2) *_1 A_3 \iso A_1 *_1 (A_2 *_1 A_3)
\end{equation*}
can be re-written
\begin{equation*}
((A_1 \otimes A_2) \otimes J_{1,2}) *_1 A_3 \iso A_1 *_1 ((A_2 \otimes A_3) \otimes J_{2,3}).
\end{equation*}
As $J_{1,2}$ and $J_{2,3}$ are geometrically constant, by (\ref{global unit4}) (and associativity and commutativity of $*_1$), the left hand side is isomorphic to
\begin{equation*}
((A_1 \otimes A_2) *_1 A_3) \otimes J_{1,2} \iso ((A_1 \otimes A_2) \otimes A_3) \otimes J_{12,3} \otimes J_{1,2},
\end{equation*}
and the right hand side to 
\begin{equation*}
((A_1 \otimes (A_2 \otimes A_3) \otimes J_{1,23})) *_1 J_{2,3} \iso A_1 \otimes (A_2 \otimes A_3) \otimes J_{1,23} \otimes J_{2,3}.
\end{equation*}
As $A_1 \otimes A_2 \otimes A_3$ is locally free of rank 1 over $\mathbf{Z}_{\ell}[\mu \times \mu \times \mu]$ (where $\mu := (\mu_{\bullet}$), (\ref{Deligne associativity2}) follows.
\end{proof}

\begin{remark}\label{Deligne associativity3} By \ref{J and Jacobi sum}, (\ref{Deligne associativity2}) implies the following identity between Jacobi sums :
\begin{equation}\label{Deligne associativity4}
J(\chi_1 \otimes \chi_2,\chi_3) J(\chi_1,\chi_2) = J(\chi_1, \chi_2 \otimes \chi_3)J(\chi_2,\chi_3),
\end{equation}
for non-trivial characters $\chi_i : \mathbf{F}_q^* \to E_{\lambda}^*$ ($i = 1, 2, 3)$, with the notations of \ref{J and Jacobi sum}.  ((\ref{Deligne associativity4}) follows from the classical identity between Gauss and Jacobi sums, see e. g. (\cite{SGA41/2}, Sommes trigonométriques, (4.15.2)).)
\end{remark}

\begin{remark}\label{unipotent} For $i = 1, 2$ let $M_i$ be a lisse $\overline{\mathbf{Q}}_{\ell}$-sheaf on $\mathbf{G}_{m,k}$, tamely ramified at $\{0\}$ and $\{\infty\}$, such that the image of the representation $\rho_i : \pi_1(\mathbf{G}_{m,k},\overline{1}) \to \mathrm{GL}(V_i)$, where $V_i = (M_i)_{\overline{1}}$, is contained in a unipotent subgroup of $\mathrm{GL}(V_i)$. Using (\ref{Deligne associativity2}) Deligne has shown (unpublished) that under these assumptions there is an isomorphism $M_1 *_1 M_2 \iso M_1 \otimes M_2$. 
\end{remark}

\begin{remark}\label{obstruction loc}Constructions and results from \ref{obstruction} to \ref{unipotent} imply local variants (on $A_h$), with $\pi_n$ replaced by $\pi^{\mathrm{loc}}_n$, $\mathcal{B}$ by $\mathcal{B}^{\mathrm{loc}}$, etc. 
\end{remark}

\subsection{}\label{variation} Let us now return to the situation considered in \ref{isolated4}, and assume $k$ algebraically closed. Let us write $(S,s)$ for $(A_h = A_{sh},\{0\})$, $(T,t)$ for $((A^2)_h,\{(0,0)\})$, $f : X \to T$ for $f_h : X_h \to (A^2)_h$, and denote by $a : T \to S$ the morphism defined by the sum map.  

Assume that $R^{n_i}\Psi_{f_i}(\Lambda)$ is tame at $x_i$ (so that, by (\ref{Thom-Sebastiani1}) and \ref{contracted product2 local}), $R\Psi_{af}(\Lambda)$ is tame at $x$). Then, by (\cite{Illusie03}, 3.5), if $\sigma$ is a topological generator of $I_t$, the variation morphism (induced by $\sigma - 1$)
\begin{equation}\label{variation1}
\mathrm{Var}(\sigma)_i : R^{n_i}\Phi_{f_i}(\Lambda)_{x_i} \to H^{n_i}_{\{x_i\}}((X_i)_s, R\Psi_{f_i}(\Lambda))
\end{equation}
is an isomorphism. Similarly, 
\begin{equation}\label{variation2}
\mathrm{Var}(\sigma) : R^{n+1}\Phi_{af}(\Lambda)_{x} \to H^{n+1}_{\{x\}}(X_s,R\Psi_{af}(\Lambda)),
\end{equation}
where $X_s := (af)^{-1}(s)$, is an isomorphism. Note that $\mathrm{Var}(\sigma)$ commutes with $I_t$, hence defines an isomorphism $\mathrm{Var}(\sigma)_i^*$ between the sheaves (on $S$) $R^{n_i}\Phi_{f_i}(\Lambda)_{x_i}$ and $H^{n_i}_{\{x_i\}}((X_i)_s, R\Psi_{f_i}(\Lambda))$ (resp. an isomorphism $\mathrm{Var}(\sigma)^*$ between $R^{n+1}\Phi_{af}(\Lambda)_{x}$ and $H^{n+1}_{\{x\}}(X_s, R\Psi_{af}(\Lambda))$). 

Recall that, by (\cite{Illusie03}, 3.7), we have a commutative diagram
\begin{equation}\label{variation6}
\xymatrix{R\Psi_{f_i}(\Lambda)[n_i] \ar[d]_{(\sigma-1)_i} \ar[r] & i_{x_i*}R^{n_i}\Phi_{f_i}(\Lambda)_{x_i} \ar[d]_{\mathrm{Var}(\sigma)_i} \\
R\Psi_{f_i}(\Lambda)[n_i] & i_{x_i*}H^{n_i}_{x_i}((X_i)_s,R\Psi_{f_i}(\Lambda)) \ar[l]},
\end{equation}
where the horizontal maps are adjunction maps, which are dual to each other (up to a Tate twist). Moreover, by (\cite{Illusie03}, 3.8), $i_{x_i*}R^{n_i}\Phi_{f_i}(\Lambda)_{x_i} = R\Phi_{f_i}(\Lambda)[n_i]$, the entries of (\ref{variation6}) are perverse sheaves, and $\mathrm{Var}(\sigma)_i$ is the canonical isomorphism from the co-image of the left vertical map to its image (in the category of perverse sheaves).  In particular, $(\sigma -1)_i : R\Psi_{f_i}(\Lambda)[n_i] \to R\Psi_{f_i}(\Lambda)[n_i]$ factors uniquely through $R\Phi_{f_i}(\Lambda)[n_i]$ :
\begin{equation}\label{variation6a}
\xymatrix{R\Psi_{f_i}(\Lambda)[n_i] \ar[d]_{(\sigma-1)_i} \ar[r] & R\Phi_{f_i}(\Lambda)[n_i] \ar[dl]^{\mathrm{Var}(\sigma)_i} \\
R\Psi_{f_i}(\Lambda)[n_i] & {}}
\end{equation}

The isomorphisms (\ref{Thom-Sebastiani1}) and (\ref{Thom-Sebastiani2}) define an isomorphism from the triangle
\begin{equation}
\xymatrix{R\Psi_{f_1}(\Lambda)[n_1] *^L R\Psi_{f_2}(\Lambda)[n_2] \ar[d]_{(\sigma-1)_1*^L (\sigma-1)_2} \ar[r] & R\Phi_{f_1}(\Lambda)[n_1] *^L R\Phi_{f_2}(\Lambda)[n_2] \ar[dl]^{\mathrm{Var}(\sigma)_1 *^L \mathrm{Var}(\sigma)_2} \\
R\Psi_{f_1}(\Lambda)[n_1] *^L R\Psi_{f_2}(\Lambda)[n_2] & {}}
\end{equation}
to the triangle
\begin{equation}
\xymatrix{R\Psi_{af}(\Lambda)[n+1] \ar[d]_{\sigma-1} \ar[r] & R\Phi_{af}(\Lambda)[n+1] \ar[dl]^{\mathrm{Var}(\sigma)} \\
R\Psi_{af}(\Lambda)[n+1] & {}}
\end{equation}

Using now the factorizations given by the lower triangles of (\ref{variation6}), we find that the isomorphisms (\ref{Thom-Sebastiani1}) and (\ref{Thom-Sebastiani2}) define an isomorphism from the triangle
\begin{equation}
\xymatrix{{} & R\Phi_{f_1}(\Lambda)[n_1] *^L R\Phi_{f_2}(\Lambda)[n_2] \ar[d]^{\mathrm{Var}(\sigma)_1 *^L \mathrm{Var}(\sigma)_2} \ar[dl]_{\mathrm{Var}(\sigma)_1 *^L \mathrm{Var}(\sigma)_2}\\
R\Psi_{f_1}(\Lambda)[n_1] *^L R\Psi_{f_2}(\Lambda)[n_2] & i_{x_1*}H^{n_i}_{x_1}((X_i)_s,R\Psi_{f_1}(\Lambda)) *^L i_{x_2*}H^{n_2}_{x_2}((X_2)_s,R\Psi_{f_2}(\Lambda)) \ar[l]}
\end{equation}
to the triangle
\begin{equation}
\xymatrix{{} & R\Phi_{f}(\Lambda)[n+1] \ar[d]^{\mathrm{Var}(\sigma)} \ar[dl]_{\mathrm{Var(\sigma)}} \\
R\Psi_{af}(\Lambda)[n+1] & i_{x*}H^{n+1}_{x}(X_s,R\Psi_{af}(\Lambda)) \ar[l]}.
\end{equation}

In particular, we have a commutative square of isomorphisms
\begin{equation}\label{variation4}
\xymatrix{R^{n_1}\Phi_{f_1}(\Lambda)_{x_1} *_1 R^{n_2}\Phi_{f_2}(\Lambda)_{x_2} \ar[r] \ar[d]_{\mathrm{Var}(\sigma)^*_1 *_1 \mathrm{Var}(\sigma)^*_2}& R^{n+1}\Phi_{af}(\Lambda)_x \ar[d]_{\mathrm{Var}(\sigma)^*} \\
H^{n_1}_{\{x_1\}}((X_1)_s, R\Psi_{f_1}(\Lambda)) *_1 H^{n_2}_{\{x_2\}}((X_2)_s, R\Psi_{f_1}(\Lambda)) \ar[r] &H^{n+1}_{\{x\}}(X_s, R\Psi_{af}(\Lambda))},
\end{equation}
where the upper horizontal map is the isomorphism (\ref{Thom-Sebastiani2}), and  the vertical maps are the isomorphisms defined by (\ref{variation1}) and (\ref{variation2}). In other words, with the above identifications, we have 
\begin{equation}\label{variation5}
\mathrm{Var}(\sigma)^* = \mathrm{Var}(\sigma)^*_1 *_1 \mathrm{Var}(\sigma)^*_2.
\end{equation}
\begin{remark}
(a) Presumably, the isomorphism
\begin{equation}\label{variation3}
H^{n+1}_{\{x\}}(X_s, R\Psi_{af}(\Lambda)) \iso H^{n_1}_{\{x_1\}}((X_1)_s, R\Psi_{f_1}(\Lambda)) *_1 H^{n_2}_{\{x_2\}}((X_2)_s, R\Psi_{f_1}(\Lambda)),
\end{equation}
in (\ref{variation4}) is the inverse of the isomorphism deduced by duality from the upper horizontal isomorphism, using the perfect pairings 
$$
R^{n_i}\Phi_{f_i}(\Lambda)_{x_i}  \otimes H^{n_i}_{x_i}((X_i)_s,R\Psi_{f_i}(\Lambda)) \to H^{2n_i}_{\{x_i\}}((X_i)_s,K_{(X_i)_s}[-2n_i]) = \Lambda
$$
(resp.
$$
R^{n+1}\Phi_{af}(\Lambda)_{x}  \otimes H^{n+1}_{x}((X_i)_s,R\Psi_{af}(\Lambda)) \to H^{2n+2}_{\{x\}}((X)_s,K_{(X)_s}[-2n-2]) = \Lambda),
$$
where $n = n_1 +n_2$, and $K$ denotes a dualizing complex ($= g^!\Lambda$, with $g$ the projection to $\SP \,k$). I haven't checked it.

(b) In view of \ref{* = otimes over kbar3}, we could replace $*_1$ by $\otimes$ in (\ref{variation5}). For $k = \mathbf{C}$, we recover (0.3).

\end{remark}

\medskip
\textit{Acknowledgements.} {\small I am grateful to Pierre Deligne, Ofer Gabber and Gérard Laumon for illuminating conversations, and Fabrice Orgogozo for a careful reading of a first draft of this paper and helpful comments. I thank T. Saito for helping me correct an argument in (\ref{psi-good-examples} (b)), suggesting \ref{locally acyclic ter}, statements and proofs of \ref{tameness at infinity}, \ref{dimtot}, and W. Sawin for a remark on (\ref{remark on Thom-Sebastiani} (c)) and the communication of references \cite{Rojas-Leon:2013aa} and \cite{Sawin14}. I thank Alexander Beilinson, Vladimir Drinfeld, Lei Fu, Takeshi Saito, and Weizhe Zheng for precious comments and criticism. I am especially grateful to Weizhe Zheng for sending me his note \cite{Zheng16}, and accepting to submit it as an appendix to this paper.  I thank Tung-Yang Lee for the realization of the figure illustrating the proof of 5.11. Finally, I wish to heartily thank the referee for a meticulous examination of the manuscript and the request of many corrections and improvements. Parts of this paper were written or revised while I was visiting the Department of Mathematics of the University of Chicago, the Morningside Center of Mathematics of the AMSS in Beijing, the Graduate School of Mathematical Sciences at the University of Tokyo, and the Department of Mathematics at Academia Sinica, NCTS, Taipei.  I warmly thank these institutions for their hospitality and support.}

\end{document}


\appendix
\section{K\"unneth formula for nearby cycles}

Weizhe Zheng\footnote{Morningside Center of Mathematics, Academy of
Mathematics and Systems Science, Chinese Academy of Sciences, Beijing
100190, China; University of the Chinese Academy of Sciences, Beijing
100049, China; email: \texttt{wzheng@math.ac.cn}. Partially supported by
National Natural Science Foundation of China Grant 11321101.}

\bigskip
The goal of this appendix is to give a short proof of Illusie's K\"unneth
formula for the nearby cycle functor $R\Psi$ over general bases
\cite[Theorem 2.3]{Illusie}, by imitating the proof of the K\"unneth formula
for the direct image functor $Rf_*$ over a field [SGA5 III (1.6.4)]. We also
discuss a couple of projection formulas related to $R\Psi$, which were
proved earlier by Takeshi Saito under slightly different assumptions
\cite[Section 3.1]{Saito}. This appendix owes very much to the work of Luc
Illusie and Takeshi Saito.

Let $\Lambda$ be a commutative ring.

\begin{definition}
We say that a morphism of schemes $f\colon X\to Y$ satisfies $(*)$ if there
exists an open covering $(U_\alpha)$ of $X$ and integers $(N_\alpha)$ such
that for every geometric point $x$ of $U_\alpha$, the functor $f_{(x)*}$ has
$\Lambda$-cohomological dimension $\le N_\alpha$. Here $f_{(x)}\colon
X_{(x)}\to Y_{(x)}$ denotes the strict localization of $f$ at $x$ and
$Y_{(x)}$ denotes the strict localization of $Y$ at the geometric point
$x\to X\xrightarrow{f} Y$. We say that $f$ satisfies $(*)$ universally if
every base change of $f$ satisfies $(*)$.
\end{definition}

It follows from a result of Gabber \cite[Proposition 3.1]{Or} that, if
$\Lambda$ is torsion and $f$ is locally of finite type, then $f$ satisfies
$(*)$ universally. See \cite[Corollary 1.15]{Illusie}.

Following \cite[Definition 1.5]{Illusie}, we say that a pair $(f,K)$, where
$f\colon X\to Y$ is a morphism of schemes and $K\in D(X,\Lambda)$, is
\emph{$\Psi$-good} if $R\Psi_f K$ commutes with base change by every
morphism of schemes $Y'\to Y$. Note that, by the compatibility of base
change maps with composition, if $(f,K)$ is $\Psi$-good, then $(f',K')$ is
$\Psi$-good for any base change $f'\colon X'\to Y'$ of $f$. Here $K'$
denotes the restriction of $K$ to $X'$.

\begin{remark}
Note that condition $(*)$ implies that the description of $R\Psi_f K$ by
sections \cite[(1.12.6)]{Illusie} extends to the unbounded derived category:
for all $K\in D(X,\Lambda)$, we have
\begin{equation}\label{e.RPsisec}
\sigma_x^* R\Psi_f K\simeq Rf_{(x)*}K_{(x)}.
\end{equation}
By this description, if either $K$ is in $D^+(X,\Lambda)$ or $f$ satisfies
$(*)$ universally, then we have the following equivalences:
\begin{enumerate}
\item $(f,K)$ is $\Psi$-good if and only if for every geometric point $x$
    of $X$, $Rf_{(x)*}K_{(x)}$, where $K_{(x)}$ denotes the restriction of
    $K$ to $X_{(x)}$, commutes with \emph{localized} base change: for any
    commutative square of strictly local schemes and local morphisms
\[\xymatrix{X'\ar[r]^{h}\ar[d]_{f'} & X_{(x)}\ar[d]^{f_{(x)}}\\
Y'\ar[r]^g & Y_{(y)}}
\]
that is \emph{essentially Cartesian} in the sense that the induced
morphism $X'\to X_{(x)}\times_{Y_{(y)}} Y'$ is a strict localization, the
base change map $g^*Rf_{(x)*} K_{(x)}\to Rf'_*h^* K_{(x)}$ is an
isomorphism.

\item $R\Psi_f K$ commutes with locally quasi-finite base change if and
    only if for every geometric point $x$ of $X$, $Rf_{(x)*}K_{(x)}$
    commutes with finite base change. In the case where $Y$ is
    quasi-compact and quasi-separated, this is further equivalent, by
    Zariski's Main Theorem, to the condition that $R\Psi_f K$ commutes
    with finite base change.
\end{enumerate}
\end{remark}

The following is a slight generalization of Illusie's K\"unneth formula for
$R\Psi$ \cite[Theorem 2.3]{Illusie}.

\begin{theorem}[K\"unneth formula for $R\Psi$]\label{t.Kunn}
Let $S$ be a scheme. Let $f_i\colon X_i\to Y_i$, $i=1,2$ be morphisms of
schemes over $S$ satisfying $(*)$ universally. Let $K_i\in
D^-(X_i,\Lambda)$, $i=1,2$. Assume either
\begin{enumerate}
\item $(f_1,K_1)$ and $(f_2,K_2)$ are $\Psi$-good; or
\item $Y_1=Y_2=S$, $(f_2,K_2)$ is $\Psi$-good, and $R\Psi_{f_1} K_1$
    commutes with locally quasi-finite base change.
\end{enumerate}
Then the K\"unneth map
\[c\colon R\Psi_{f_1} K_1\boxtimes_S^L
R\Psi_{f_2} K_2\to R\Psi_{f} K
\]
is an isomorphism. Here $f=f_1\times_S
f_2\colon X=X_1\times_S X_2\to Y_1\times_S Y_2=Y$, $K=K_1\boxtimes_S^L K_2$.
\end{theorem}

Note that case (1) follows from case (2) applied to the base change of $f_1$
and $f_2$ to $Y$.

\begin{proof}
Note that $c$ is an isomorphism if and only if, for each geometric point $x$
of $X$, the map $\sigma_x^* c$ is an isomorphism. By \eqref{e.RPsisec},
$\sigma_x^* c$ is the local K\"unneth map
\begin{equation}\label{e.localK}
Rf_{1(x)*} K_{1(x)}\boxtimes_{S_{(x)}}^L Rf_{2(x)*} K_{2(x)} \to Rf_{(x)*} (K_{1(x)}\boxtimes_{S_{(x)}}^L K_{2(x)}).
\end{equation}
Here $f_{i(x)}$ denotes the strict localization of $f_i$ at (the image of)
$x$, $K_{i(x)}$ denotes the restriction of $K_i$ to $X_{i(x)}$, and
$K_{1(x)}\boxtimes_{S_{(x)}}^L K_{2(x)}$ denotes $p_1^*K_{1(x)} \otimes^L
p_2^*K_{2(x)}$, where $p_i\colon X_{(x)}\to X_{i(x)}$ is the strict
localization of the projection at $x$, and similarly for the
$\boxtimes^L_{S_{(x)}}$ on the left. We imitate the proof of the K\"unneth
formula for $Rf_*$ over a field [SGA5 III (1.6.4)] as follows. Since
$(f_{1(x)},f_{2(x)})=(f_{1(x)},\id)(\id,f_{2(x)})$, we are reduced to
showing that the local K\"unneth maps for (a) $((\id, f_{2(x)}),
(K_{1(x)},K_{2(x)}))$ and (b) $((f_{1(x)},\id),
(K_{1(x)},Rf_{2(x)*}K_{2(x)}))$ are isomorphisms. By symmetry, since
$(\id,K_1)$ is $\Psi$-good, case (a) follows from case (1b). We are thus
reduced to showing the following generalization of (b): for all $K'_2\in
D^-(Y_{2(x)},\Lambda)$, the local K\"unneth map for $((f_{1(x)},\id),
(K_{1(x)},K'_2))$ is an isomorphism. Consider the commutative diagram with
essentially Cartesian squares
\[\xymatrix{X_{(x)}\ar[r]^{f_{(x)}}\ar[d]_{p_1} & Y_{(x)}\ar[r]^{q_2}\ar[d]^{q_1} & Y_{2(x)}\ar[d]\\
X_{1(x)}\ar[r]^{f_{1(x)}} & Y_{1(x)}\ar[r] & S_{(x)}.}
\]
The local K\"unneth map is the composition
\[q_1^*Rf_{1(x)*} K_{1(x)}\otimes^L q_2^*K'_2 \xrightarrow{a} Rf_{(x)*}p_1^*K_{1(x)}\otimes^L q_2^*K'_2 \xrightarrow{b} Rf_{(x)*} (p_1^* K_{1(x)}\otimes^L f_{(x)}^*q_2^* K'_2), \]
where the projection formula map $b$ is an isomorphism by
[SGA4$\frac{1}{2}$, Th.\ finitude, App., Remarque 1.2 a)] and the assumption
that $Rf_{(x)*}p_1^*K_{1(x)}$ commutes with finite base change. The base
change map $a$ is trivial in case (2), and an isomorphism by the assumption
that $(f_1,K_1)$ is $\Psi$-good in case (1).
\end{proof}

The following is a variant of \cite[Corollary 2.5]{Illusie}.

\begin{cor}
Let $f_i\colon X_i\to S$, $i=1,2$ be morphisms of schemes satisfying $(*)$
universally. Let $K_i\in D^b(X_i,\Lambda)$, $i=1,2$ such that
$K=K_1\boxtimes_S^L K_2\in D^b(X_1\times_S X_2,\Lambda)$. Assume that
$(f_1,K_1)$ is locally acyclic and $(f_2,K_2)$ is \emph{universally} locally
acyclic. Then $(f_1\times_S f_2, K)$ is locally acyclic.
\end{cor}

\begin{proof}
This follows from Theorem \ref{t.Kunn} (2) and the fact that a pair $(f,K)$
with $K\in D^+$ is locally acyclic if and only if $R\Phi_f K$ vanishes and
commutes with locally quasi-finite base change (\cite[Proposition
1.7]{Saitocc}, \cite[Example 1.7 (b)]{Illusie}).
\end{proof}

\begin{remark}\label{r.unbounded}\leavevmode
\begin{enumerate}
\item The result of [SGA4$\frac{1}{2}$, Th.\ finitude, App., Remarque 1.2
    a)] extends in fact to the case where $F\in D(S,A^{\op})$ is
    unbounded. Indeed, $Rf_*$ commutes with small direct sums by
    \cite[Proposition 14.3.4 (ii)]{KS}, and we reduce to the case $F\in
    D^-(S,A^\op)$ by \cite[Lemma 14.4.1]{KS}.

\item Similarly, for any morphism of schemes $f\colon X\to Y$ satisfying
    $(*)$, the functor $R\Psi_f\colon D(X,\Lambda)\to D(X\atimes_Y
    Y,\Lambda)$ commutes with small direct sums by \cite[Proposition
    14.3.4 (ii)]{KS}.
\end{enumerate}
\end{remark}

Recall that $R\Psi_{\id_Y}\simeq p_2^*$, where $p_2\colon Y\atimes_{Y} Y\to
Y$ is the second projection. Thus, in the case where $Y_1=Y_2=S$ and
$f_2=\id$, the K\"unneth formula reduces to a projection formula. The latter
extends to the case where $K_2$ is unbounded as follows.

\begin{prop}[Projection formula for $R\Psi_f$]
Let $f\colon X\to Y$ be a morphism of schemes satisfying $(*)$ universally.
Let $K\in D^-(X,\Lambda)$ be such that $R\Psi_f K$ commutes with locally
quasi-finite base change. Let $M\in D(Y,\Lambda)$. Then the projection
formula map
\[b\colon R\Psi_f K\otimes^L p_{2,f}^* M\to R\Psi_f(K\otimes^L f^*M)\]
is an isomorphism. Here $p_{2,f}\colon X\atimes_Y Y\to Y$ is the second
projection.
\end{prop}

\begin{proof}
This follows from Theorem \ref{t.Kunn} (2) and Remark \ref{r.unbounded} (2).
We can also prove it directly as follows. By \eqref{e.RPsisec}, for any
geometric point $x$ of $X$, $\sigma_x^* b$ is the projection formula map
$Rf_{(x)*}K_{(x)} \otimes^L M_{(x)}\to Rf_{(x)*}(K_{(x)}\otimes^L
f_{(x)}^*M_{(x)}) $, which is an isomorphism by Remark \ref{r.unbounded}
(1). Here as before $M_{(x)}$ denotes the restriction of $M$ to $Y_{(x)}$.
\end{proof}

Next we discuss another projection formula.

\begin{prop}[Projection formula for $R(\id\atimes g)_*$]
Let $X\xrightarrow{f} Y\xrightarrow{g} Z$ be morphisms of schemes, $g$
satisfying $(*)$. Then, for all $K\in D(X\atimes_Y Y,\Lambda)$ and $M\in
D(X,\Lambda)$, the projection formula map
\[b\colon R(\id_X\atimes g)_* K\otimes^L p_{1,gf}^* M\to R(\id_X\atimes g)_*(K\otimes^L p_{1,f}^*M)\]
is an isomorphism. Here $\id_X\atimes g\colon X\atimes_Y Y \to X\atimes_Z Z$
is the morphism induced by $g$, and $p_{1,f}\colon X\atimes_Y Y\to Y$ and
$p_{1,gf}\colon X\atimes_Z Z\to X$ are the first projections.
\end{prop}

Note that if $K=p_{2,f}^* K'$ for some $K'\in D(Y,\Lambda)$, then
\begin{multline*}
R(\id_X\atimes g)_* K\simeq R(\id_X\atimes g)_* (f\atimes
\id_Y)^*R\Psi_{\id_Y} K' \\
\simeq (f\atimes
\id_Z)^* R(\id_Y\atimes g)_* R\Psi_{\id_Y} K'\simeq (f\atimes \id_Z)^*R\Psi_g K'.
\end{multline*}
Here the second isomorphism is a trivial base change (see Lemma \ref{l.tbc}
below). The case $K'=\Lambda$ of the projection formula is used in
\cite[Proposition 1.16]{Illusie}.

\begin{proof}
The assumption $(*)$ implies that the description of $R(\id_X\atimes g)_*$
by sections \cite[Proposition 1.13]{Illusie} holds in the unbounded derived
category: for any $L\in D(X\atimes_Y Y,\Lambda)$ and for any geometric point
$x$ of $X$, we have
\begin{equation}\label{e.Rgsec}
\sigma_x^*R(\id_Y\atimes g)_* L\simeq Rg_{(x)*}\sigma_x^* L.
\end{equation}
Thus $\sigma^*_xb$ is the projection formula map $Rg_{(x)*} K_{x}\otimes^L
M_{x} \to Rg_{(x)*}(K_{x}\otimes^L M_{x})$, where $K_{x}$ is the restriction
of $K$ to $x\atimes_Y Y$, and $M_{x}$ is the stalk of $M$ at $x$. It
suffices to apply the following trivial projection formula.
\end{proof}

\begin{lemma}[Trivial projection formula]
Let $f\colon X\to Y$ be a morphism of topoi such that $R^qf_*\colon
\Sh(X,\Lambda)\to \Sh(Y,\Lambda)$ commutes with small direct sums for all
$q$ and vanishes for $q$ sufficiently large. Then, for all $K\in
D(X,\Lambda)$ and $M\in D(\Lambda)$, the projection formula map
\[Rf_*K\otimes^L M\to Rf_*(K\otimes^L M)\]
is an isomorphism.
\end{lemma}

The assumption that $R^q f_*$ commutes with small direct sums is satisfied if
$X$ is algebraic, $Y$ is locally coherent, and $f$ is coherent [SGA4 VI
Th\'eor\`eme 5.1].

\begin{proof}
Again by \cite[Proposition 14.3.4 (ii)]{KS}, $Rf_*$ commutes with small
direct sums. Thus by \cite[Lemma 14.4.1]{KS} we may assume that $K,M\in
D^-$. Taking a free resolution of $M$, we reduce to the trivial case where
$M=\Lambda$.
\end{proof}

\begin{lemma}[Trivial base change]\label{l.tbc}
Let $X\xrightarrow{f} X'\to Y\xrightarrow{g} Z$ be morphisms of schemes.
Consider the square of topoi
\[\xymatrix{X\atimes_Y Y\ar[r]^{f \atimes \id_Y} \ar[d]_{\id_X\atimes g}& X'\atimes_Y Y \ar[d]^{\id_{X'}\atimes
g}\\
X\atimes_Z Z\ar[r]^{f\atimes \id_Z} & X'\atimes_Z Z,}
\]
commutative up to canonical isomorphism. Let $K\in D(X'\atimes_Y Y)$. Assume
$K\in D^+$ or that $g$ satisfies $(*)$. Then the base change map
\[a\colon (f\atimes
\id_Z)^*R(\id_{X'}\atimes g)_*K\to
R(\id_{X}\atimes g)_*(f\atimes \id_Y)^* K
\]
is an isomorphism.
\end{lemma}

\begin{proof}
By \eqref{e.Rgsec}, for any geometric point $x$ of $X$, $\sigma_x^*a$ can be
identified with the identity on $Rg_{(x)*} K_{x}$.
\end{proof}